\documentclass{article}

\usepackage{Preamble}

\author{Samuel Lewis-Monkman}

\title{Weak Reducibility for Bridge Trisections}

\date{October 2026}

\begin{document}

\maketitle

\begin{abstract}
    A bridge splitting of a knot is weakly reducible if it has a pair of disjoint shadows, one from above and one from below, in the splitting surface. The process of weak reduction exchanges local minima and maxima in the splitting, untelescoping the splitting and giving a method for detecting essential surfaces in the knot complement. Bridge trisections generalize bridge splittings to the setting of surfaces in 4-manifolds. We define weak reducibility in this context, giving diagrammatic conditions for carrying out a weak reduction operation on bridge splittings. We produce a partial classification of 4-bridge trisections of surfaces in $S^4$.
\end{abstract}

\section{Introduction}

The theory of \textit{trisections}, first described in \cite{gk}, gives a method for decomposing 4-manifolds into 1-handlebodies, analogous to Heegaard splittings in dimension three. The idea was extended to surfaces embedded in $S^4$ in \cite{mz} as \textit{bridge trisections}, by which we may view a surface as a union of three trivial disk systems, just as a bridge splitting allows us to view a classical knot as a union of two trivial tangles. This theory does not merely extend the notion of splittings in three dimensions; as it turns out, due to a foundational result of Laudenbach and Poénaru \cite{lp}, the four-dimensional 1-handlebodies and disks trivially embedded therein are determined smoothly by their 3-dimensional boundaries, collectively called the trisection's ``spine". Thus, we can fully determine a smooth 4-manifold containing an embedded surface by a triple of 3-dimensional pieces fitting certain restrictions (described in Section 2.2).

Trisections give us a natural hierarchy of complexity: \textit{genus} and \textit{bridge number}. Just as with Heegaard splittings of 3-manifolds, a 4-manifold has a well-defined minimal genus for the central surface of a trisection. In the case of bridge trisections of surfaces embedded in a 4-manifold, we minimize the number of ``bridges" or arcs appearing in the 3-dimensional data of the splitting. A classification of 4-manifolds with trisection genus at most two is given by \cite{genus2} building on the genus one case from \cite{gk}. In \cite{mz}, we see that bridge trisections with up to three bridges are standard. Beyond this, the classification problems in both the manifold trisection and bridge trisection settings remain open and challenging.

Since the theory of trisections reduces a four-dimensional object to 3-dimensional information, it is natural to attempt to use tools from classical knot theory to study these objects. In particular, geometric ideas related to Heegaard splittings and bridge splittings show promise for generalization. In \cite{a&r}, Aranda and Zupan define \textit{weak reducibility}, a diagrammatic condition on Heegaard splittings of 3-manifolds, in the setting of trisected 4-manifolds, and use it to partially classify genus-three trisected 4-manifolds.

In this article, we define weak reducibility for embedded surfaces and exhibit a weak reduction operation that realizes in the diagram one of three elementary operations on the surface splitting:

\theoremstyle{plain}
\newtheorem*{thm:ElemOp}{Theorem \ref{thm:ElemOp}}
\begin{thm:ElemOp}
Let $\mathcal{T} = (\S; T_\a, T_\b, T_\g)$ be bridge trisection of a surface knot $K$ in $S^4$ such that $T_\a \cup \overline{T_\b}$ has a 1-bridge component $U_\e$ disjoint from some arc $\d \in T_\g$. Suppose that $\mathcal{T}_\d$, the result of forgetting $\d$, is a bridge trisection. Then $\mathcal{T}$ is weakly reducible with weak-reducing arcs $\d$, $\e_\a$, and  $\e_\b$, where $\e_\a \cup \overline{\e_\b} = U_\e$. Further, $\mathcal{T}_\d$ is related to $\mathcal{T}$ by a detubing, a deperturbation, or the removal of an $\R P^2$ summand.
\end{thm:ElemOp}

We prove this result by noting that weak reduction, which is a particular arc surgery on the tangles in the bridge trisections spine, admits a dual operation recovering the original splitting from the surgered one and then showing that if the surgered splitting consists of trivial tangles, then it satisfies conditions laid out in \cite{C&R} so that the dual surgery is a perturbation, tubing, or attachment of an $\R P^2$ summand. With this result in hand, we are able to classify all 4-bridge weakly-reducible splittings of surfaces in $S^4$. Such a surface is either unknotted (and possibly perturbed) or  a\textit{ribbon sphere}, part of a family of knotted $S^2$'s described in Section 2.2.

\theoremstyle{plain}
\newtheorem*{thm:4Br}{Theorem \ref{thm:4Br}}
\begin{thm:4Br}
A non-trivially knotted surface $K \subset S^4$ admits a 4-bridge weakly-reducible trisection if and only if it is a ribbon sphere with fusion number one.

\end{thm:4Br}

The classification comes by first showing (by way of a technical result regarding classical unlink splittings) that every weakly-reducible 4-bridge trisection weakly reduces to a 3-bridge trisection, applying Theorem \ref{thm:ElemOp}, and appealing to the existing classification of 3-bridge trisections. As an immediate result of Theorem \ref{thm:4Br}, we show that there exist 4-bridge trisections which are not weakly reducible.

\subsection{Organization}

Section 2 gives an overview of bridge splittings in dimension three and the classical notion of weak reducibility, followed by a parallel discussion of bridge trisections in dimension four and our new definition of weak reducibility in that setting. In Section 3, we prove some technical lemmas relating to the triples of tangles we deal with in bridge trisections. In particular, we show that every weakly-reducible 4-bridge splitting of a classical unlink satisfies certain conditions that become important in the following discussion. In Section 4, we prove our main results.

\subsection{Tools used}

The author made no use of generative AI in the preparation of this paper. All figures were made by the author, mostly using the excellent \href{http://www.mathcha.io}{Mathcha} editor.

\subsection{Acknowledgements}

The author would like to thank his advisor, Alex Zupan, for introducing him to the problem and offering endless advice and guidance throughout the project, as well as Román Aranda for his excellent mentorship, Jeffrey Meier for his willingness to answer a stream of questions, and Scott Taylor for his helpful and patient correspondence.

\section{Preliminaries}

\subsection{Splittings in Three Dimensions}

In classical (i.e., 3-dimensional) knot theory, we often study embeddings of 1-manifolds in 3-manifolds by decomposing them into simpler pieces. In particular, we can build knots and links by gluing together \emph{tangles}, collections of arcs embedded in 3-manifolds with boundary. In the following, tangles will always be properly and smoothly embedded, usually in $B^3$ or $S^2 \x I$. We will mostly care about \emph{trivial tangles}, which in $B^3$ are collections of arcs such that all arcs are simultaneously ambient isotopic (rel $\del$) into $\del B^3 = S^3$. In $S^2 \x I$, we allow arcs to be vertical (meaning isotopic to $\{x\} \x I$ for some $x \in S^2$), and require that all non-vertical arcs are simultaneously isotopic into a single boundary component, denoted $\del_+ $. Note that such an isotopy for a tangle $T$ in a manifold $M$ traces out a disk $B_a$ associated to each arc $a \in T$. When $B_a$ intersects $\del M$ transversely, we call it a \emph{bridge disk} for the arc $a$, and call its intersection $a^*$ with $\del M$ a \emph{shadow} for $B_a$.

Classical knots and links are often examined via bridge splittings, analogously to Heegaard splittings of 3-manifolds.

\begin{definition}[Bridge Splitting]
    Let $L$ be a 1-manifold properly embedded in $S^3$, $B^3$, or $S^2 \x I$. Let $\del_+ (S^2 \x I) = S^2 \x \{1\}$, $\del_- (S^2 \x I) = S^2 \x \{0\}$, and $\del_+ B^3 = \del B^3$. A \textit{bridge splitting} of $L$ is a decomposition of the form 
    \begin{itemize}
        \item $(S^3,L) = (B_1,T_1) \cup_{\del_+} (B_2,T_2)$,
        \item $(B^3,L) = (B_1, T_1) \cup_{\del_+} (C_2,T_2)$,
        \item $(S^2 \x I,L) = (C_1,T_1) \cup_{\del_+} (C_2,T_2)$, or
        \item $(S^2 \x I,L) = (C_1,T_1) \cup_{\del_-} (C_2,T_2)$,
    \end{itemize}
    where each $B_i$ is diffeomorphic to $B^3$ and each $C_i$ is diffeomorphic to $S^2 \x I$, and where each of $T_1$ and $T_2$ is a trivial tangle. We denote the \emph{splitting surface}, the shared boundary of the 3-manifolds, as $\S$.
\end{definition}

A full discussion of Heegaard splittings is omitted here, but we will note that one can obtain a Heegaard splitting for the double cover of $S^3$ branched over a link $L$ by taking, for each tangle $T_i$, the double cover of $B^3$ branched over $T_i$, which is a genus-$(n-1)$ 1-handlebody. We specify a gluing map between the boundaries by lifting a bridge disk for each arc to a compressing curve on $\S_{n-1}$, the boundary of the 1-handlebody.

We will often wish to examine a system of shadows for a bridge splitting.

\begin{definition}[Shadow Diagram]
    Let $\mathcal{B}_1$ and $\mathcal{B}_2$ be systems of bridge disks for $T_1$ and $T_2$, respectively, in a bridge splitting $\mathcal{H} = (\S;T_1,T_2)$ of a 1-manifold $L$. Record a shadow $a^*$ in $\S$ for each $B_a \in \mathcal{B}_1$, denoting by $T_1^*$ a collection of one shadow for each arc. Likewise, record a shadow $b^*$ for each $B_b \in \mathcal{B}_2$, denoting by $T_2^*$ the collection of one shadow for each arc. In the diagram, we mark the endpoints $\S \cap L$ of the shadows. We call the collection $(\S;T_1^*,T_2^*)$ of all these shadows a \emph{shadow diagram} for the splitting.
\end{definition}

Note that shadow diagrams are not, in general, unique. In particular, non-isotopic bridge disks will give rise to non-isotopic (rel $\del$) shadows. We can define a move on a system of shadows that allows us to go between diagrams for non-equivalent bridge disk systems.

\begin{definition}[Shadow Slide]
    Let $T^*$ be a system of shadows in a sphere $\S$ which is a component of the boundary of some 3-manifold $M$ containing a trivial tangle $T$. Let $a,b \in T^*$ be disjoint shadows, and let $\w$ be an arc in $\S$ connecting $a$ to $\del \nu(b)$, which is a curve in $\S$, such that $\w \cap (a \cup \nu(b))$ consists only of the endpoints of $\w$. By banding $a$ to $\del \nu(b)$ along $\w$ with surface framing, we obtain a new shadow $a'$ and a new system of shadows $T^{*\prime} = (T^* \bs \{a\}) \cup \{a'\}$ for the tangle $T$.
    We call this procedure a \emph{shadow slide} of $a$ over $b$.
\end{definition}

Shadow slides are quite useful. Indeed, by Proposition 5.2 of \cite{BraidedLinks}, we see that sliding a shadow in $T^*$ corresponds exactly to a different choice of bridge disk for the corresponding arc in $T$.

\

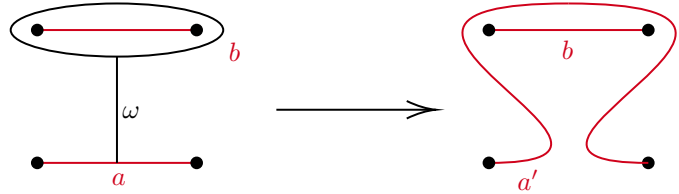
\begin{figure}[H]
    \centering

\tikzset{every picture/.style={line width=0.75pt}} 

\tikzset{every picture/.style={line width=0.75pt}} 

\begin{tikzpicture}[x=1pt,y=1pt,yscale=-1,xscale=1]
\path (120,100); 

\draw [color={rgb, 255:red, 208; green, 2; blue, 27 }  ,draw opacity=1 ]   (360,130) .. controls (440,130.17) and (280.33,69.83) .. (390,70) ;
\draw [color={rgb, 255:red, 208; green, 2; blue, 27 }  ,draw opacity=1 ]   (190,80) -- (250,80) ;
\draw [color={rgb, 255:red, 208; green, 2; blue, 27 }  ,draw opacity=1 ]   (190,130) -- (250,130) ;
\draw  [fill={rgb, 255:red, 0; green, 0; blue, 0 }  ,fill opacity=1 ] (188,130) .. controls (188,128.9) and (188.9,128) .. (190,128) .. controls (191.1,128) and (192,128.9) .. (192,130) .. controls (192,131.1) and (191.1,132) .. (190,132) .. controls (188.9,132) and (188,131.1) .. (188,130) -- cycle ;
\draw  [fill={rgb, 255:red, 0; green, 0; blue, 0 }  ,fill opacity=1 ] (248,130) .. controls (248,128.9) and (248.9,128) .. (250,128) .. controls (251.1,128) and (252,128.9) .. (252,130) .. controls (252,131.1) and (251.1,132) .. (250,132) .. controls (248.9,132) and (248,131.1) .. (248,130) -- cycle ;
\draw  [fill={rgb, 255:red, 0; green, 0; blue, 0 }  ,fill opacity=1 ] (188,80) .. controls (188,78.9) and (188.9,78) .. (190,78) .. controls (191.1,78) and (192,78.9) .. (192,80) .. controls (192,81.1) and (191.1,82) .. (190,82) .. controls (188.9,82) and (188,81.1) .. (188,80) -- cycle ;
\draw  [fill={rgb, 255:red, 0; green, 0; blue, 0 }  ,fill opacity=1 ] (248,80) .. controls (248,78.9) and (248.9,78) .. (250,78) .. controls (251.1,78) and (252,78.9) .. (252,80) .. controls (252,81.1) and (251.1,82) .. (250,82) .. controls (248.9,82) and (248,81.1) .. (248,80) -- cycle ;
\draw   (180,80) .. controls (180,74.48) and (197.91,70) .. (220,70) .. controls (242.09,70) and (260,74.48) .. (260,80) .. controls (260,85.52) and (242.09,90) .. (220,90) .. controls (197.91,90) and (180,85.52) .. (180,80) -- cycle ;
\draw    (220,90) -- (220,130) ;
\draw    (280,110) -- (338,110) ;
\draw [shift={(340,110)}, rotate = 180] [color={rgb, 255:red, 0; green, 0; blue, 0 }  ][line width=0.75]    (10.93,-3.29) .. controls (6.95,-1.4) and (3.31,-0.3) .. (0,0) .. controls (3.31,0.3) and (6.95,1.4) .. (10.93,3.29)   ;
\draw [color={rgb, 255:red, 208; green, 2; blue, 27 }  ,draw opacity=1 ]   (360,80) -- (420,80) ;
\draw  [fill={rgb, 255:red, 0; green, 0; blue, 0 }  ,fill opacity=1 ] (358,130) .. controls (358,128.9) and (358.9,128) .. (360,128) .. controls (361.1,128) and (362,128.9) .. (362,130) .. controls (362,131.1) and (361.1,132) .. (360,132) .. controls (358.9,132) and (358,131.1) .. (358,130) -- cycle ;
\draw  [fill={rgb, 255:red, 0; green, 0; blue, 0 }  ,fill opacity=1 ] (418,130) .. controls (418,128.9) and (418.9,128) .. (420,128) .. controls (421.1,128) and (422,128.9) .. (422,130) .. controls (422,131.1) and (421.1,132) .. (420,132) .. controls (418.9,132) and (418,131.1) .. (418,130) -- cycle ;
\draw  [fill={rgb, 255:red, 0; green, 0; blue, 0 }  ,fill opacity=1 ] (358,80) .. controls (358,78.9) and (358.9,78) .. (360,78) .. controls (361.1,78) and (362,78.9) .. (362,80) .. controls (362,81.1) and (361.1,82) .. (360,82) .. controls (358.9,82) and (358,81.1) .. (358,80) -- cycle ;
\draw  [fill={rgb, 255:red, 0; green, 0; blue, 0 }  ,fill opacity=1 ] (418,80) .. controls (418,78.9) and (418.9,78) .. (420,78) .. controls (421.1,78) and (422,78.9) .. (422,80) .. controls (422,81.1) and (421.1,82) .. (420,82) .. controls (418.9,82) and (418,81.1) .. (418,80) -- cycle ;
\draw [color={rgb, 255:red, 208; green, 2; blue, 27 }  ,draw opacity=1 ]   (390,70) .. controls (500.33,70.17) and (340.33,130.17) .. (420,130) ;

\draw (217,133) node [anchor=north west][inner sep=0.75pt]    {$\textcolor[rgb]{0.82,0.01,0.11}{a}$};
\draw (261,83.4) node [anchor=north west][inner sep=0.75pt]    {$\textcolor[rgb]{0.82,0.01,0.11}{b}$};
\draw (221,107.5) node [anchor=north west][inner sep=0.75pt]    {$\omega $};
\draw (369.67,131.4) node [anchor=north west][inner sep=0.75pt]    {$\textcolor[rgb]{0.82,0.01,0.11}{a'}$};
\draw (386.33,83) node [anchor=north west][inner sep=0.75pt]    {$\textcolor[rgb]{0.82,0.01,0.11}{b}$};

\end{tikzpicture}
    \caption{A shadow slide of $a$ over $b$}
    \label{fig:ShadowSlide}
\end{figure}

Bridge splittings are used to stratify classical knots and links by a natural notion of complexity: the \emph{bridge number} of a splitting is the number $n$ of arcs in either side of the splitting. Since each arc has its boundary in the splitting surface, the surface appears in the link exterior with $2n$ punctures. Often, a particular splitting of a link in $S^3$, for example, fails to minimize intersections with level surfaces. Doing so can help us detect geometry in the exterior of the link. \textit{Weak reducibility} is a a condition on shadows with an associated operation that allows us to ``untelescope" a splitting.

\begin{definition}[Weakly-Reducible Bridge Splitting]
    Let $\mathcal{H} = (\S;T_1,T_2)$ be a bridge splitting of a link $L$ in $S^3$. Suppose the splitting admits a shadow diagram with shadows $\d^* \in T_1^*$ and $\e^* \in T_2^*$ for arcs $\d \in T_1$ and $\e \in T_2$ such that $\d^* \cap \e^* = \0$. Then we say that the surface $\S$ is \emph{weakly-reducible} with \textit{weak-reducing arcs} $\d$ and $\e$. A surface that is not weakly-reducible and that does not admit a pair of shadows $a \in T_1^*$ and $b \in T_2^*$ with $a$ and $b$ isotopic is called \emph{strongly irreducible}.
\end{definition}

Weak reducibility for a bridge trisection in $S^3$ as defined here (i.e., with respect to the ambient genus $0$ Heegaard splitting of $S^3$) is exactly weak $T$-irreducibility in the sense of \cite{HS}, where $T$ is in our case a link $L$. The untelescoping operation discussed therein can be realized as follows in our special case.

\begin{figure}[H]
    \centering

\tikzset{every picture/.style={line width=0.75pt}} 

\begin{tikzpicture}[x=1pt,y=1pt,yscale=-1,xscale=1]

\draw [color={rgb, 255:red, 0; green, 0; blue, 0 }  ,draw opacity=1 ]   (280,110.15) .. controls (279.76,130.06) and (259.76,129.99) .. (260,110.07) ;
\draw [color={rgb, 255:red, 0; green, 0; blue, 0 }  ,draw opacity=1 ]   (300,110.22) .. controls (299.69,150.14) and (239.69,150.08) .. (240,110) ;
\draw [color={rgb, 255:red, 0; green, 0; blue, 0 }  ,draw opacity=1 ]   (60,80) .. controls (60.17,60.08) and (80.17,60.08) .. (80,80) ;
\draw [color={rgb, 255:red, 0; green, 0; blue, 0 }  ,draw opacity=1 ]   (120,80) .. controls (120.17,60.08) and (140.17,60.08) .. (140,80) ;
\draw [color={rgb, 255:red, 0; green, 0; blue, 0 }  ,draw opacity=1 ]   (100,80) .. controls (100.17,40.08) and (160.17,39.92) .. (160,80) ;
\draw [color={rgb, 255:red, 0; green, 0; blue, 0 }  ,draw opacity=1 ]   (160,80.37) .. controls (159.76,100.28) and (139.76,100.21) .. (140,80.29) ;
\draw [color={rgb, 255:red, 0; green, 0; blue, 0 }  ,draw opacity=1 ]   (100,80.15) .. controls (99.76,100.06) and (79.76,99.99) .. (80,80.07) ;
\draw [color={rgb, 255:red, 0; green, 0; blue, 0 }  ,draw opacity=1 ]   (120,80.22) .. controls (119.69,120.14) and (59.69,120.08) .. (60,80) ;
\draw [color={rgb, 255:red, 0; green, 0; blue, 0 }  ,draw opacity=1 ]   (50,80) -- (170,80) ;
\draw    (180,80) -- (218,80) ;
\draw [shift={(220,80)}, rotate = 180] [color={rgb, 255:red, 0; green, 0; blue, 0 }  ][line width=0.75]    (10.93,-3.29) .. controls (6.95,-1.4) and (3.31,-0.3) .. (0,0) .. controls (3.31,0.3) and (6.95,1.4) .. (10.93,3.29)   ;
\draw [color={rgb, 255:red, 0; green, 0; blue, 0 }  ,draw opacity=1 ]   (250,80) -- (310,80) ;
\draw [color={rgb, 255:red, 0; green, 0; blue, 0 }  ,draw opacity=1 ]   (250,50) -- (330,50) ;
\draw [color={rgb, 255:red, 0; green, 0; blue, 0 }  ,draw opacity=1 ]   (230,110) -- (310,110) ;
\draw [color={rgb, 255:red, 0; green, 0; blue, 0 }  ,draw opacity=1 ]   (280,50) .. controls (280.17,30.08) and (300.17,30.08) .. (300,50) ;
\draw [color={rgb, 255:red, 0; green, 0; blue, 0 }  ,draw opacity=1 ]   (260,50) .. controls (260.17,10.08) and (320.17,9.92) .. (320,50) ;
\draw [color={rgb, 255:red, 0; green, 0; blue, 0 }  ,draw opacity=1 ]   (280,110) .. controls (280.29,79.86) and (260.29,79.86) .. (260,50) ;
\draw [color={rgb, 255:red, 0; green, 0; blue, 0 }  ,draw opacity=1 ]   (300,110) .. controls (300.29,79.86) and (280.29,79.86) .. (280,50) ;
\draw [color={rgb, 255:red, 0; green, 0; blue, 0 }  ,draw opacity=1 ]   (240,110) .. controls (240.17,90.08) and (260.17,90.08) .. (260,110) ;
\draw [color={rgb, 255:red, 0; green, 0; blue, 0 }  ,draw opacity=1 ]   (320,50.07) .. controls (319.76,69.99) and (299.76,69.92) .. (300,50) ;

\draw (63,53) node [anchor=north west][inner sep=0.75pt]  [color={rgb, 255:red, 0; green, 0; blue, 0 }  ,opacity=1]  {$a_{\delta }$};
\draw (143,98) node [anchor=north west][inner sep=0.75pt]  [color={rgb, 255:red, 0; green, 0; blue, 0 }  ,opacity=1]  {$a_{\varepsilon }$};
\draw (36,72.4) node [anchor=north west][inner sep=0.75pt] {$\Sigma $};
\draw (184,62) node [anchor=north west][inner sep=0.75pt] [align=left] {{ W.R.}};
\draw (320,55) node [anchor=north west][inner sep=0.75pt]  [color={rgb, 255:red, 0; green, 0; blue, 0 }  ,opacity=1 ]  {$a_{\varepsilon }$};
\draw (230,90) node [anchor=north west][inner sep=0.75pt]  [color={rgb, 255:red, 0; green, 0; blue, 0 }  ,opacity=1 ]  {$a_{\delta }$};
\draw (312,75) node [anchor=north west][inner sep=0.75pt] {$\Sigma _{\delta \varepsilon }$};
\draw (233,44) node [anchor=north west][inner sep=0.75pt]  {$\Sigma _{\delta }$};
\draw (312,105) node [anchor=north west][inner sep=0.75pt]  {$\Sigma _{\varepsilon }$};

\end{tikzpicture}
    \caption{Weak reduction via a Morse function}
    \label{fig:WRex}
\end{figure}
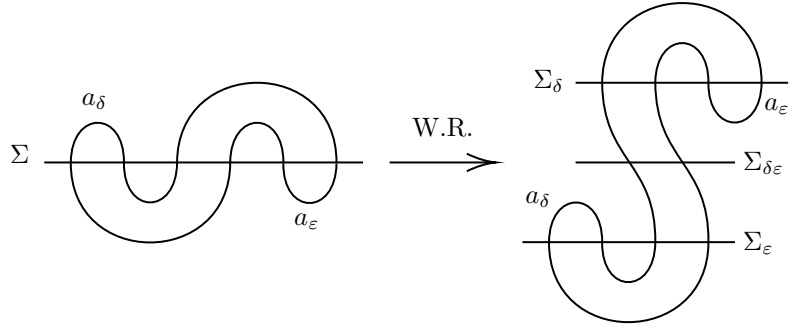

\begin{definition}
\label{def:WRClassical}
Let $\S$ be a surface with sets of shadows $T_1^*$ and $T_2^*$ for tangles $T_1$ and $T_2$ in a bridge splitting, with $T_1 \subset H_1$ and $T_2 \subset H_2$ for $H_1$ and $H_2$ each either $B^3$ or $S_2 \x I$. Suppose $\S$ is weakly-reducible along shadows $\d^* \in T_1^*$ and $\e^* \in T_2^*$ coming from bridge disks $B_\d$ and $B_\e$ for $\d$ and $\e$, respectively. We define $\S_\d$ as $(\S \bs \nu(\d^*)) \cup D_\d$, where $D_\d = (\del \nu(B_\d)) \cap H_1$, and $\S_\e$ likewise. By $\S_{\d\e}$ we denote $(\S_\d \bs \nu(\e)) \cup D_\e$.
We call this system of surgeries a \emph{weak reduction} of the splitting surface, as in Figure \ref{fig:WRex}.
    
\end{definition}

Quite useful in the study of smooth manifolds are Morse functions, smooth maps from a smooth manifold to $\R$ with isolated, non-degenerate critical levels. A thorough discussion of Morse functions and related theory can be found in \cite{Morse}. The relevant Morse function here is the natural height function, with respect to which a bridge splitting of a knot places all local minima below all local maxima, separated by the splitting surface.

On the level of this Morse function, weak reduction can be viewed as an isotopy of the splitting exchanging the heights of $a_\d$ and $a_\e$, in which case $\S_\d$, $\S_\e$, and $\S_{\d\e}$ are level surfaces separating critical points, which are no longer self-indexed. Notably, this process can be applied more than once when surfaces inside the new multiple splitting are themselves weakly-reducible.

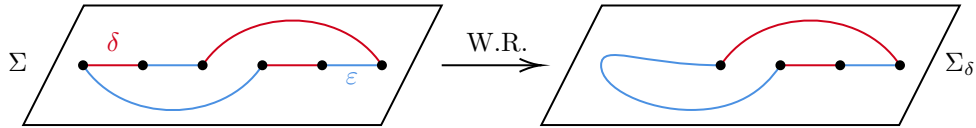
\begin{figure}[h]
    \centering

\tikzset{every picture/.style={line width=0.75pt}} 

\begin{tikzpicture}[x=0.75pt,y=0.75pt,yscale=-1,xscale=1]

\draw [color={rgb, 255:red, 74; green, 144; blue, 226 }  ,draw opacity=1 ]   (60,130) .. controls (80.13,160.07) and (130.13,159.93) .. (150,130) ;
\draw [color={rgb, 255:red, 208; green, 2; blue, 27 }  ,draw opacity=1 ]   (120,130) .. controls (140.15,100.08) and (190.17,100) .. (210,130) ;
\draw [color={rgb, 255:red, 74; green, 144; blue, 226 }  ,draw opacity=1 ]   (90,130) -- (120,130) ;
\draw [color={rgb, 255:red, 208; green, 2; blue, 27 }  ,draw opacity=1 ]   (150,130) -- (180,130) ;
\draw [color={rgb, 255:red, 74; green, 144; blue, 226 }  ,draw opacity=1 ]   (180,130) -- (210,130) ;
\draw [color={rgb, 255:red, 208; green, 2; blue, 27 }  ,draw opacity=1 ]   (60,130) -- (90,130) ;
\draw  [fill={rgb, 255:red, 0; green, 0; blue, 0 }  ,fill opacity=1 ] (58,130) .. controls (58,128.9) and (58.9,128) .. (60,128) .. controls (61.1,128) and (62,128.9) .. (62,130) .. controls (62,131.1) and (61.1,132) .. (60,132) .. controls (58.9,132) and (58,131.1) .. (58,130) -- cycle ;
\draw  [fill={rgb, 255:red, 0; green, 0; blue, 0 }  ,fill opacity=1 ] (88,130) .. controls (88,128.9) and (88.9,128) .. (90,128) .. controls (91.1,128) and (92,128.9) .. (92,130) .. controls (92,131.1) and (91.1,132) .. (90,132) .. controls (88.9,132) and (88,131.1) .. (88,130) -- cycle ;
\draw  [fill={rgb, 255:red, 0; green, 0; blue, 0 }  ,fill opacity=1 ] (118,130) .. controls (118,128.9) and (118.9,128) .. (120,128) .. controls (121.1,128) and (122,128.9) .. (122,130) .. controls (122,131.1) and (121.1,132) .. (120,132) .. controls (118.9,132) and (118,131.1) .. (118,130) -- cycle ;
\draw  [fill={rgb, 255:red, 0; green, 0; blue, 0 }  ,fill opacity=1 ] (148,130) .. controls (148,128.9) and (148.9,128) .. (150,128) .. controls (151.1,128) and (152,128.9) .. (152,130) .. controls (152,131.1) and (151.1,132) .. (150,132) .. controls (148.9,132) and (148,131.1) .. (148,130) -- cycle ;
\draw  [fill={rgb, 255:red, 0; green, 0; blue, 0 }  ,fill opacity=1 ] (178,130) .. controls (178,128.9) and (178.9,128) .. (180,128) .. controls (181.1,128) and (182,128.9) .. (182,130) .. controls (182,131.1) and (181.1,132) .. (180,132) .. controls (178.9,132) and (178,131.1) .. (178,130) -- cycle ;
\draw  [fill={rgb, 255:red, 0; green, 0; blue, 0 }  ,fill opacity=1 ] (208,130) .. controls (208,128.9) and (208.9,128) .. (210,128) .. controls (211.1,128) and (212,128.9) .. (212,130) .. controls (212,131.1) and (211.1,132) .. (210,132) .. controls (208.9,132) and (208,131.1) .. (208,130) -- cycle ;
\draw [color={rgb, 255:red, 74; green, 144; blue, 226 }  ,draw opacity=1 ]   (380,130) .. controls (340.11,129.98) and (320.11,118.43) .. (320,130) .. controls (319.89,141.57) and (342.78,152.43) .. (365.22,152.43) .. controls (387.67,152.43) and (401.29,143.15) .. (410,130) ;
\draw [color={rgb, 255:red, 208; green, 2; blue, 27 }  ,draw opacity=1 ]   (380,130) .. controls (400.15,100.08) and (450.17,100) .. (470,130) ;
\draw [color={rgb, 255:red, 208; green, 2; blue, 27 }  ,draw opacity=1 ]   (410,130) -- (440,130) ;
\draw [color={rgb, 255:red, 74; green, 144; blue, 226 }  ,draw opacity=1 ]   (440,130) -- (470,130) ;
\draw  [fill={rgb, 255:red, 0; green, 0; blue, 0 }  ,fill opacity=1 ] (378,130) .. controls (378,128.9) and (378.9,128) .. (380,128) .. controls (381.1,128) and (382,128.9) .. (382,130) .. controls (382,131.1) and (381.1,132) .. (380,132) .. controls (378.9,132) and (378,131.1) .. (378,130) -- cycle ;
\draw  [fill={rgb, 255:red, 0; green, 0; blue, 0 }  ,fill opacity=1 ] (408,130) .. controls (408,128.9) and (408.9,128) .. (410,128) .. controls (411.1,128) and (412,128.9) .. (412,130) .. controls (412,131.1) and (411.1,132) .. (410,132) .. controls (408.9,132) and (408,131.1) .. (408,130) -- cycle ;
\draw  [fill={rgb, 255:red, 0; green, 0; blue, 0 }  ,fill opacity=1 ] (438,130) .. controls (438,128.9) and (438.9,128) .. (440,128) .. controls (441.1,128) and (442,128.9) .. (442,130) .. controls (442,131.1) and (441.1,132) .. (440,132) .. controls (438.9,132) and (438,131.1) .. (438,130) -- cycle ;
\draw  [fill={rgb, 255:red, 0; green, 0; blue, 0 }  ,fill opacity=1 ] (468,130) .. controls (468,128.9) and (468.9,128) .. (470,128) .. controls (471.1,128) and (472,128.9) .. (472,130) .. controls (472,131.1) and (471.1,132) .. (470,132) .. controls (468.9,132) and (468,131.1) .. (468,130) -- cycle ;
\draw   (320.4,100) -- (500,100) -- (469.6,160) -- (290,160) -- cycle ;
\draw   (60.4,100) -- (240,100) -- (209.6,160) -- (30,160) -- cycle ;
\draw    (240,130) -- (288,130) ;
\draw [shift={(290,130)}, rotate = 180] [color={rgb, 255:red, 0; green, 0; blue, 0 }  ][line width=0.75]    (10.93,-3.29) .. controls (6.95,-1.4) and (3.31,-0.3) .. (0,0) .. controls (3.31,0.3) and (6.95,1.4) .. (10.93,3.29)   ;

\draw (70.4,112.2) node [anchor=north west][inner sep=0.75pt]    {$\textcolor[rgb]{0.82,0.01,0.11}{\delta }$};
\draw (190.4,132) node [anchor=north west][inner sep=0.75pt]    {$\textcolor[rgb]{0.29,0.56,0.89}{\varepsilon }$};
\draw (251,112) node [anchor=north west][inner sep=0.75pt]   [align=left] {W.R.};
\draw (21,122.4) node [anchor=north west][inner sep=0.75pt]    {$\Sigma $};
\draw (491,122.4) node [anchor=north west][inner sep=0.75pt]    {$\Sigma _{\delta }$};

\end{tikzpicture}
    \caption{A weak reduction modeled with shadows}
    \label{fig:WRExShadows}
\end{figure}

To realize this operation using shadows, we will need to understand how to resolve a system of arcs (for example, shadows) along a separate arc. For a system of arcs $A$ and an arc $b$ which intersects $A$ only at endpoints it shares with two distinct arcs $a,c \in A$, the \textit{resolution} of $A$ along $b$ is $A_b = (A \bs \{a,c\}) \cup \{a \cup b \cup c\}$. If a system of arcs $A$ and an arc $b \not\in A$ meet the criteria for resolution, we say that $b$ is \textit{primitive} with respect to $A$. In fact, in the setting of shadows, it is enough to know the shadow $b$ along which we resolve intersects one shadow $a$ only at a single endpoint, as we may slide any other shadows from the same tangle as $a$ which intersect $b$ in its interior over $a$, thus removing all interior intersections. In this case (if, say, $\d^* \in T_\a^*$ is primitive with respect to $T_\b^*$), we obtain a shadow diagram $T_{1\d}^{*}$ by deleting the arc $\d^*$ from $T_1^*$ and $T_{2\d}^*$ by resolving $T_2^*$ along $\d^*$, as in Figure \ref{fig:WRExShadows}. We define $T_{1\e}^*$ and $T_{2\e}^*$ similarly.

\

In order to understand untelescoped splittings, we consider tangles in a slightly more general setting.



\begin{definition}
     A \emph{multiple bridge splitting} of a 1-manifold $L$ in a 3-manifold $M$ diffeomorphic to $S^3$, $B^3$, or $S^2 \x I$ is a decomposition $(M, L) = (H_1, T_1) \cup (H_2, T_2) \cup \cdots \cup (H_\l, T_\l)$, where each $T_i \subset H_i$ is a trivial tangle, and where boundaries are glued with matching sign. We say that a surface $\S_+ = \del_+ H_i = \del_+ H_{i+1}$ in the splitting is \textit{thick} and $\S_- = \del_- H_i = \del_- H_{i+1}$ is \textit{thin}.
\end{definition}

Note that since we weakly reduce a particular splitting surface, each thick level in a multiple bridge splitting may admit a weak reduction. We record multiply-weakly-reduced surfaces by a list of arcs along which they have been weakly reduced. For example, the surface $\S_{\d\d'}$ has been weakly reduced first along an arc $\d$ and then along an arc $\d'$ in $\S_\d$.

\begin{rmk}
    We here restrict ourselves to multiple splittings where each level surface is a sphere. In the more general case, we can allow \textit{compression-bodies} (a generalization of 1-handlebodies) beyond the trivial $S^2 \x I$.
\end{rmk}

\subsection{Splittings in Four Dimensions}

In four dimensions, surfaces, rather than curves, may come with interesting embeddings. Surface knots and links have been studied via numerous related constructions; our discussion here is limited to \emph{bridge trisections} of surfaces, defined below.

The family of examples of surface knots most relevant to our discussion is that of \emph{ribbon surfaces}, first described in \cite{yaj}.

\begin{definition}[Ribbon Surface]
    A surface $S \subset S^4$ is \emph{ribbon} if it can be obtained by attaching finitely many \textit{tubes} to a collection $\mathcal{S}$ of unlinked spheres, each bounding a $B^3$, where tube attachment is a surgery replacing the neighborhood of a pair of points $a,b$ on one or two spheres in $\mathcal{S}$ with $\del (\a \x D^2)$ for an arc $\a$ disjoint in its interior from each sphere of $\mathcal{S}$ and with endpoints at $a$ and $b$. We call this decomposition a \textit{ball-tube presentation} of the surface.
\end{definition}

One construction that provides examples of ribbon spheres is that of \emph{spinning} a classical knot. We remove a small open ball in $S^3$ containing an unknotted arc of some classical knot $K$, fix the boundary, and sweep out a circle with the interior. This construction was originally described by Artin \cite{Artin} and expanded by Zeeman \cite{Zeeman}, adding a possible ``twisting" to the spin. In \cite{Fox}, Fox describes an additional ``rolling" deformation as part of spinning a knot.

\begin{figure}[h]
    \centering

\tikzset{every picture/.style={line width=0.75pt}} 
\begin{tikzpicture}[x=0.75pt,y=0.75pt,yscale=-1,xscale=1]

\draw  [draw opacity=0] (139.35,161.21) .. controls (128.52,184.65) and (110.45,200) .. (90,200) .. controls (56.86,200) and (30,159.71) .. (30,110) .. controls (30,60.29) and (56.86,20) .. (90,20) .. controls (109.37,20) and (126.59,33.77) .. (137.56,55.13) -- (90,110) -- cycle ; \draw   (139.35,161.21) .. controls (128.52,184.65) and (110.45,200) .. (90,200) .. controls (56.86,200) and (30,159.71) .. (30,110) .. controls (30,60.29) and (56.86,20) .. (90,20) .. controls (109.37,20) and (126.59,33.77) .. (137.56,55.13) ;  
\draw  [draw opacity=0] (180.55,59.01) .. controls (191.37,35.45) and (209.48,20) .. (230,20) .. controls (263.14,20) and (290,60.29) .. (290,110) .. controls (290,159.71) and (263.14,200) .. (230,200) .. controls (210.6,200) and (193.35,186.19) .. (182.39,164.77) -- (230,110) -- cycle ; \draw   (180.55,59.01) .. controls (191.37,35.45) and (209.48,20) .. (230,20) .. controls (263.14,20) and (290,60.29) .. (290,110) .. controls (290,159.71) and (263.14,200) .. (230,200) .. controls (210.6,200) and (193.35,186.19) .. (182.39,164.77) ;  
\draw    (105,75) .. controls (143.62,38.15) and (173.64,80) .. (210,50) ;
\draw    (110,80) .. controls (146.36,50) and (173.64,90) .. (210,60) ;
\draw    (215,120) .. controls (176.38,82.92) and (146.36,125) .. (110,95) ;
\draw    (210,125) .. controls (173.64,95) and (146.36,135) .. (110,105) ;
\draw    (110,160) .. controls (146.36,130) and (173.64,170) .. (210,140) ;
\draw    (110,170) .. controls (146.36,140) and (173.64,180) .. (210,150) ;
\draw  [draw opacity=0] (114.55,79.17) .. controls (113.76,82.61) and (112.02,85) .. (110,85) .. controls (107.24,85) and (105,80.52) .. (105,75) .. controls (105,69.48) and (107.24,65) .. (110,65) .. controls (110.79,65) and (111.53,65.36) .. (112.2,66.01) -- (110,75) -- cycle ; \draw   (114.55,79.17) .. controls (113.76,82.61) and (112.02,85) .. (110,85) .. controls (107.24,85) and (105,80.52) .. (105,75) .. controls (105,69.48) and (107.24,65) .. (110,65) .. controls (110.79,65) and (111.53,65.36) .. (112.2,66.01) ;  
\draw    (110,105) .. controls (106.87,101.47) and (107.8,97.33) .. (110,95) ;
\draw  [dash pattern={on 0.84pt off 2.51pt}]  (157.98,70) .. controls (161.12,66.47) and (160.18,62.33) .. (157.98,60) ;
\draw  [dash pattern={on 0.84pt off 2.51pt}]  (110,105) .. controls (113.13,101.47) and (112.2,97.33) .. (110,95) ;
\draw    (157.98,70) .. controls (154.85,66.47) and (155.78,62.33) .. (157.98,60) ;
\draw    (105,100) .. controls (103.16,97.93) and (102.72,95.65) .. (103.12,93.65) .. controls (103.4,92.24) and (104.09,90.97) .. (105,90) ;
\draw    (105,100) .. controls (108.13,96.47) and (107.2,92.33) .. (105,90) ;
\draw  [dash pattern={on 4.5pt off 4.5pt}]  (105,100) .. controls (91.73,95.27) and (89.67,85.67) .. (105,75) ;
\draw  [dash pattern={on 4.5pt off 4.5pt}]  (105,90) .. controls (100,90.73) and (101.8,85.53) .. (110,80) ;
\draw  [dash pattern={on 0.84pt off 2.51pt}]  (160,115) .. controls (163.13,111.47) and (162.2,107.33) .. (160,105) ;
\draw    (160,115) .. controls (156.87,111.47) and (157.8,107.33) .. (160,105) ;
\draw  [draw opacity=0] (207.68,111.13) .. controls (208.38,110.41) and (209.16,110) .. (210,110) .. controls (212.76,110) and (215,114.48) .. (215,120) .. controls (215,125.52) and (212.76,130) .. (210,130) .. controls (207.96,130) and (206.21,127.56) .. (205.43,124.06) -- (210,120) -- cycle ; \draw   (207.68,111.13) .. controls (208.38,110.41) and (209.16,110) .. (210,110) .. controls (212.76,110) and (215,114.48) .. (215,120) .. controls (215,125.52) and (212.76,130) .. (210,130) .. controls (207.96,130) and (206.21,127.56) .. (205.43,124.06) ;  
\draw    (210,150) .. controls (213.13,146.47) and (212.2,142.33) .. (210,140) ;
\draw  [dash pattern={on 0.84pt off 2.51pt}]  (210,150) .. controls (206.87,146.47) and (207.8,142.33) .. (210,140) ;
\draw    (215,145) .. controls (216.84,142.93) and (217.28,140.65) .. (216.88,138.65) .. controls (216.6,137.24) and (215.91,135.97) .. (215,135) ;
\draw    (215,145) .. controls (211.87,141.47) and (212.8,137.33) .. (215,135) ;
\draw  [dash pattern={on 4.5pt off 4.5pt}]  (215,145) .. controls (228.27,140.27) and (230.33,130.67) .. (215,120) ;
\draw  [dash pattern={on 4.5pt off 4.5pt}]  (215,135) .. controls (220,135.73) and (218.2,130.53) .. (210,125) ;
\draw    (210,60) .. controls (213,57.54) and (215,55.38) .. (220,55) ;
\draw    (210,50) .. controls (213,47.54) and (214.85,44.46) .. (215,40) ;
\draw    (110,160) .. controls (107,162.46) and (105,164.62) .. (100,165) ;
\draw    (110,170) .. controls (107,172.46) and (105.15,175.54) .. (105,180) ;
\draw  [draw opacity=0] (149.63,120.03) .. controls (148.99,128.76) and (147.51,137.11) .. (145.32,144.9) -- (90,110) -- cycle ; \draw   (149.63,120.03) .. controls (148.99,128.76) and (147.51,137.11) .. (145.32,144.9) ;  
\draw  [draw opacity=0] (144.43,72.07) .. controls (147.19,80.97) and (149.01,90.66) .. (149.69,100.84) -- (90,110) -- cycle ; \draw   (144.43,72.07) .. controls (147.19,80.97) and (149.01,90.66) .. (149.69,100.84) ;  
\draw  [draw opacity=0] (170.16,116.71) .. controls (170.72,127.99) and (172.66,138.69) .. (175.73,148.43) -- (230,110) -- cycle ; \draw   (170.16,116.71) .. controls (170.72,127.99) and (172.66,138.69) .. (175.73,148.43) ;  
\draw  [draw opacity=0] (174.82,74.61) .. controls (172.49,82.77) and (170.94,91.54) .. (170.32,100.72) -- (230,110) -- cycle ; \draw   (174.82,74.61) .. controls (172.49,82.77) and (170.94,91.54) .. (170.32,100.72) ;  
\draw  [dash pattern={on 0.84pt off 2.51pt}]  (160,160) .. controls (163.13,156.47) and (162.2,152.33) .. (160,150) ;
\draw    (160,160) .. controls (156.87,156.47) and (157.8,152.33) .. (160,150) ;

\end{tikzpicture}
    \caption{A ball-tube presentation of the spun trefoil}
    \label{fig:BTSpunTref}
\end{figure}
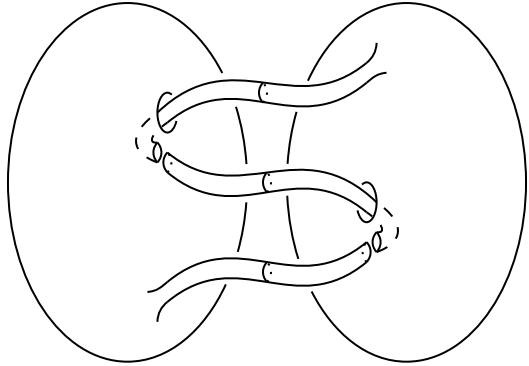

From Theorem 1 of \cite{yaj} and recognizing that spinning a classical knot produces a sphere whose projection into the equatorial $S^3$ in $S^4$ has as its singularity set a simple closed curve for each crossing in a chosen diagram of the knot, we can see that 0-twist, 0-roll spun knots are always ribbon. Figure \ref{fig:BTSpunTref} is a ball-tube presentation of the (0-twist, 0-roll) spun trefoil.

Knotted surfaces can be decomposed via bridge trisections, a construction due to Meier and Zupan in \cite{mz}, analogous to bridge splittings of classical knots. Bridge trisections allow us to study surface knots via triples of 3-dimensional data, just as trisections of 4-manifolds (detailed in \cite{gk}) do for those objects.

Bridge trisections make use of several ideas analogous to those in three dimensions. Instead of decomposing a knotted $S^1$ into arcs, we view a surface as a union of disks. A disk $D$ in a 4-dimensional handlebody $H$ is \textit{trivial} if it is isotopic to a disk in $\del H$. A system of disks is trivial if such an isotopy for each disk can all be realized via a single ambient isotopy, just as in the classical case.

In all discussions hereafter, we work in the smooth category.

\begin{definition}
    A \emph{bridge trisection} of a surface $S$ in $S^4$ with respect to the unique genus $0$ trisection of the ambient space is a decomposition $(S^4, S) = (X_1,D_1) \cup (X_2, D_2) \cup (X_3, D_3)$ where:
    \begin{itemize}
        \item $X_i \cong B^4$,
        \item $D_i = \sqcup_{c_i} D^2$ is the union of all disks in a trivial disk system $\mathcal{D}_i$,
        \item $X_i \cap X_{i+1} \cong B^3$,
        \item $D_i \cap D_{i+1}$ is a trivial tangle with $b$ bridges,
        \item $X_1 \cap X_2 \cap X_3 = \S \cong S^2$,
        \item $D_1 \cap D_2 \cap D_3$ is $2b$ marked points in $\S$,
    \end{itemize}
taking indices mod $3$. We call the tuple $\mathcal{T} := (\S;\mathcal{D}_1,\mathcal{D}_2, \mathcal{D}_3)$ a $(b;c_1,c_2,c_3)$-bridge trisection, or sometimes a $b$-bridge trisection of $S$.
\end{definition}

A few facts fall out of this definition immediately. We will state the most relevant here, but \cite{gk} contains the details foundational to the study of the trisected 4-manifolds, and \cite{mz} discusses bridge trisections of surfaces in $S^4$ in depth.

Firstly, notice that since each of the pairwise intersections $(H_\a,T_\a) = (X_1, D_1) \cap (X_2,D_2)$, $(H_\b,T_\b) = (X_2, D_2) \cap (X_3,D_3)$, and $(H_\g,T_\g) = (X_3, D_3) \cap (X_3,D_3)$ is a trivial tangle inside a 3-ball, their pairwise unions form $b$-bridge splittings of links in $S^3$. Further, since each such pair caps off to form a trivial disk system, it must be an unlink. This is helpful because, due to a result of Laudenbach and Po\'enaru \cite{lp}, the tuple $(\S;T_\a,T_\b,T_\g)$ of this 3-dimensional information (called the \emph{spine} of the trisection) determines the bridge trisection uniquely (as described in \cite{mz}), and we often denote the splitting using this tuple. As such, we may study bridge trisections, and therefore surface knots, by triples of trivial tangles that pairwise form unlinks.

Because we can describe a bridge trisection in this way, we may study a bridge trisection using a shadow diagram for the spine. A shadow diagram for a bridge-trisected surface consists of three systems of shadows, one for each tangle $T_\a$, $T_\b$, and $T_\g$ in the surface $\S$. A shadow diagram is sufficient to recover a bridge trisection, but not unique. As in the 3-dimensional case, we can easily produce non-isotopic shadow diagrams by sliding shadows from any particular tangle over one another.

We consider surfaces in $S^4$ up to isotopy, but the same surface may yield bridge trisections that are not obviously equivalent. In \cite{mz}, Meier and Zupan describe a set of moves on tangles necessary and sufficient to determine when two bridge trisections describe the same surface. The moves relevant to our discussion are perturbation (``stabilization" in language of \cite{mz}), which increases the bridge number of the trisection, and deperturbation (``destabilization" in language of \cite{mz}), which decreases the bridge number of the trisection. We say that a bridge trisection is \emph{perturbed} if it admits a deperturbation.

As in the case of a classical bridge splitting, we may build a trisection of the double cover of $S^4$ branched over $K$ by taking, for each tangle in the spine, the double cover of $B^4$ branched over the tangle and lifting a bridge disk for each arc to a compressing disk for the resulting 1-handlebody. The set of three handlebodies equipped with gluing maps that we obtain is called a ``Heegaard triple" in general, and in the case that each pairwise union of handlebodies is $\#_{k_i} S^1 \x S^2$ for some $k_i$, it is a trisection.

Equipped with diagrams of bridge-trisected surfaces consisting only of 3-dimensional information, we can define weak reduction for bridge trisections analogously to the classical notion.

\begin{definition}[Weakly-Reducible Bridge Trisection]
\label{def:WR4D}

Let $ \mathcal{T} = (\S, T_\a, T_\b, T_\g) $ be a bridge trisection of a surface link $ L $ in $ S^4 $, where each of $ T_\a $, $ T_\b $, and $ T_\g $ is a trivial tangle in a 3-ball $H_\a$, $H_\b$, and $H_\g$, respectively, and $ \S = \del H_\a = \del H_\b = \del H_\g $ is the splitting surface. 

Consider a shadow diagram $ (T_\a^*, T_\b^*, T_\g^*) $ for $ \mathcal{T} $. If this contains a disjoint pair of shadows $ \e^* \subset T_\a^* \cap T_\b^* $ and $ \d^* \subset T_\g^* $, then $ \mathcal{T} $ is \emph{weakly-reducible}. We call the shadows $\e^*$ and $\d^*$ \textit{weak-reducing shadows} and the arcs $\e_\a \in T_\a$, $\e_\b \in T_\b$, and $\d \in T_\g$ for which $\e^*$ and $\d^*$ are shadows the \textit{weak-reducing arcs} for the splitting.
    
\end{definition}

As in the 3-dimensional case, we perform weak reduction  by ``forgetting'' or surgering along the isolated reducing arc $\d$.

\begin{definition}[Forgetting]
    
Let $\mathcal{T} = (\S; T_\a, T_\b, T_\g) $ be a bridge trisection of a surface link $ L $ and $ \d \in T_\g $ a bridge arc.

The result $\mathcal{T}_\d$ of \emph{forgetting} $\d$ (sometimes \emph{weakly reducing along} $\d$) is a triple of tangles $(\S_\d; f_\d(T_\a), f_\d(T_\b), f_\d(T_\g))$, obtained as follows:

\begin{itemize}
    \item $f_\d(T_\g)$ is $T_\g \bs \d$.
    \item $f_\d(T_\a)$ is the result of resolving $T_\a$ along $\d$.
    \item $f_\d(T_\b)$ is the result of resolving $T_\b$ along $\d$.
\end{itemize}

\end{definition}

Once again, we can realize this forgetting on the level of shadows when $\d$ is primitive with respect to each of the other tangles.

A well-known bridge trisection of the spun trefoil is weakly-reducible.

\begin{figure}[h]
    \centering

\tikzset{every picture/.style={line width=0.75pt}} 

\begin{tikzpicture}[x=1.25pt,y=1.25pt,yscale=-1,xscale=1]

\draw [color={rgb, 255:red, 126; green, 211; blue, 33 }  ,draw opacity=1 ]   (300,80.02) .. controls (300.13,70.27) and (310.13,70.27) .. (310,80.02) ;
\draw [color={rgb, 255:red, 126; green, 211; blue, 33 }  ,draw opacity=1 ]   (250,80.02) .. controls (249.83,47.35) and (244.83,18.83) .. (258.56,28) ;
\draw [color={rgb, 255:red, 126; green, 211; blue, 33 }  ,draw opacity=1 ]   (320,80.02) .. controls (319.89,59.39) and (311.76,41.69) .. (298.8,30.09) .. controls (285.84,18.48) and (251.67,26) .. (260.56,34.22) ;
\draw [color={rgb, 255:red, 126; green, 211; blue, 33 }  ,draw opacity=1 ]   (261.22,41.56) .. controls (253.18,36.78) and (269.44,35.33) .. (262.33,30.22) ;
\draw [color={rgb, 255:red, 126; green, 211; blue, 33 }  ,draw opacity=1 ]   (264.56,62.67) .. controls (268.11,60) and (262.11,56.22) .. (260.56,53.56) .. controls (259,50.89) and (259.44,46.44) .. (260.78,44.67) .. controls (262.11,42.89) and (267.89,39.78) .. (263.2,36.89) ;
\draw [color={rgb, 255:red, 126; green, 211; blue, 33 }  ,draw opacity=1 ]   (290,80) .. controls (290.11,40.44) and (271,50) .. (264.11,43.78) ;
\draw [color={rgb, 255:red, 126; green, 211; blue, 33 }  ,draw opacity=1 ]   (280,79.8) .. controls (280.78,48.78) and (269.62,51.14) .. (264.56,56) ;
\draw [color={rgb, 255:red, 126; green, 211; blue, 33 }  ,draw opacity=1 ]   (264.56,69.11) .. controls (271.89,64.44) and (254.56,62.67) .. (261.67,57.56) ;
\draw [color={rgb, 255:red, 126; green, 211; blue, 33 }  ,draw opacity=1 ]   (270,80.02) .. controls (270.11,70.89) and (254.11,70.44) .. (262.11,65.56) ;
\draw [color={rgb, 255:red, 126; green, 211; blue, 33 }  ,draw opacity=1 ]   (260,80) .. controls (260.11,76.67) and (259,73.78) .. (260.78,72) ;
\draw    (20,80) -- (110,80) ;
\draw [color={rgb, 255:red, 255; green, 0; blue, 0 }  ,draw opacity=1 ]   (30,80) .. controls (30.13,70.25) and (40.13,70.25) .. (40,80) ;
\draw [color={rgb, 255:red, 255; green, 0; blue, 0 }  ,draw opacity=1 ]   (70,80) .. controls (70.13,70.25) and (80.13,70.25) .. (80,80) ;
\draw [color={rgb, 255:red, 255; green, 0; blue, 0 }  ,draw opacity=1 ]   (50,80) .. controls (50.13,40) and (100.13,40) .. (100,80) ;
\draw [color={rgb, 255:red, 255; green, 0; blue, 0 }  ,draw opacity=1 ]   (60,80) .. controls (60,56) and (90.5,55.5) .. (90,80) ;
\draw    (130,80) -- (220,80) ;
\draw [color={rgb, 255:red, 0; green, 0; blue, 255 }  ,draw opacity=1 ]   (140,80) .. controls (140.13,70.25) and (150.13,70.25) .. (150,80) ;
\draw [color={rgb, 255:red, 0; green, 0; blue, 255 }  ,draw opacity=1 ]   (170,80) .. controls (170.13,70.25) and (180.13,70.25) .. (180,80) ;
\draw [color={rgb, 255:red, 0; green, 0; blue, 255 }  ,draw opacity=1 ]   (160,80) .. controls (160,56) and (190.5,55.5) .. (190,80) ;
\draw [color={rgb, 255:red, 0; green, 0; blue, 255 }  ,draw opacity=1 ]   (200,80) .. controls (200.13,70.25) and (210.13,70.25) .. (210,80) ;
\draw    (240,80) -- (330,80) ;

\draw (65,90) node {$T_{\alpha }$};
\draw (175,90) node  {$T_{\beta }$};
\draw (285,90) node {$T_{\gamma }$};
\draw (32.5,63) node [anchor=north west][inner sep=0.75pt]  {$\textcolor[rgb]{1,0,0}{\varepsilon_\a }$};
\draw (142,63) node [anchor=north west][inner sep=0.75pt]  [font=\small]  {$\textcolor[rgb]{0,0,1}{\varepsilon_\b}$};
\draw (302,63) node [anchor=north west][inner sep=0.75pt]  [font=\small]  {$\textcolor[rgb]{0.49,0.83,0.13}{\delta }$};

\end{tikzpicture}

    \caption{A tri-plane diagram for a weakly-reducible bridge trisection of $\operatorname{spin}(3_1)$}
    \label{fig:SpunTrefTri}
\end{figure}
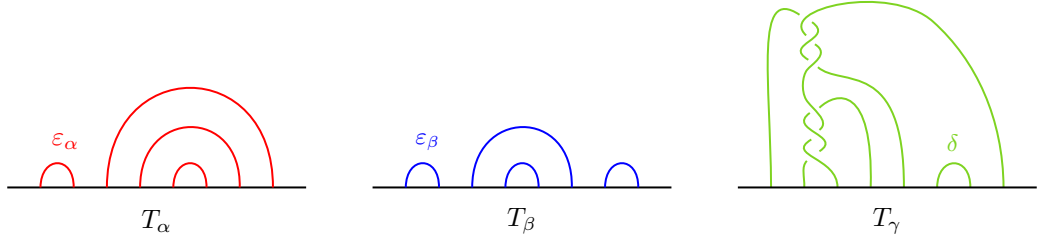

\begin{eg}

The 4-bridge trisection of $\operatorname{spin}(3_1)$ (due to \cite{wolf}) represented in Figure \ref{fig:SpunTrefTri} is weakly-reducible. Forgetting $\d$ yields an unlink of two spheres.
    
\end{eg}

\begin{figure}[h]
    \centering

\tikzset{every picture/.style={line width=0.75pt}} 

\begin{tikzpicture}[x=1pt,y=1pt,yscale=-1,xscale=1]

\draw  [fill={rgb, 255:red, 0; green, 0; blue, 0 }  ,fill opacity=1 ] (27.86,100) .. controls (27.86,98.9) and (28.75,98) .. (29.86,98) .. controls (30.96,98) and (31.86,98.9) .. (31.86,100) .. controls (31.86,101.1) and (30.96,102) .. (29.86,102) .. controls (28.75,102) and (27.86,101.1) .. (27.86,100) -- cycle ;
\draw  [fill={rgb, 255:red, 0; green, 0; blue, 0 }  ,fill opacity=1 ] (58.33,100) .. controls (58.33,98.9) and (59.23,98) .. (60.33,98) .. controls (61.44,98) and (62.33,98.9) .. (62.33,100) .. controls (62.33,101.1) and (61.44,102) .. (60.33,102) .. controls (59.23,102) and (58.33,101.1) .. (58.33,100) -- cycle ;
\draw  [fill={rgb, 255:red, 0; green, 0; blue, 0 }  ,fill opacity=1 ] (88,99.67) .. controls (88,98.56) and (88.9,97.67) .. (90,97.67) .. controls (91.1,97.67) and (92,98.56) .. (92,99.67) .. controls (92,100.77) and (91.1,101.67) .. (90,101.67) .. controls (88.9,101.67) and (88,100.77) .. (88,99.67) -- cycle ;
\draw  [fill={rgb, 255:red, 0; green, 0; blue, 0 }  ,fill opacity=1 ] (118,99.67) .. controls (118,98.56) and (118.9,97.67) .. (120,97.67) .. controls (121.1,97.67) and (122,98.56) .. (122,99.67) .. controls (122,100.77) and (121.1,101.67) .. (120,101.67) .. controls (118.9,101.67) and (118,100.77) .. (118,99.67) -- cycle ;
\draw  [fill={rgb, 255:red, 0; green, 0; blue, 0 }  ,fill opacity=1 ] (148,100) .. controls (148,98.9) and (148.9,98) .. (150,98) .. controls (151.1,98) and (152,98.9) .. (152,100) .. controls (152,101.1) and (151.1,102) .. (150,102) .. controls (148.9,102) and (148,101.1) .. (148,100) -- cycle ;
\draw  [fill={rgb, 255:red, 0; green, 0; blue, 0 }  ,fill opacity=1 ] (178,100) .. controls (178,98.9) and (178.9,98) .. (180,98) .. controls (181.1,98) and (182,98.9) .. (182,100) .. controls (182,101.1) and (181.1,102) .. (180,102) .. controls (178.9,102) and (178,101.1) .. (178,100) -- cycle ;
\draw  [fill={rgb, 255:red, 0; green, 0; blue, 0 }  ,fill opacity=1 ] (208,100) .. controls (208,98.9) and (208.9,98) .. (210,98) .. controls (211.1,98) and (212,98.9) .. (212,100) .. controls (212,101.1) and (211.1,102) .. (210,102) .. controls (208.9,102) and (208,101.1) .. (208,100) -- cycle ;
\draw  [fill={rgb, 255:red, 0; green, 0; blue, 0 }  ,fill opacity=1 ] (238.33,100) .. controls (238.33,98.9) and (239.23,98) .. (240.33,98) .. controls (241.44,98) and (242.33,98.9) .. (242.33,100) .. controls (242.33,101.1) and (241.44,102) .. (240.33,102) .. controls (239.23,102) and (238.33,101.1) .. (238.33,100) -- cycle ;
\draw [color={rgb, 255:red, 208; green, 2; blue, 27 }  ,draw opacity=1 ]   (29.86,98) .. controls (40.17,87.98) and (50.69,88) .. (60.33,98) ;
\draw [color={rgb, 255:red, 208; green, 2; blue, 27 }  ,draw opacity=1 ]   (150,98) .. controls (160.31,87.98) and (170.36,88) .. (180,98) ;
\draw [color={rgb, 255:red, 208; green, 2; blue, 27 }  ,draw opacity=1 ]   (120,97.67) .. controls (140.4,68) and (190.2,68) .. (210,98) ;
\draw [color={rgb, 255:red, 208; green, 2; blue, 27 }  ,draw opacity=1 ]   (90,97.67) .. controls (110.22,48.22) and (220.22,48.22) .. (240.33,98) ;
\draw [color={rgb, 255:red, 74; green, 144; blue, 226 }  ,draw opacity=1 ]   (120,101.67) .. controls (130.07,112.14) and (140.36,112) .. (150,102) ;
\draw [color={rgb, 255:red, 74; green, 144; blue, 226 }  ,draw opacity=1 ]   (90,101.67) .. controls (110.21,131.86) and (160.07,131.86) .. (180,102) ;
\draw [color={rgb, 255:red, 74; green, 144; blue, 226 }  ,draw opacity=1 ]   (29.86,102) .. controls (39.93,112.48) and (50.69,112) .. (60.33,102) ;
\draw [color={rgb, 255:red, 74; green, 144; blue, 226 }  ,draw opacity=1 ]   (210,102) .. controls (220.07,112.48) and (230.69,112) .. (240.33,102) ;
\draw [color={rgb, 255:red, 126; green, 211; blue, 33 }  ,draw opacity=1 ]   (182,100) -- (208,100) ;
\draw [color={rgb, 255:red, 126; green, 211; blue, 33 }  ,draw opacity=1 ]   (60.33,98) .. controls (60.3,83.1) and (100.1,85.1) .. (100,100) .. controls (100.14,120.43) and (50.14,120.14) .. (50,100) .. controls (50.3,64.9) and (120.5,72.97) .. (120,97.67) ;
\draw [color={rgb, 255:red, 126; green, 211; blue, 33 }  ,draw opacity=1 ]   (90,101.67) .. controls (89.86,114.86) and (55.14,114.71) .. (55,100) .. controls (55.14,75.14) and (110.14,80.14) .. (110,100) .. controls (110,125.14) and (160.14,120.29) .. (160,100) .. controls (160.29,55.14) and (35.1,50.1) .. (35,100) .. controls (35.11,150.22) and (240.22,150) .. (240.33,102) ;
\draw [color={rgb, 255:red, 126; green, 211; blue, 33 }  ,draw opacity=1 ]   (150,98) .. controls (150,68.9) and (40.2,55.1) .. (40,100) .. controls (40,135.1) and (169.8,137.3) .. (170,100) .. controls (170.4,45.5) and (30.2,45.1) .. (29.86,98) ;

\draw (191,102.4) node [anchor=north west][inner sep=0.75pt] {$\textcolor[rgb]{0.49,0.83,0.13}{\delta^* }$};
\draw (30,105) node [anchor=north east][inner sep=0.75pt]    {$\textcolor[rgb]{0.56,0.07,1}{\varepsilon^* }$};

\end{tikzpicture}
    
    \caption{A shadow diagram for the bridge trisection in Figure \ref{fig:SpunTrefTri}}
    \label{fig:SpunTrefShadows}
\end{figure}
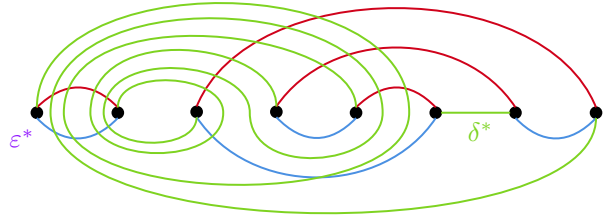

To see that this bridge trisection is weakly-reducible, we examine a shadow diagram for the trisection (Figure \ref{fig:SpunTrefShadows}), marking the weak-reducing shadows $\d^*$ and $\e^*$.

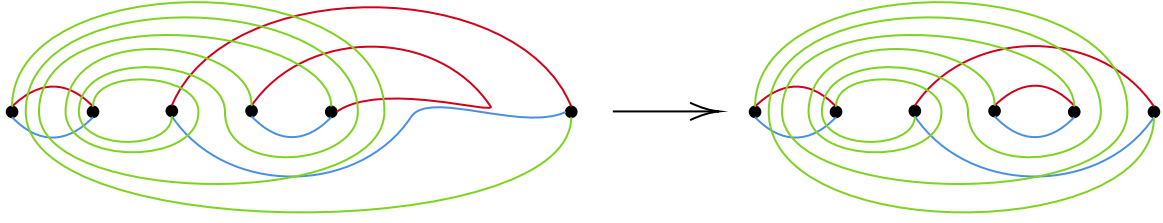
\begin{figure}[h]
    \centering

\tikzset{every picture/.style={line width=0.75pt}} 

\begin{tikzpicture}[x=1pt,y=1pt,yscale=-1,xscale=1]

\draw  [fill={rgb, 255:red, 0; green, 0; blue, 0 }  ,fill opacity=1 ] (45.86,82.16) .. controls (45.86,81.06) and (46.75,80.16) .. (47.86,80.16) .. controls (48.96,80.16) and (49.86,81.06) .. (49.86,82.16) .. controls (49.86,83.27) and (48.96,84.16) .. (47.86,84.16) .. controls (46.75,84.16) and (45.86,83.27) .. (45.86,82.16) -- cycle ;
\draw  [fill={rgb, 255:red, 0; green, 0; blue, 0 }  ,fill opacity=1 ] (76.33,82.16) .. controls (76.33,81.06) and (77.23,80.16) .. (78.33,80.16) .. controls (79.44,80.16) and (80.33,81.06) .. (80.33,82.16) .. controls (80.33,83.27) and (79.44,84.16) .. (78.33,84.16) .. controls (77.23,84.16) and (76.33,83.27) .. (76.33,82.16) -- cycle ;
\draw  [fill={rgb, 255:red, 0; green, 0; blue, 0 }  ,fill opacity=1 ] (106,81.83) .. controls (106,80.73) and (106.9,79.83) .. (108,79.83) .. controls (109.1,79.83) and (110,80.73) .. (110,81.83) .. controls (110,82.93) and (109.1,83.83) .. (108,83.83) .. controls (106.9,83.83) and (106,82.93) .. (106,81.83) -- cycle ;
\draw  [fill={rgb, 255:red, 0; green, 0; blue, 0 }  ,fill opacity=1 ] (136,81.83) .. controls (136,80.73) and (136.9,79.83) .. (138,79.83) .. controls (139.1,79.83) and (140,80.73) .. (140,81.83) .. controls (140,82.93) and (139.1,83.83) .. (138,83.83) .. controls (136.9,83.83) and (136,82.93) .. (136,81.83) -- cycle ;
\draw  [fill={rgb, 255:red, 0; green, 0; blue, 0 }  ,fill opacity=1 ] (166,82.16) .. controls (166,81.06) and (166.9,80.16) .. (168,80.16) .. controls (169.1,80.16) and (170,81.06) .. (170,82.16) .. controls (170,83.27) and (169.1,84.16) .. (168,84.16) .. controls (166.9,84.16) and (166,83.27) .. (166,82.16) -- cycle ;
\draw  [fill={rgb, 255:red, 0; green, 0; blue, 0 }  ,fill opacity=1 ] (256.33,82.16) .. controls (256.33,81.06) and (257.23,80.16) .. (258.33,80.16) .. controls (259.44,80.16) and (260.33,81.06) .. (260.33,82.16) .. controls (260.33,83.27) and (259.44,84.16) .. (258.33,84.16) .. controls (257.23,84.16) and (256.33,83.27) .. (256.33,82.16) -- cycle ;
\draw [color={rgb, 255:red, 208; green, 2; blue, 27 }  ,draw opacity=1 ]   (47.86,80.16) .. controls (58.17,70.15) and (68.69,70.16) .. (78.33,80.16) ;
\draw [color={rgb, 255:red, 208; green, 2; blue, 27 }  ,draw opacity=1 ]   (138,79.83) .. controls (158.4,50.16) and (208.2,50.16) .. (228,80.16) .. controls (229.7,83.3) and (186.5,70.5) .. (170,82.16) ;
\draw [color={rgb, 255:red, 208; green, 2; blue, 27 }  ,draw opacity=1 ]   (108,79.83) .. controls (128.22,30.39) and (238.22,30.39) .. (258.33,80.16) ;
\draw [color={rgb, 255:red, 74; green, 144; blue, 226 }  ,draw opacity=1 ]   (138,83.83) .. controls (148.07,94.31) and (158.36,94.16) .. (168,84.16) ;
\draw [color={rgb, 255:red, 74; green, 144; blue, 226 }  ,draw opacity=1 ]   (108,83.83) .. controls (128.21,114.02) and (178.07,114.02) .. (198,84.16) .. controls (204.7,73.5) and (240.1,89.9) .. (256.33,82.16) ;
\draw [color={rgb, 255:red, 74; green, 144; blue, 226 }  ,draw opacity=1 ]   (47.86,84.16) .. controls (57.93,94.64) and (68.69,94.16) .. (78.33,84.16) ;
\draw [color={rgb, 255:red, 126; green, 211; blue, 33 }  ,draw opacity=1 ]   (78.33,80.16) .. controls (78.3,65.26) and (118.1,67.26) .. (118,82.16) .. controls (118.14,102.59) and (68.14,102.31) .. (68,82.16) .. controls (68.3,47.06) and (138.5,55.13) .. (138,79.83) ;
\draw [color={rgb, 255:red, 126; green, 211; blue, 33 }  ,draw opacity=1 ]   (108,83.83) .. controls (107.86,97.02) and (73.14,96.88) .. (73,82.16) .. controls (73.14,57.31) and (128.14,62.31) .. (128,82.16) .. controls (128,107.31) and (178.14,102.45) .. (178,82.16) .. controls (178.29,37.31) and (53.1,32.26) .. (53,82.16) .. controls (53.11,132.39) and (258.22,132.16) .. (258.33,84.16) ;
\draw [color={rgb, 255:red, 126; green, 211; blue, 33 }  ,draw opacity=1 ]   (168,80.16) .. controls (168,51.06) and (58.2,37.26) .. (58,82.16) .. controls (58,117.26) and (187.8,119.46) .. (188,82.16) .. controls (188.4,27.66) and (48.2,27.26) .. (47.86,80.16) ;
\draw  [fill={rgb, 255:red, 0; green, 0; blue, 0 }  ,fill opacity=1 ] (325.52,82.16) .. controls (325.52,81.06) and (326.42,80.16) .. (327.52,80.16) .. controls (328.63,80.16) and (329.52,81.06) .. (329.52,82.16) .. controls (329.52,83.27) and (328.63,84.16) .. (327.52,84.16) .. controls (326.42,84.16) and (325.52,83.27) .. (325.52,82.16) -- cycle ;
\draw  [fill={rgb, 255:red, 0; green, 0; blue, 0 }  ,fill opacity=1 ] (356,82.16) .. controls (356,81.06) and (356.9,80.16) .. (358,80.16) .. controls (359.1,80.16) and (360,81.06) .. (360,82.16) .. controls (360,83.27) and (359.1,84.16) .. (358,84.16) .. controls (356.9,84.16) and (356,83.27) .. (356,82.16) -- cycle ;
\draw  [fill={rgb, 255:red, 0; green, 0; blue, 0 }  ,fill opacity=1 ] (385.67,81.83) .. controls (385.67,80.73) and (386.56,79.83) .. (387.67,79.83) .. controls (388.77,79.83) and (389.67,80.73) .. (389.67,81.83) .. controls (389.67,82.93) and (388.77,83.83) .. (387.67,83.83) .. controls (386.56,83.83) and (385.67,82.93) .. (385.67,81.83) -- cycle ;
\draw  [fill={rgb, 255:red, 0; green, 0; blue, 0 }  ,fill opacity=1 ] (415.67,81.83) .. controls (415.67,80.73) and (416.56,79.83) .. (417.67,79.83) .. controls (418.77,79.83) and (419.67,80.73) .. (419.67,81.83) .. controls (419.67,82.93) and (418.77,83.83) .. (417.67,83.83) .. controls (416.56,83.83) and (415.67,82.93) .. (415.67,81.83) -- cycle ;
\draw  [fill={rgb, 255:red, 0; green, 0; blue, 0 }  ,fill opacity=1 ] (445.67,82.16) .. controls (445.67,81.06) and (446.56,80.16) .. (447.67,80.16) .. controls (448.77,80.16) and (449.67,81.06) .. (449.67,82.16) .. controls (449.67,83.27) and (448.77,84.16) .. (447.67,84.16) .. controls (446.56,84.16) and (445.67,83.27) .. (445.67,82.16) -- cycle ;
\draw  [fill={rgb, 255:red, 0; green, 0; blue, 0 }  ,fill opacity=1 ] (475.67,82.16) .. controls (475.67,81.06) and (476.56,80.16) .. (477.67,80.16) .. controls (478.77,80.16) and (479.67,81.06) .. (479.67,82.16) .. controls (479.67,83.27) and (478.77,84.16) .. (477.67,84.16) .. controls (476.56,84.16) and (475.67,83.27) .. (475.67,82.16) -- cycle ;
\draw [color={rgb, 255:red, 208; green, 2; blue, 27 }  ,draw opacity=1 ]   (327.52,80.16) .. controls (337.83,70.15) and (348.36,70.16) .. (358,80.16) ;
\draw [color={rgb, 255:red, 208; green, 2; blue, 27 }  ,draw opacity=1 ]   (417.67,79.83) .. controls (427.98,69.81) and (438.02,69.83) .. (447.67,79.83) ;
\draw [color={rgb, 255:red, 208; green, 2; blue, 27 }  ,draw opacity=1 ]   (387.67,79.83) .. controls (407.63,49.88) and (457.38,49.88) .. (477.67,80.16) ;
\draw [color={rgb, 255:red, 74; green, 144; blue, 226 }  ,draw opacity=1 ]   (417.67,83.83) .. controls (427.74,94.31) and (438.02,94.16) .. (447.67,84.16) ;
\draw [color={rgb, 255:red, 74; green, 144; blue, 226 }  ,draw opacity=1 ]   (387.67,83.83) .. controls (407.88,114.02) and (457.74,114.02) .. (477.67,84.16) ;
\draw [color={rgb, 255:red, 74; green, 144; blue, 226 }  ,draw opacity=1 ]   (327.52,84.16) .. controls (337.6,94.64) and (348.36,94.16) .. (358,84.16) ;
\draw [color={rgb, 255:red, 126; green, 211; blue, 33 }  ,draw opacity=1 ]   (358,80.16) .. controls (357.97,65.26) and (397.77,67.26) .. (397.67,82.16) .. controls (397.81,102.59) and (347.81,102.31) .. (347.67,82.16) .. controls (347.97,47.06) and (418.17,55.13) .. (417.67,79.83) ;
\draw [color={rgb, 255:red, 126; green, 211; blue, 33 }  ,draw opacity=1 ]   (387.67,83.83) .. controls (387.52,97.02) and (352.81,96.88) .. (352.67,82.16) .. controls (352.81,57.31) and (407.81,62.31) .. (407.67,82.16) .. controls (407.67,107.31) and (457.81,102.45) .. (457.67,82.16) .. controls (457.95,37.31) and (332.77,32.26) .. (332.67,82.16) .. controls (332.78,132.39) and (477.56,132.16) .. (477.67,84.16) ;
\draw [color={rgb, 255:red, 126; green, 211; blue, 33 }  ,draw opacity=1 ]   (447.67,80.16) .. controls (447.67,51.06) and (337.87,37.26) .. (337.67,82.16) .. controls (337.67,117.26) and (467.47,119.46) .. (467.67,82.16) .. controls (468.07,27.66) and (327.87,27.26) .. (327.52,80.16) ;
\draw    (274,82) -- (312,82) ;
\draw [shift={(314,82)}, rotate = 180] [color={rgb, 255:red, 0; green, 0; blue, 0 }  ][line width=0.75]    (10.93,-3.29) .. controls (6.95,-1.4) and (3.31,-0.3) .. (0,0) .. controls (3.31,0.3) and (6.95,1.4) .. (10.93,3.29)   ;

\end{tikzpicture}

    \caption{The shadow diagram $f_\d(\S;T_\a^*,T_\b^*,T_\g^*)$, with an isotopy}
    \label{fig:SpunTrefForget}
\end{figure}

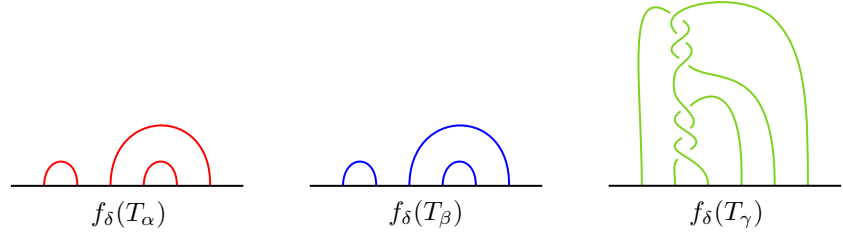
\begin{figure}[h]
    \centering

\tikzset{every picture/.style={line width=0.75pt}} 

\begin{tikzpicture}[x=1.25pt,y=1.25pt,yscale=-1,xscale=1]

\draw [color={rgb, 255:red, 126; green, 211; blue, 33 }  ,draw opacity=1 ]   (230,80) .. controls (229.83,47.33) and (224.83,18.81) .. (238.56,27.98) ;
\draw [color={rgb, 255:red, 126; green, 211; blue, 33 }  ,draw opacity=1 ]   (280,79.98) .. controls (279.89,59.35) and (280.17,40.15) .. (270,29.98) .. controls (259.83,19.81) and (231.67,25.98) .. (240.56,34.2) ;
\draw [color={rgb, 255:red, 126; green, 211; blue, 33 }  ,draw opacity=1 ]   (241.22,41.54) .. controls (233.18,36.77) and (249.44,35.31) .. (242.33,30.2) ;
\draw [color={rgb, 255:red, 126; green, 211; blue, 33 }  ,draw opacity=1 ]   (244.56,62.65) .. controls (248.11,59.98) and (242.11,56.2) .. (240.56,53.54) .. controls (239,50.87) and (239.44,46.43) .. (240.78,44.65) .. controls (242.11,42.87) and (247.89,39.76) .. (243.2,36.87) ;
\draw [color={rgb, 255:red, 126; green, 211; blue, 33 }  ,draw opacity=1 ]   (270,79.98) .. controls (270.11,40.43) and (251,49.98) .. (244.11,43.76) ;
\draw [color={rgb, 255:red, 126; green, 211; blue, 33 }  ,draw opacity=1 ]   (260,79.78) .. controls (260.78,48.76) and (249.62,51.12) .. (244.56,55.98) ;
\draw [color={rgb, 255:red, 126; green, 211; blue, 33 }  ,draw opacity=1 ]   (244.56,69.09) .. controls (251.89,64.43) and (234.56,62.65) .. (241.67,57.54) ;
\draw [color={rgb, 255:red, 126; green, 211; blue, 33 }  ,draw opacity=1 ]   (250,80) .. controls (250.11,70.87) and (234.11,70.43) .. (242.11,65.54) ;
\draw [color={rgb, 255:red, 126; green, 211; blue, 33 }  ,draw opacity=1 ]   (240,79.98) .. controls (240.11,76.65) and (239,73.76) .. (240.78,71.98) ;
\draw [color={rgb, 255:red, 255; green, 0; blue, 0 }  ,draw opacity=1 ]   (49.99,80) .. controls (50.12,70.25) and (60.12,70.25) .. (59.99,80) ;
\draw [color={rgb, 255:red, 255; green, 0; blue, 0 }  ,draw opacity=1 ]   (79.99,80) .. controls (80.12,70.25) and (90.12,70.25) .. (89.99,80) ;
\draw [color={rgb, 255:red, 255; green, 0; blue, 0 }  ,draw opacity=1 ]   (69.99,80) .. controls (69.99,56) and (100.49,55.5) .. (99.99,80) ;
\draw [color={rgb, 255:red, 0; green, 0; blue, 255 }  ,draw opacity=1 ]   (140,80) .. controls (140.13,70.25) and (150.13,70.25) .. (150,80) ;
\draw [color={rgb, 255:red, 0; green, 0; blue, 255 }  ,draw opacity=1 ]   (170,80) .. controls (170.13,70.25) and (180.13,70.25) .. (180,80) ;
\draw [color={rgb, 255:red, 0; green, 0; blue, 255 }  ,draw opacity=1 ]   (160,80) .. controls (160,56) and (190.5,55.5) .. (190,80) ;
\draw    (220,79.98) -- (290,79.98) ;
\draw    (40,80) -- (110,80) ;
\draw    (130,80) -- (200,80) ;

\draw (75.5,88.5) node   {$f_{\delta }( T_{\alpha })$};
\draw (165,88.5) node   {$f_{\delta }( T_{\beta })$};
\draw (255.5,88.5) node  {$f_{\delta }( T_{\gamma })$};

\end{tikzpicture}
    \caption{A tri-plane diagram for $f_\d(\S;T_\a,T_\b,T_\g)$}
    \label{fig:SpunTrefForgTri}
\end{figure}

As we can see, $\d^*$ is primitive with respect to each of $T_\a$ and $T_\b$, and so we can forget $\d^*$ and obtain a new shadow diagram $f_\d(\S;T_\a^*,T_\b^*,T_\g^*)$, Figure \ref{fig:SpunTrefForget}. This new shadow diagram induces a tri-plane diagram, and once we verify that all pairwise unions of the tangles form unlinks, we see that this is in fact a bridge trisection. In Figure \ref{fig:SpunTrefForgTri}, we can verify that each pair of tangles does in fact form an unlink, and thus gives a bridge trisection for some surface link $L$. The bridge trisection this diagram represents has parameters $(3;3,2,2)$, and therefore $\chi(L) = 4$. We can then apply Proposition 5.1 of \cite{jmmzsolids} to see that $L$ is an unlink of two $S^2$'s. This is consistent with $\operatorname{spin}(3_1)$'s structure as a ribbon sphere, as it can be obtained by attaching a single tube to a pair of unlinked spheres. This example illustrates and anticipates our main finding regarding 4-bridge weakly-reducible trisections (Theorem \ref{thm:4Br}).

\section{Technical Lemmas}

We will assemble the proofs of our main results by way of a few technical lemmas concerning arc surgery on on bridge trisections, the uniqueness of bridge disks within certain building blocks of bridge splittings, and weakly-reducible splittings of classical unlinks. In the end, we are able to conclude that a 4-bridge trisection in $S^4$ satisfies the needed conditions to apply Theorem \ref{thm:ElemOp}.






First, we will examine bridge splittings of trivial tangles. The following theorem of Hayashi and Shimokawa (here paraphrased to match our notation), when applied in the case of a $(b+1)$-bridge splitting of a $b$-component tangle in $S^2 \x I$, tells us that the splitting is perturbed (``stabilized" in the language of \cite{TrivialHS}).

\begin{thm}[\cite{TrivialHS} Theorem 1.1]
\label{thm:HS}
    Let $C$ be a compression body, $T$ a trivial tangle in $C$, and $H$ a Heegaard splitting of $(C,T)$. Then $H$ is either the trivial Heegaard splitting or else is perturbed or contains a 1-bridge unlinked component.
\end{thm}

When our splitting is perturbed, we would like to guarantee that if one side contains only one trivial arc (and $b$ vertical arcs), then that arc is primitive with respect to the opposite tangle. This becomes important in realizing weak reduction on trisections geometrically.

\begin{lemma}\label{lemma:vert}
    Let $T$ be a tangle in $S^2 \x I$ (viewed as a compression-body with boundary components $\S_- \cong \S_+ \cong S^2$) consisting of a single arc $a$ with endpoints in $\S_+$ and $n$ vertical arcs $v_1, ..., v_n$. Then there is a unique bridge disk $B_a$ for $a$.
\end{lemma}

\begin{proof}
    Let $B_a$ and $B'_a$ be two bridge disks for $a$, and position them to minimize $|B_a \cap B'_a|$. Note that $B_a \cap B'_a$ is a collection of simple closed curves and arcs with endpoints in $\S_+$. If $B_a \cap B'_a$ contains any closed curves, we take such a curve $c$ bounding an innermost disk $D$ in $B_a$ and note that $c$ bounds a disk $D'$ in $B'_a$ as well. The sphere $D \cup D'$ bounds a ball away from all vertical arcs of $T$, as $(S^2 \x I) \bs (\cup_{i=1}^n \nu(v_i))$ contains no essential spheres. We isotope $D'$ across this ball until it lies in a $\nu(D)$, thus removing all intersections with $B_a$ in the interior of $D'$. We extend this isotopy to $B'_a$. 
    
    Then, since $\nu(D)$ is a 3-ball, we further isotope $B'_a$ across $D$, removing the curve $c$ of intersection between $B'_a$ and $B_a$. In this way, we could decrease $|B_a \cap B'_a|$, so $B_a \cap B'_a$ contains no curves.

    Suppose now that $B_a \cap B'_a$ contains an arc, which must have endpoints in $\S_+$. Take such an arc $b$ cobounding with some arc $b^*$ in $B_a \cap \S_+$ an outermost disk $D$ in $B_a$. As above, $b$ also cobounds a disk $D'$ in $B'_a$ with an arc $b'^*$ of $B'_a \cap \S_+$. Since $\S_+$ is a sphere, $b^* \cup b'^*$ bounds a disk $D^*$ in $\S_+$. The union $D \cup D' \cup D^* = S$ is a sphere.
    
    Note that if $D^*$ contains any punctures (i.e., an endpoint of any $v_i$), then $v_i$ intersects $S$ away from $\S_+$, as its other endpoint lies in $\S_-$. However, each of $D$ and $D'$ is contained in a bridge disk for $a$, and these are by definition disjoint from $T \bs a$. So $S$ is a sphere in $(S^2 \x I) \bs (\cup_{i=1}^n \nu(v_i))$, and therefore bounds a ball. We carry out a similar process to that for curves of intersection and remove the arc $b$ of $B_a \cap B'_a$.

    In this way, we see that $B_a$ and $B'_a$ intersect only in their boundary. The disk $B_a \cup B'_a$ has its boundary in $\S_+$, and so there is a disk $B^* \subset \S_+$ that contains no punctures, as in the case of $D^*$, above. This disk cobounds a ball with $B_a \cup B'_a$. We can then isotopy $B_a$ across this ball to coincide with $B'_a$, and so $B_a$ and $B'_a$ are isotopic in $(S^2 \x I) \bs (\cup_{i=1}^n \nu(v_i))$.
    
\end{proof}

These two results together let us conclude that, if we forget some arc $\d$ in a bridge trisection and obtain a bridge trisection, then our original trisection is weakly-reducible, and we can realize the forgetting in the shadow diagram.

\begin{cor}
\label{cor:prim}
    Let $\mathcal{T} = (\S;T_\a,T_\b,T_\g)$ be a bridge trisection where $T_\a \cup \overline{T_\b}$ contains a 1-bridge component disjoint from some arc $\d$ in $T_\g$. Suppose that after forgetting $\d$, the resulting triple of tangles $\mathcal{T}_\d = (\S_\d;f_\d(T_\a),f_\d(T_\b),f_\d(T_\g))$ is a bridge trisection. Then $\mathcal{T}$ is weakly-reducible, and the weak-reducing shadow $\d^*$ for $\d$ is primitive with respect to each of $T_\a$ and $T_\b$.
\end{cor}

\begin{proof}
    After forgetting $\d$, which we can think of as isotoping $\d$ across $\S$ in each the pairwise unions $T_\b \cup \overline{T_\g}$ and $T_\g \cup \overline{T_\b}$, we have in each of $f_\d(T_\a)$ and $f_\d(T_\b)$ a $(b+1)$-bridge splitting of a $b$-bridge tangle. Since $\mathcal{T}_\d$ is a bridge trisection by assumption, each of these tangles is trivial. Then by Theorem \ref{thm:HS}, each of $f_\d(T_\a)$ and $f_\d(T_\b)$ is perturbed, since the splitting is non-trivial and the tangles contain no closed components. The shadow $\d^*$ for $\d$ is unique by Lemma \ref{lemma:vert}, so $\d^*$ is primitive with respect to each of $T_\a$ and $T_\b$. Thus, $\d^*$ is disjoint in its interior from every $\a$ and $\b$ shadow. Since $\d$ does not intersect the 1-bridge component $\e_\a \cup \e_\b$ of $T_\a \cup \overline{T_\b}$, $\d^*$ does not intersect a shadow $\e^*$ either of $\e_\a$ or $\e_\b$ at an endpoint, and so $\mathcal{T}$ is weakly-reducible.
\end{proof}

We now turn to the specific setting of 4-bridge splittings. A technical result regarding weakly-reducible bridge splittings of classical unlinks will allow us to confirm that 4-bridge weakly-reducible trisections satisfy an important condition.

\begin{lemma}\label{lemma:Classical}

Let $\mathcal{U}$ be a weakly-reducible splitting of an unlink with four bridges and at most two components, with reducing arcs $\d$ and $\e$. Then one of the following is true:

\begin{itemize}

    \item $\mathcal{U}$ is perturbed, with one of $\d$ or $\e$ primitive with respect to the opposite tangle.

    \item One of $\d^*$ or $\e^*$ co-bounds a bridge disk on each side of the splitting.

\end{itemize}

\end{lemma}

\begin{proof}

Consider any four-bridge splitting $T_1 \cup T_2$ of an unlink $U$ (with at most two components) for tangles $T_1 \subset H_1$ and $T_2 \subset H_2$, where $H_1 \cong H_2 \cong B^3$. Let $\d \subset T_\a$ and $\e \subset T_\b$ bound bridge disks $B_\d$ and $B_\e$, respectively, such that $\S \cap B_\d = \d^*$ and $\S \cap B_\e = \e^*$. Consider $D_\d$ and $D_\e$ as in Definition \ref{def:WRClassical}. There are multiple cases to consider, where $c(U)$ is the number of components and $b_i$ is the bridge number of component $i$.

\underline{Case 1:} Suppose $c(U) = 1$. After weakly reducing along $\d$ and then $\e$, we have level surfaces $\S_\d$, $\S_\e$, and $\S_{\d\e}$ (again, as in Definition \ref{def:WRClassical}).

\

\begin{figure}[H]
    \centering
\begin{tikzpicture}
    \draw [black!50] (0,0) -- (7,0) node at (0,0) [anchor = east] {$\S$};
    \draw [thick, black] plot [smooth, tension=2] coordinates {(0.5,1) (1,1.5) (1.5,1)};
    \draw [thick, black] plot [smooth, tension=2] coordinates {(2,1) (2.5,1.5) (3,1)};
    \draw [thick, black] plot [smooth, tension=2] coordinates {(3.5,1) (4,1.5) (4.5,1)};
    \draw [thick, black] plot [smooth, tension=2] coordinates {(5,1) (5.5,1.5) (6,1)} node at (5.5,1.5) [anchor = south] {$\d$};
    \draw [thick, black] plot [smooth, tension=2] coordinates {(0.5,-1) (1,-1.5) (1.5,-1)} node at (1,-1.5) [anchor = north] {$\e$}; 
    \draw [thick, black] plot [smooth, tension=2] coordinates {(2,-1) (2.5,-1.5) (3,-1)};
    \draw [thick, black] plot [smooth, tension=2] coordinates {(3.5,-1) (4,-1.5) (4.5,-1)};
    \draw [thick, black] plot [smooth, tension=2] coordinates {(5,-1) (5.5,-1.5) (6,-1)};
    \draw [thick, black, dashed] plot [smooth, tension = 1.5] coordinates {(0.5,1) (0.5,-1)};
    \draw [thick, black, dashed] plot [smooth, tension = 1.5] coordinates {(1.5,1) (1.6,0.5) (1.9,-0.5) (2,-1)};
    \draw [thick, black, dashed] plot [smooth, tension = 1.5] coordinates {(2,1) (1.9,0.5) (1.6, -0.5) (1.5,-1)};
    \draw [thick, black, dashed] plot [smooth, tension = 1.5] coordinates {(3,1) (3.1,0.5) (3.4,-0.5) (3.5,-1)};
    \draw [thick, black, dashed] plot [smooth, tension = 1.5] coordinates {(3.5,1) (3.4,0.5) (3.1,-0.5) (3,-1)};
    \draw [thick, black, dashed] plot [smooth, tension = 1.5] coordinates {(4.5,1) (4.6,0.5) (4.9,-0.5) (5,-1)};
    \draw [thick, black, dashed] plot [smooth, tension = 1.5] coordinates {(5,1) (4.9,0.5) (4.6,-0.5) (4.5,-1)};
    \draw [thick, black, dashed] plot [smooth, tension = 1.5] coordinates {(6,1) (6,-1)};
    \draw [->] [thick, black] (7.5,0) -- (8.5,0) node at (8,0) [anchor = south] {W.R.};
    \draw [black!50] (9,1.5) -- (16,1.5) node at (16,1.5) [anchor = west] {$\S_\d$};
    \draw [black!50] (10,0) -- (15,0) node at (15,0) [anchor = west] {$\S_{\d\e}$};
    \draw [black!50] (9,-1.5) -- (16,-1.5)node at (16,-1.5) [anchor = west] {$\S_\e$};
    \draw [thick, black] plot [smooth, tension=2] coordinates {(9.5,2) (10,2.5) (10.5,2)};
    \draw [thick, black] plot [smooth, tension=2] coordinates {(11,2) (11.5,2.5) (12,2)};
    \draw [thick, black] plot [smooth, tension=2] coordinates {(12.5,2) (13,2.5) (13.5,2)};
    \draw [thick, black] plot [smooth, tension=2] coordinates {(14,-1) (14.5,-0.5) (15,-1)};
    \draw [thick, black] plot [smooth, tension=2] coordinates {(9.5,1) (10,0.5) (10.5,1)}; 
    \draw [thick, black] plot [smooth, tension=2] coordinates {(11,-2) (11.5,-2.5) (12,-2)};
    \draw [thick, black] plot [smooth, tension=2] coordinates {(12.5,-2) (13,-2.5) (13.5,-2)};
    \draw [thick, black] plot [smooth, tension=2] coordinates {(14,-2) (14.5,-2.5) (15,-2)};
    \draw [thick, black, dashed] plot [smooth, tension = 1.5] coordinates {(9.5,2) (9.5,1)};
    \draw [thick, black, dashed] plot [smooth, tension = 1.5] coordinates {(10.5,2) (10.7,0.5) (10.9,-0.5) (11,-2)};
    \draw [thick, black, dashed] plot [smooth, tension = 1.5] coordinates {(11,2) (10.9,1.6) (10.6, 1.4) (10.5,1)};
    \draw [thick, black, dashed] plot [smooth, tension = 1.5] coordinates {(12,2) (12.1,1) (12.4,0) (12.5,-2)};
    \draw [thick, black, dashed] plot [smooth, tension = 1.5] coordinates {(12.5,2) (12.4,1) (12.1,0) (12,-2)};
    \draw [thick, black, dashed] plot [smooth, tension = 1.5] coordinates {(13.5,2) (13.6,0.5) (13.8,-0.5) (14,-2)};
    \draw [thick, black, dashed] plot [smooth, tension = 1.5] coordinates {(14,-1) (13.9,-1.4) (13.6,-1.6) (13.5,-2)};
    \draw [thick, black, dashed] plot [smooth, tension = 1.5] coordinates {(15,-1) (15,-2)};
\end{tikzpicture}
    \caption{Weakly reducing $\S$ in Case 1}
    \label{fig:4BrL1.1}
\end{figure}
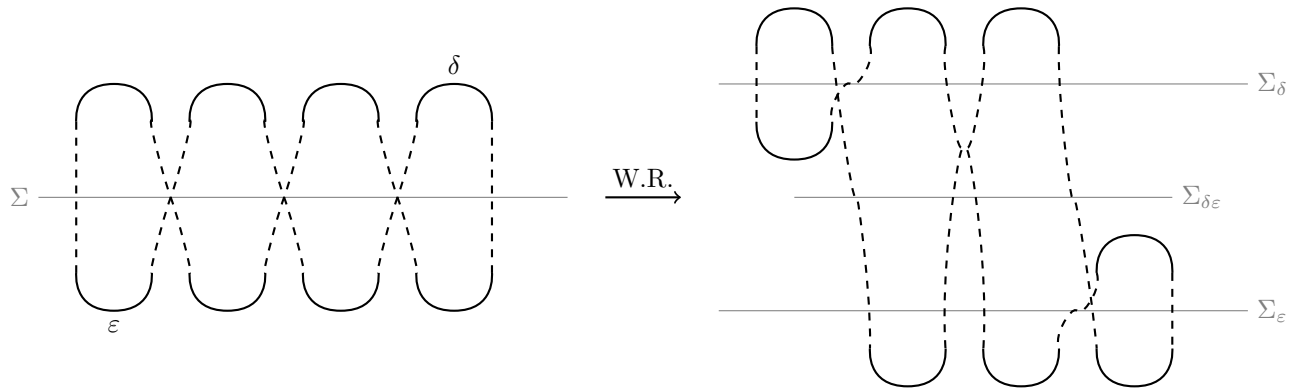

\

Note that if both thick surfaces are strongly irreducible, then the thin surface is essential, as per Theorem 1.2 of \cite{HS}. Since $U$ is unknotted, its complement contains no essential punctured surfaces, so at least one of $\S_\d$ and $\S_\e$ is weakly-reducible, say $\S_\d$. We see that $\e^*$ is still a shadow in $\S_\d$, since the surgery from $\S$ to $\S_\d$ occurs away from $\e^*$ and $B_\e$. Indeed, $\e^*$ is the shadow for the unique bridge disk for the tangle between $\S_\d$ and $\S_{\d\e}$ by Lemma \ref{lemma:vert}, as the tangle in the compression-body bounded by $\S_{\d\e}$ and $\S_\d$ consists of $\e$ and four vertical arcs. Let $\d'^*$ be a shadow in $\S_\d$ for an arc $\d'$ in the tangle opposite $\e^*$ such that $\d'^*$ is disjoint from $\e^*$.

We now weakly reduce $\S_\d$ along $\e$ and $\d'$, obtaining two new surfaces, $\S_{\d\d'}$ and $\S_{\d\d'\e}$. Again, $e^*$ is untouched by the surgery along $\d'$ and remains a shadow in $\S_{\d\d'}$. Indeed, the tangle between $\S_{\d\d'}$ and $\S_{\d\d'\e}$ consists of a single arc trivial with respect to $\S_{\d\d'}$ and two vertical arcs, so $\e^*$ is the unique shadow for $\e$. Note also that $\S_{\d\d'}$ has four marked points: the endpoints of $\e$ and of the two vertical arcs. As such, $\e^*$ has a well-defined rational (or infinite) slope in $\S_{\d\d'}$. Further, since $\S_{\d\d'\e}$ intersects $U$ in two points, it decomposes $U$ into connect summands, and since $U$ is unknotted, each of these is a punctured unknot. We may normalize such that $e^*$ has slope $\frac10$ in $\S_{\d\d'}$. Then the union of $\e^*$ with the other $\frac10$ arc in $\S_{\d\d'}$ and the opposite shadows in $\S_{\d\d'}$ is an unknot. Let $\frac{p}{q}$ denote the slope of the shadows in $\S_{\d\d'}$ opposite $\e^*$. We have an unknot exactly when $\det \begin{bmatrix}
    1 & p \\
    0 & q
\end{bmatrix} = q = 1$ \cite{Schu}. When $q = 1$, since $\e^*$ is a $\frac10$ curve, we have that $\e^*$ intersects the other shadows only in an endpoint, and therefore $\e$ is primitive with respect to the tangle with endpoints in $\S_{\d_\d'}$.

\begin{figure}[H]
    \centering

\begin{tikzpicture}

\draw [black!50] (0,1.5) -- (7,1.5) node at (0,1.5) [anchor = east] {$\S_\d$};
    \draw [black!50] (1,0) -- (6,0);
    \draw [black!50] (0,-1.5) -- (7,-1.5);
    \draw [thick, black] plot [smooth, tension=2] coordinates {(0.5,2) (1,2.5) (1.5,2)};
    \draw [thick, black] plot [smooth, tension=2] coordinates {(2,2) (2.5,2.5) (3,2)};
    \draw [thick, black] plot [smooth, tension=2] coordinates {(3.5,2) (4,2.5) (4.5,2)} node at (4,2.5) [anchor = south] {$\d'$};
    \draw [thick, black] plot [smooth, tension=2] coordinates {(5,-1) (5.5,-0.5) (6,-1)};
    \draw [thick, black] plot [smooth, tension=2] coordinates {(0.5,1) (1,0.5) (1.5,1)} node at (1,0.5) [anchor = north] {$\e$}; 
    \draw [thick, black] plot [smooth, tension=2] coordinates {(2,-2) (2.5,-2.5) (3,-2)};
    \draw [thick, black] plot [smooth, tension=2] coordinates {(3.5,-2) (4,-2.5) (4.5,-2)};
    \draw [thick, black] plot [smooth, tension=2] coordinates {(5,-2) (5.5,-2.5) (6,-2)};
    \draw [thick, black, dashed] plot [smooth, tension = 1.5] coordinates {(0.5,2) (0.5,1)};
    \draw [thick, black, dashed] plot [smooth, tension = 1.5] coordinates {(1.5,2) (1.7,0.5) (1.9,-0.5) (2,-2)};
    \draw [thick, black, dashed] plot [smooth, tension = 1.5] coordinates {(2,2) (1.9,1.6) (1.6, 1.4) (1.5,1)};
    \draw [thick, black, dashed] plot [smooth, tension = 1.5] coordinates {(3,2) (3.1,1) (3.4,0) (3.5,-2)};
    \draw [thick, black, dashed] plot [smooth, tension = 1.5] coordinates {(3.5,2) (3.4,1) (3.1,0) (3,-2)};
    \draw [thick, black, dashed] plot [smooth, tension = 1.5] coordinates {(4.5,2) (4.6,0.5) (4.8,-0.5) (5,-2)};
    \draw [thick, black, dashed] plot [smooth, tension = 1.5] coordinates {(5,-1) (4.9,-1.4) (4.6,-1.6) (4.5,-2)};
    \draw [thick, black, dashed] plot [smooth, tension = 1.5] coordinates {(6,-1) (6,-2)};
    \draw [->] [thick, black] (7.5,1) -- (8.5,1) node at (8,1) [anchor = south] {W.R.};
    \draw [black!50] (9,3) -- (16,3) node at (16,3) [anchor = west] {$\S_{\d\d'}$};
    \draw [black!50] (10,1.5) -- (14,1.5) node at (14,1.5) [anchor = west] {$\S_{\d\d'\e}$};
    \draw [black!50] (10,0) -- (15,0);
    \draw [black!50] (9,-1.5) -- (16,-1.5);
    \draw [thick, black] plot [smooth, tension=2] coordinates {(9.5,3.5) (10,4) (10.5,3.5)};
    \draw [thick, black] plot [smooth, tension=2] coordinates {(11,3.5) (11.5,4) (12,3.5)};
    \draw [thick, black] plot [smooth, tension=2] coordinates {(12.5,0.5) (13,1) (13.5,0.5)};
    \draw [thick, black] plot [smooth, tension=2] coordinates {(14,-1) (14.5,-0.5) (15,-1)};
    \draw [thick, black] plot [smooth, tension=2] coordinates {(9.5,2.5) (10,2) (10.5,2.5)}  node at (10,2) [anchor = north] {$\e$}; 
    \draw [thick, black] plot [smooth, tension=2] coordinates {(11,-2) (11.5,-2.5) (12,-2)};
    \draw [thick, black] plot [smooth, tension=2] coordinates {(12.5,-2) (13,-2.5) (13.5,-2)};
    \draw [thick, black] plot [smooth, tension=2] coordinates {(14,-2) (14.5,-2.5) (15,-2)};
    \draw [thick, black, dashed] plot [smooth, tension = 1.5] coordinates {(9.5,3.5) (9.5,2.5)};
    \draw [thick, black, dashed] plot [smooth, tension = 1.5] coordinates {(10.5,3.5) (10.7,2) (10.9,-0.5) (11,-2)};
    \draw [thick, black, dashed] plot [smooth, tension = 1.5] coordinates {(11,3.5) (10.9,3.1) (10.6, 2.9) (10.5,2.5)};
    \draw [thick, black, dashed] plot [smooth, tension = 1.5] coordinates {(12,3.5) (12.1,2) (12.4,-0.5) (12.5,-2)};
    \draw [thick, black, dashed] plot [smooth, tension = 1.5] coordinates {(12.5,0.5) (12.4,0) (12.1,-1) (12,-2)};
    \draw [thick, black, dashed] plot [smooth, tension = 1.5] coordinates {(13.5,0.5) (13.6,0) (13.8,-1) (14,-2)};
    \draw [thick, black, dashed] plot [smooth, tension = 1.5] coordinates {(14,-1) (13.9,-1.4) (13.6,-1.6) (13.5,-2)};
    \draw [thick, black, dashed] plot [smooth, tension = 1.5] coordinates {(15,-1) (15,-2)};
    
\end{tikzpicture}

    \caption{Weakly reducing $\S_\d$ in Case 1}
    \label{fig:4BRL1.2}
\end{figure}

We now have a unique bridge disk for $\e$ represented by $\e^*$ at each of the levels $\S_\d$ and $\S_{\d\d'}$, and we know that $\e$ is primitive with respect to the opposite tangle in $\S_{\d\d'}$. Let $\d''^*$ be a shadow in $\S_{\d\d'}$ for an arc $\d''$ in the tangle opposite $\e$ such that $\d''^*$ intersects $\e^*$ only at an endpoint. Note that any isotopy of $\d''^*$ in $\S_{\d\d'}$ may be realized in $\S_\d$ by isotopy and possibly shadow slides over $\d'^*$, corresponding to isotopy across the ``scar" $D_{\d'}$. Therefore, we can see that there is a shadow $\hat{\d}'^*$ in $\S_{\d}$ that intersects $\e^*$ only at an endpoint, as well, and $\hat{\d}'^*$ is the shadow of some bridge disk for the an arc $\hat{\d}'$ in the tangle opposite $\e$ in $\S_{\d}$.

Applying the same argument with regard to the scar $D_\d$ for $\d$ in $\S_\d$, we see that $\e^*$ intersects some arc $\hat{\d^*}$ in $\S$ only at a single endpoint, and again we have that $\hat{\d}^*$ is in fact the shadow of some bridge disk for an arc $\hat{\d} \in T_1$. We can subsequently slide shadows intersecting $\e^*$ in $\S$ over $\hat{\d}^*$, noting that $\d^*$ will remain unchanged, as it is disjoint from $\e^*$ by assumption, until we obtain a system of shadows for $T_1$ that intersects $\e^*$ only at its endpoints. Finally, we have a system of bridge disks for each of $T_1$ and $T_2$ that include $B_\d$ and $B_\e$, respectively, and via which $\e^*$ is primitive with respect to $T_1$.

\underline{Case 2:} Suppose $c(U) = 2$ and $b_1 = b_2 = 2$. Note first that the weak-reducing shadows $\d^*$ and $\e^*$ must arise from different components, as shadows from opposite sides of a 2-bridge splitting of an unknot will always intersect at their endpoints. We weakly reduce, obtaining new level surfaces $\S_\d$, $\S_\e$, and $\S_{\d\e}$.

\begin{figure}[h]
   \centering

\begin{tikzpicture}
    \draw [black!50] (0,0) -- (7,0) node at (0,0) [anchor = east] {$\S$};
    \draw [thick, black] plot [smooth, tension=2] coordinates {(0.5,1) (1,1.5) (1.5,1)};
    \draw [thick, purple] plot [smooth, tension=2] coordinates {(2,1) (2.5,1.5) (3,1)};
    \draw [thick, purple] plot [smooth, tension=2] coordinates {(3.5,1) (4,1.5) (4.5,1)};
    \draw [thick, black] plot [smooth, tension=2] coordinates {(5,1) (5.5,1.5) (6,1)} node at (5.5,1.5) [anchor = south] {$\d$};
    \draw [thick, purple] plot [smooth, tension=2] coordinates {(0.5,-1) (1,-1.5) (1.5,-1)} node at (1,-1.5) [anchor = north] {$\e$}; 
    \draw [thick, black] plot [smooth, tension=2] coordinates {(2,-1) (2.5,-1.5) (3,-1)};
    \draw [thick, black] plot [smooth, tension=2] coordinates {(3.5,-1) (4,-1.5) (4.5,-1)};
    \draw [thick, purple] plot [smooth, tension=2] coordinates {(5,-1) (5.5,-1.5) (6,-1)};
    \draw [thick, black, dashed] plot [smooth, tension = 0.5] coordinates {(0.5,1) (0.7,0.3) (1.8,-0.3) (2,-1)};
    \draw [thick, black, dashed] plot [smooth, tension = 0.5] coordinates {(1.5,1) (1.7,0.3) (3.3,-0.3) (3.5,-1)};
    \draw [thick, purple, dashed] plot [smooth, tension = 0.5] coordinates {(2,1) (1.8,0.5) (0.7, 0) (0.5,-1)};
    \draw [thick, purple, dashed] plot [smooth, tension = 0.5] coordinates {(3,1) (3.2,0.5) (4.8,0) (5,-1)};
    \draw [thick, purple, dashed] plot [smooth, tension = 0.5] coordinates {(3.5,1) (3.3,0.5) (1.7,0) (1.5,-1)};
    \draw [thick, purple, dashed] plot [smooth, tension = 0.5] coordinates {(4.5,1) (4.7,0.5) (5.8,0) (6,-1)};
    \draw [thick, black, dashed] plot [smooth, tension = 0.5] coordinates {(5,1) (4.8,0.3) (3.2,-0.3) (3,-1)};
    \draw [thick, black, dashed] plot [smooth, tension = 0.5] coordinates {(6,1) (5.8,0.3) (4.7,-0.3) (4.5,-1)};
    \draw [->] [thick, black] (7.5,0) -- (8.5,0) node at (8,0) [anchor = south] {W.R.};
    \draw [black!50] (9,1.5) -- (16,1.5) node at (16,1.5) [anchor = west] {$\S_\d$};
    \draw [black!50] (10,0) -- (15,0) node at (15,0) [anchor = west] {$\S_{\d\e}$};
    \draw [black!50] (9,-1.5) -- (16,-1.5) node at (16,-1.5) [anchor = west] {$\S_\e$};
    \draw [thick, black] plot [smooth, tension=2] coordinates {(9.5,2.5) (10,3) (10.5,2.5)};
    \draw [thick, purple] plot [smooth, tension=2] coordinates {(11,2.5) (11.5,3) (12,2.5)};
    \draw [thick, purple] plot [smooth, tension=2] coordinates {(12.5,2.5) (13,3) (13.5,2.5)};
    \draw [thick, black] plot [smooth, tension=2] coordinates {(14,-1) (14.5,-0.5) (15,-1)};
    \draw [thick, purple] plot [smooth, tension=2] coordinates {(9.5,1) (10,0.5) (10.5,1)}; 
    \draw [thick, black] plot [smooth, tension=2] coordinates {(11,-2.5) (11.5,-3) (12,-2.5)};
    \draw [thick, black] plot [smooth, tension=2] coordinates {(12.5,-2.5) (13,-3) (13.5,-2.5)};
    \draw [thick, purple] plot [smooth, tension=2] coordinates {(14,-2.5) (14.5,-3) (15,-2.5)};
    \draw [thick, black, dashed] plot [smooth, tension = 0.5] coordinates {(9.5,2.5) (9.8,2) (10.5,1.5) (11,-2.5)};
    \draw [thick, black, dashed] plot [smooth, tension = 0.5] coordinates {(10.5,2.5) (10.5,2.2) (11,1.5) (12,-1) (12.5,-2.5)};
    \draw [thick, purple, dashed] plot [smooth, tension = 0.5] coordinates {(11,2.5) (10.8,2.2) (10, 1.7) (9.5,1)};
    \draw [thick, purple, dashed] plot [smooth, tension = 0.5] coordinates {(12,2.5) (12.2,-1) (13.3,-1.5) (14,-2.5)};
    \draw [thick, purple, dashed] plot [smooth, tension = 0.5] coordinates {(12.5,2.5) (12.3,2.2) (10.7,1.8) (10.5,1)};
    \draw [thick, purple, dashed] plot [smooth, tension = 0.5] coordinates {(13.5,2.5) (13.7,-1) (14.8,-1.7) (15,-2.5)};
    \draw [thick, black, dashed] plot [smooth, tension = 0.5] coordinates {(14,-1) (13.8,-1.4) (12.2,-1.6) (12,-2.5)};
    \draw [thick, black, dashed] plot [smooth, tension = 0.5] coordinates {(15,-1) (14.8,-1.5) (13.7,-1.7) (13.5,-2.5)};
\end{tikzpicture}
    \caption{Weakly reducing $\S$ in Case 2}
    \label{fig:4BrL2.1}
\end{figure}
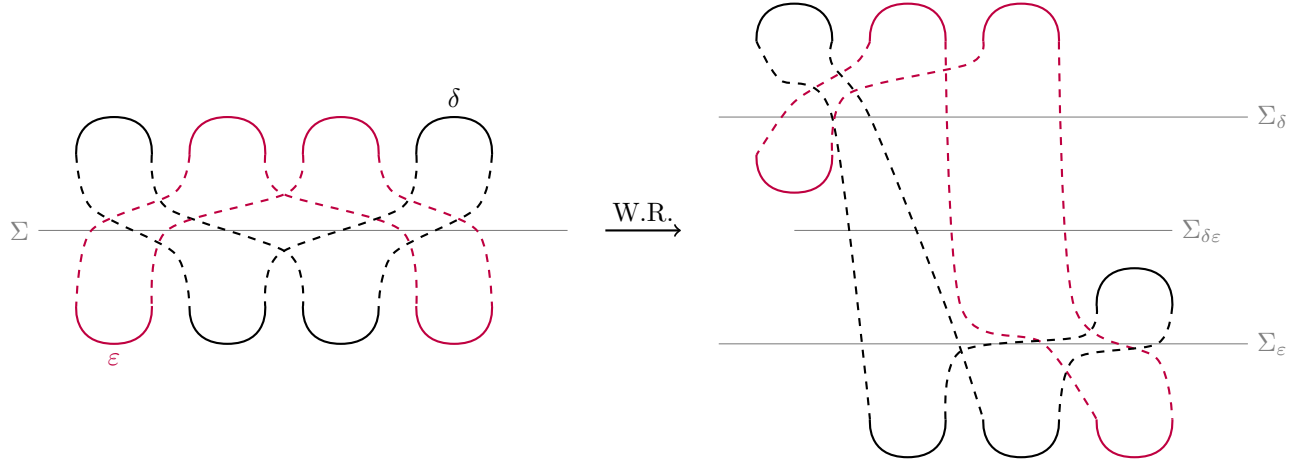

\

As in Case 1, we appeal to \cite{HS} and see that since $U$ is an unlink, $\overline{S^3 \bs\nu (U)}$ contains no essential surfaces other than spheres disjoint from the boundary, and so one of $\S_\d$ or $\S_\e$ must be weakly-reducible. We will assume $\S_\d$. Again, since the surgery from $\S$ to $\S_\d$ occurs away from $\e^*$, and by Lemma \ref{lemma:vert}, $\e^*$ appears in $\S_\d$ as the shadow of the unique bridge disk for $\e$. Thus, there must be a shadow $\d'*$ for the opposite tangle such that $\d'*$ is disjoint from $\e^*$.

Further, $\d'*$ cannot be a shadow for the same component of $U$ as $\e^*$, since they would then intersect at an endpoint. So $\d'^*$ is the shadow of a bridge disk for an arc $\d'$ in the component not containing $\e$. We weakly reduce along the arcs $\e$ and $\d'$, obtaining new level surfaces $\S_{\d\d'}$ and $\S_{\d\d'\e}$.

Notice now that $\S_{\d\d'\e}$ is a splitting surface for the (unknotted) component of $U$ containing $\e$, and so above $\S_{\d\d'\e}$, we have a situation identical to that in Case 1. We apply the same argument, and see that there is a system of shadows for $T_1$ containing $\d^*$ and intersecting $\e^*$ only at its endpoints. 

\begin{figure}[H]
    \centering

\begin{tikzpicture}
    \draw [black!50] (0,1.5) -- (7,1.5) node at (0,1.5) [anchor = east] {$\S_\d$};
    \draw [black!50] (1,0) -- (6,0);
    \draw [black!50] (0,-1.5) -- (7,-1.5);
    \draw [thick, black] plot [smooth, tension=2] coordinates {(0.5,2.5) (1,3) (1.5,2.5)} node at (1,3) [anchor = south] {$\d'$};
    \draw [thick, purple] plot [smooth, tension=2] coordinates {(2,2.5) (2.5,3) (3,2.5)};
    \draw [thick, purple] plot [smooth, tension=2] coordinates {(3.5,2.5) (4,3) (4.5,2.5)};
    \draw [thick, black] plot [smooth, tension=2] coordinates {(5,-1) (5.5,-0.5) (6,-1)};
    \draw [thick, purple] plot [smooth, tension=2] coordinates {(0.5,1) (1,0.5) (1.5,1)} node at (1,0.5) [anchor = north] {$\e$}; 
    \draw [thick, black] plot [smooth, tension=2] coordinates {(2,-2.5) (2.5,-3) (3,-2.5)};
    \draw [thick, black] plot [smooth, tension=2] coordinates {(3.5,-2.5) (4,-3) (4.5,-2.5)};
    \draw [thick, purple] plot [smooth, tension=2] coordinates {(5,-2.5) (5.5,-3) (6,-2.5)};
    \draw [thick, black, dashed] plot [smooth, tension = 0.5] coordinates {(0.5,2.5) (0.8,2) (1.5,1.5) (2,-2.5)};
    \draw [thick, black, dashed] plot [smooth, tension = 0.5] coordinates {(1.5,2.5) (1.5,2.2) (2,1.5) (3,-1) (3.5,-2.5)};
    \draw [thick, purple, dashed] plot [smooth, tension = 0.5] coordinates {(2,2.5) (1.8,2.2) (1, 1.7) (0.5,1)};
    \draw [thick, purple, dashed] plot [smooth, tension = 0.5] coordinates {(3,2.5) (3.2,-1) (4.3,-1.5) (5,-2.5)};
    \draw [thick, purple, dashed] plot [smooth, tension = 0.5] coordinates {(3.5,2.5) (3.3,2.2) (1.7,1.8) (1.5,1)};
    \draw [thick, purple, dashed] plot [smooth, tension = 0.5] coordinates {(4.5,2.5) (4.7,-1) (5.8,-1.7) (6,-2.5)};
    \draw [thick, black, dashed] plot [smooth, tension = 0.5] coordinates {(5,-1) (4.8,-1.4) (3.2,-1.6) (3,-2.5)};
    \draw [thick, black, dashed] plot [smooth, tension = 0.5] coordinates {(6,-1) (5.8,-1.5) (4.7,-1.7) (4.5,-2.5)};
    \draw [->] [thick, black] (7.5,0) -- (8.5,0) node at (8,0) [anchor = south] {W.R.};
    \draw [black!50] (9,3) -- (16,3) node at (16,3) [anchor = west] {$\S_{\d\d'}$};
    \draw [black!50] (10,1.5) -- (14,1.5) node at (14,1.5) [anchor = west] {$\S_{\d\d'\e}$};
    \draw [black!50] (9.5,0) -- (15,0);
    \draw [black!50] (9,-1.5) -- (16,-1.5);
    \draw [thick, black] plot [smooth, tension=2] coordinates {(9.5,0.5) (10,1) (10.5,0.5)};
    \draw [thick, purple] plot [smooth, tension=2] coordinates {(11,3.5) (11.5,4) (12,3.5)};
    \draw [thick, purple] plot [smooth, tension=2] coordinates {(12.5,3.5) (13,4) (13.5,3.5)};
    \draw [thick, black] plot [smooth, tension=2] coordinates {(14,-1) (14.5,-0.5) (15,-1)};
    \draw [thick, purple] plot [smooth, tension=2] coordinates {(9.5,2.5) (10,2) (10.5,2.5)} node at (10,2) [anchor = north] {$\e$}; 
    \draw [thick, black] plot [smooth, tension=2] coordinates {(11,-2.5) (11.5,-3) (12,-2.5)};
    \draw [thick, black] plot [smooth, tension=2] coordinates {(12.5,-2.5) (13,-3) (13.5,-2.5)};
    \draw [thick, purple] plot [smooth, tension=2] coordinates {(14,-2.5) (14.5,-3) (15,-2.5)};
    \draw [thick, black, dashed] plot [smooth, tension = 0.5] coordinates {(9.5,0.5) (9.8,0) (10.5,-0.5) (11,-2.5)};
    \draw [thick, black, dashed] plot [smooth, tension = 0.5] coordinates {(10.5,0.5) (10.5,0.2) (11,-0.5) (12,-1) (12.5,-2.5)};
    \draw [thick, purple, dashed] plot [smooth, tension = 0.5] coordinates {(11,3.5) (10.8,3.2) (10,3) (9.5,2.5)};
    \draw [thick, purple, dashed] plot [smooth, tension = 0.5] coordinates {(12,3.5) (12.2,-1) (13.3,-1.5) (14,-2.5)};
    \draw [thick, purple, dashed] plot [smooth, tension = 0.5] coordinates {(12.5,3.5) (12.3,3.2) (10.7,2.8) (10.5,2.5)};
    \draw [thick, purple, dashed] plot [smooth, tension = 0.5] coordinates {(13.5,3.5) (13.7,-1) (14.8,-1.7) (15,-2.5)};
    \draw [thick, black, dashed] plot [smooth, tension = 0.5] coordinates {(14,-1) (13.8,-1.4) (12.2,-1.6) (12,-2.5)};
    \draw [thick, black, dashed] plot [smooth, tension = 0.5] coordinates {(15,-1) (14.8,-1.5) (13.7,-1.7) (13.5,-2.5)};
\end{tikzpicture}

    \caption{Weakly reducing $\S_\d$ in Case 2}
    \label{fig:4BrL2.2}
\end{figure}

\

\underline{Case 3:} Suppose $c(U) = 2$ with $b_1 = 3$ and $b_2 = 1$. Suppose further that the weak-reducing arcs $\d$ and $\e$ are are part of the same component of $U$. They must both come from the 3-bridge component, as shadows for a 1-bridge splitting of $S^1$ will meet at both endpoints. We weakly reduce, obtaining level surfaces $\S_\d$, $\S_\e$, and $\S_{\d\e}$.

\

\begin{figure}[H]
    \centering

\begin{tikzpicture}
    \draw [black!50] (0,0) -- (7,0) node at (0,0) [anchor = east] {$\S$};
    \draw [thick, black] plot [smooth, tension=2] coordinates {(0.5,1) (1,1.5) (1.5,1)};
    \draw [thick, purple] plot [smooth, tension=2] coordinates {(2,1) (2.5,1.5) (3,1)};
    \draw [thick, black] plot [smooth, tension=2] coordinates {(3.5,1) (4,1.5) (4.5,1)};
    \draw [thick, black] plot [smooth, tension=2] coordinates {(5,1) (5.5,1.5) (6,1)} node at (5.5,1.5) [anchor = south] {$\d$};
    \draw [thick, black] plot [smooth, tension=2] coordinates {(0.5,-1) (1,-1.5) (1.5,-1)} node at (1,-1.5) [anchor = north] {$\e$}; 
    \draw [thick, black] plot [smooth, tension=2] coordinates {(2,-1) (2.5,-1.5) (3,-1)};
    \draw [thick, purple] plot [smooth, tension=2] coordinates {(3.5,-1) (4,-1.5) (4.5,-1)};
    \draw [thick, black] plot [smooth, tension=2] coordinates {(5,-1) (5.5,-1.5) (6,-1)};
    \draw [thick, black, dashed] plot [smooth, tension = 0.5] coordinates {(0.5,1) (0.5,-1)};
    \draw [thick, black, dashed] plot [smooth, tension = 0.5] coordinates {(1.5,1) (1.6,0.5) (1.9,-0.5) (2,-1)};
    \draw [thick, purple, dashed] plot [smooth, tension = 0.5] coordinates {(2,1) (2.2,0.5) (3.3, -0.5) (3.5,-1)};
    \draw [thick, purple, dashed] plot [smooth, tension = 0.5] coordinates {(3,1) (3.2,0.5) (4.3,-0.5) (4.5,-1)};
    \draw [thick, black, dashed] plot [smooth, tension = 0.5] coordinates {(3.5,1) (3.3,0.5) (1.7,-0.5) (1.5,-1)};
    \draw [thick, black, dashed] plot [smooth, tension = 0.5] coordinates {(4.5,1) (4.6,0.5) (4.9,-0.5) (5,-1)};
    \draw [thick, black, dashed] plot [smooth, tension = 0.5] coordinates {(5,1) (4.8,0.5) (3.2,-0.5) (3,-1)};
    \draw [thick, black, dashed] plot [smooth, tension = 0.5] coordinates {(6,1) (6,-1)};
    \draw [->] [thick, black] (7.5,0) -- (8.5,0) node at (8,0) [anchor = south] {W.R.};
    \draw [black!50] (9,1.5) -- (16,1.5) node at (16,1.5) [anchor = west] {$\S_\d$};
    \draw [black!50] (10,0) -- (15,0) node at (15,0) [anchor = west] {$\S_{\d\e}$};
    \draw [black!50] (9,-1.5) -- (16,-1.5) node at (16,-1.5) [anchor = west] {$\S_\e$};
    \draw [thick, black] plot [smooth, tension=2] coordinates {(9.5,2) (10,2.5) (10.5,2)};
    \draw [thick, purple] plot [smooth, tension=2] coordinates {(11,2) (11.5,2.5) (12,2)};
    \draw [thick, black] plot [smooth, tension=2] coordinates {(12.5,2) (13,2.5) (13.5,2)};
    \draw [thick, black] plot [smooth, tension=2] coordinates {(14,-1) (14.5,-0.5) (15,-1)};
    \draw [thick, black] plot [smooth, tension=2] coordinates {(9.5,1) (10,0.5) (10.5,1)}; 
    \draw [thick, black] plot [smooth, tension=2] coordinates {(11,-2) (11.5,-2.5) (12,-2)};
    \draw [thick, purple] plot [smooth, tension=2] coordinates {(12.5,-2) (13,-2.5) (13.5,-2)};
    \draw [thick, black] plot [smooth, tension=2] coordinates {(14,-2) (14.5,-2.5) (15,-2)};
    \draw [thick, black, dashed] plot [smooth, tension = 0.5] coordinates {(9.5,2) (9.5,1)};
    \draw [thick, black, dashed] plot [smooth, tension = 0.5] coordinates {(10.5,2) (10.7,0.5) (10.9,-0.5) (11,-2)};
    \draw [thick, black, dashed] plot [smooth, tension = 0.5] coordinates {(12.5,2) (12.3,1.6) (10.7, 1.4) (10.5,1)};
    \draw [thick, purple, dashed] plot [smooth, tension = 0.5] coordinates {(11,2) (11.2,1) (12.3,0) (12.5,-2)};
    \draw [thick, purple, dashed] plot [smooth, tension = 0.5] coordinates {(12,2) (12.2,1) (13.3,0) (13.5,-2)};
    \draw [thick, black, dashed] plot [smooth, tension = 0.5] coordinates {(13.5,2) (13.6,0.5) (13.8,-0.5) (14,-2)};
    \draw [thick, black, dashed] plot [smooth, tension = 0.5] coordinates {(14,-1) (13.8,-1.4) (12.2,-1.6) (12,-2)};
    \draw [thick, black, dashed] plot [smooth, tension = 0.5] coordinates {(15,-1) (15,-2)};
\end{tikzpicture}
    \caption{Weakly reducing $\S$ in Case 3}
    \label{fig:4BrL3.1}
\end{figure}
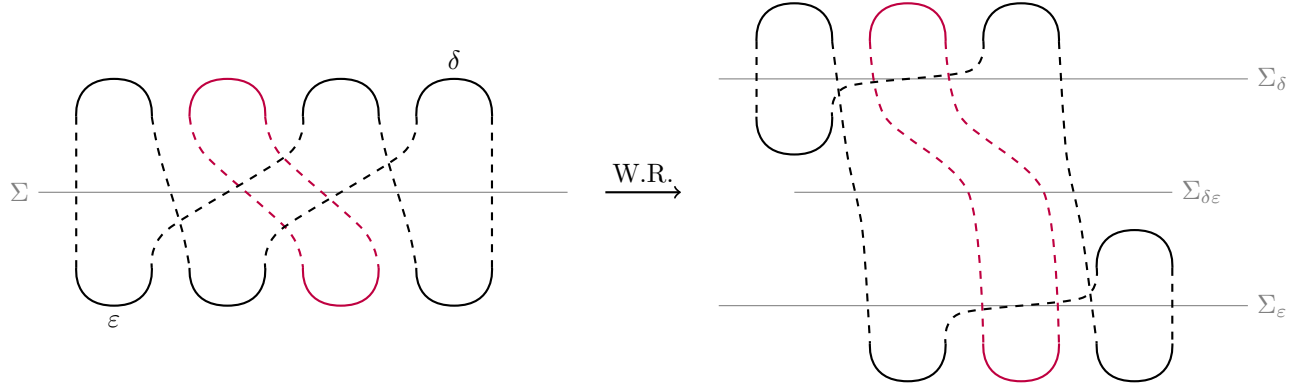

\

As in the previous cases, we must have one of the thick surfaces $\S_\d$ or $\S_\e$ weakly-reducible since $U$ is an unlink. We suppose $\S_\d$. In $\S_\d$, the tangle opposite $\e$ consists of two arcs each sharing an endpoint with $\e$ (they come from the same component of $U$) and one arc from the other component. The weak-reducing pair $\d'^*$ for $\e^*$ must be this last arc, as $\e^*$ intersects each of the other two shadows from the opposite tangle at an endpoint.

We weakly reduce along $\d'$ (the arc for which $\d'^*$ is a shadow) and $\e$, as in the previous cases, and obtain once more a surface $\S_{\d\d'\e}$ above which we have a punctured unknot split over a sphere with four marked points. We apply the argument from Case 1 and see that $\e$ is primitive with respect to $T_1$, as exhibited by a shadow diagram including $\d^*$ and $\e^*$.

\

\begin{figure}[H]
    \centering

\begin{tikzpicture}
    \draw [black!50] (0,1.5) -- (7,1.5) node at (0,1.5) [anchor = east] {$\S_\d$};
    \draw [black!50] (1,0) -- (6,0);
    \draw [black!50] (0,-1.5) -- (7,-1.5);
    \draw [thick, black] plot [smooth, tension=2] coordinates {(0.5,2) (1,2.5) (1.5,2)};
    \draw [thick, purple] plot [smooth, tension=2] coordinates {(2,2) (2.5,2.5) (3,2)} node at (2.5,2.5) [anchor = south] {$\d'$};
    \draw [thick, black] plot [smooth, tension=2] coordinates {(3.5,2) (4,2.5) (4.5,2)};
    \draw [thick, black] plot [smooth, tension=2] coordinates {(5,-1) (5.5,-0.5) (6,-1)};
    \draw [thick, black] plot [smooth, tension=2] coordinates {(0.5,1) (1,0.5) (1.5,1)} node at (1,0.5) [anchor = north] {$\e$}; 
    \draw [thick, black] plot [smooth, tension=2] coordinates {(2,-2) (2.5,-2.5) (3,-2)};
    \draw [thick, purple] plot [smooth, tension=2] coordinates {(3.5,-2) (4,-2.5) (4.5,-2)};
    \draw [thick, black] plot [smooth, tension=2] coordinates {(5,-2) (5.5,-2.5) (6,-2)};
    \draw [thick, black, dashed] plot [smooth, tension = 0.5] coordinates {(0.5,2) (0.5,1)};
    \draw [thick, black, dashed] plot [smooth, tension = 0.5] coordinates {(1.5,2) (1.7,0.5) (1.9,-0.5) (2,-2)};
    \draw [thick, black, dashed] plot [smooth, tension = 0.5] coordinates {(3.5,2) (3.3,1.6) (1.7, 1.4) (1.5,1)};
    \draw [thick, purple, dashed] plot [smooth, tension = 0.5] coordinates {(2,2) (2.2,1) (3.3,0) (3.5,-2)};
    \draw [thick, purple, dashed] plot [smooth, tension = 0.5] coordinates {(3,2) (3.2,1) (4.3,0) (4.5,-2)};
    \draw [thick, black, dashed] plot [smooth, tension = 0.5] coordinates {(4.5,2) (4.6,0.5) (4.8,-0.5) (5,-2)};
    \draw [thick, black, dashed] plot [smooth, tension = 0.5] coordinates {(5,-1) (4.8,-1.4) (3.2,-1.6) (3,-2)};
    \draw [thick, black, dashed] plot [smooth, tension = 0.5] coordinates {(6,-1) (6,-2)};
    \draw [->] [thick, black] (7.5,0) -- (8.5,0) node at (8,0) [anchor = south] {W.R.};
    \draw [black!50] (9,3) -- (16,3) node at (16,3) [anchor = west] {$\S_{\d\d'}$};
    \draw [black!50] (10,1.5) -- (14,1.5) node at (14,1.5) [anchor = west] {$\S_{\d\d'\e}$};
    \draw [black!50] (10,0) -- (15,0);
    \draw [black!50] (9,-1.5) -- (16,-1.5);
    \draw [thick, black] plot [smooth, tension=2] coordinates {(9.5,3.5) (10,4) (10.5,3.5)};
    \draw [thick, purple] plot [smooth, tension=2] coordinates {(11,0.5) (11.5,1) (12,0.5)};
    \draw [thick, black] plot [smooth, tension=2] coordinates {(12.5,3.5) (13,4) (13.5,3.5)};
    \draw [thick, black] plot [smooth, tension=2] coordinates {(14,-1) (14.5,-0.5) (15,-1)};
    \draw [thick, black] plot [smooth, tension=2] coordinates {(9.5,2.5) (10,2) (10.5,2.5)} node at (10,2) [anchor = north] {$\e$}; 
    \draw [thick, black] plot [smooth, tension=2] coordinates {(11,-2) (11.5,-2.5) (12,-2)};
    \draw [thick, purple] plot [smooth, tension=2] coordinates {(12.5,-2) (13,-2.5) (13.5,-2)};
    \draw [thick, black] plot [smooth, tension=2] coordinates {(14,-2) (14.5,-2.5) (15,-2)};
    \draw [thick, black, dashed] plot [smooth, tension = 0.5] coordinates {(9.5,3.5) (9.5,2.5)};
    \draw [thick, black, dashed] plot [smooth, tension = 0.5] coordinates {(10.5,3.5) (10.55,3.1) (10.7,2.5) (10.9,-0.5) (11,-2)};
    \draw [thick, black, dashed] plot [smooth, tension = 0.5] coordinates {(12.5,3.5) (12.3,3.1) (10.7, 2.9) (10.5,2.5)};
    \draw [thick, purple, dashed] plot [smooth, tension = 0.5] coordinates {(11,0.5) (11.2,-0.5) (12.3,-1.5) (12.5,-2)};
    \draw [thick, purple, dashed] plot [smooth, tension = 0.5] coordinates {(12,0.5) (12.2,-0.5) (13.3,-1.5) (13.5,-2)};
    \draw [thick, black, dashed] plot [smooth, tension = 0.5] coordinates {(13.5,3.5) (13.6,0.5) (13.8,-0.5) (14,-2)};
    \draw [thick, black, dashed] plot [smooth, tension = 0.5] coordinates {(14,-1) (13.8,-1.4) (12.2,-1.6) (12,-2)};
    \draw [thick, black, dashed] plot [smooth, tension = 0.5] coordinates {(15,-1) (15,-2)};
\end{tikzpicture}

\caption{Weakly reducing $\S_\d$ in Case 3}
    \label{fig:4BrL3.2}
\end{figure}

\

\underline{Case 4:} Suppose once more that $c(U) = 2$, with $b_1 = 3$ and $b_2 = 1$. Let the reducing arcs $\d$ and $\e$ lie in different components of $U$, supposing $\d$ lies in the one-bridge component. We weakly reduce as usual, noting our level surfaces $\S_\d$, $\S_\e$, and $\S_{\d\e}$

\

\begin{figure}[H]
    \centering

\begin{tikzpicture}
    \draw [black!50] (0,0) -- (7,0) node at (0,0) [anchor = east] {$\S$};
    \draw [thick, black] plot [smooth, tension=2] coordinates {(0.5,1) (1,1.5) (1.5,1)};
    \draw [thick, purple] plot [smooth, tension=2] coordinates {(2,1) (2.5,1.5) (3,1)} node at (2.5,1.5) [anchor = south] {$\d$};
    \draw [thick, black] plot [smooth, tension=2] coordinates {(3.5,1) (4,1.5) (4.5,1)};
    \draw [thick, black] plot [smooth, tension=2] coordinates {(5,1) (5.5,1.5) (6,1)};
    \draw [thick, black] plot [smooth, tension=2] coordinates {(0.5,-1) (1,-1.5) (1.5,-1)}; 
    \draw [thick, black] plot [smooth, tension=2] coordinates {(2,-1) (2.5,-1.5) (3,-1)};
    \draw [thick, purple] plot [smooth, tension=2] coordinates {(3.5,-1) (4,-1.5) (4.5,-1)};
    \draw [thick, black] plot [smooth, tension=2] coordinates {(5,-1) (5.5,-1.5) (6,-1)} node at (5.5,-1.5) [anchor = north] {$\e$};
    \draw [thick, black, dashed] plot [smooth, tension = 0.5] coordinates {(0.5,1) (0.5,-1)};
    \draw [thick, black, dashed] plot [smooth, tension = 0.5] coordinates {(1.5,1) (1.6,0.5) (1.9,-0.5) (2,-1)};
    \draw [thick, purple, dashed] plot [smooth, tension = 0.5] coordinates {(2,1) (2.2,0.5) (3.3, -0.5) (3.5,-1)};
    \draw [thick, purple, dashed] plot [smooth, tension = 0.5] coordinates {(3,1) (3.2,0.5) (4.3,-0.5) (4.5,-1)};
    \draw [thick, black, dashed] plot [smooth, tension = 0.5] coordinates {(3.5,1) (3.3,0.5) (1.7,-0.5) (1.5,-1)};
    \draw [thick, black, dashed] plot [smooth, tension = 0.5] coordinates {(4.5,1) (4.6,0.5) (4.9,-0.5) (5,-1)};
    \draw [thick, black, dashed] plot [smooth, tension = 0.5] coordinates {(5,1) (4.8,0.5) (3.2,-0.5) (3,-1)};
    \draw [thick, black, dashed] plot [smooth, tension = 0.5] coordinates {(6,1) (6,-1)};
    \draw [->] [thick, black] (7.5,0) -- (8.5,0) node at (8,0) [anchor = south] {W.R.};
    \draw [black!50] (9,1.5) -- (16,1.5) node at (16,1.5) [anchor = west] {$\S_\d$};
    \draw [black!50] (9.25,0) -- (15,0) node at (15,0) [anchor = west] {$\S_{\d\e}$};
    \draw [black!50] (9,-1.5) -- (16,-1.5) node at (16,-1.5) [anchor = west] {$\S_\e$};
    \draw [thick, black] plot [smooth, tension=2] coordinates {(9.5,2) (10,2.5) (10.5,2)};
    \draw [thick, purple] plot [smooth, tension=2] coordinates {(11,-1) (11.5,-0.5) (12,-1)};
    \draw [thick, black] plot [smooth, tension=2] coordinates {(12.5,2) (13,2.5) (13.5,2)};
    \draw [thick, black] plot [smooth, tension=2] coordinates {(14,2) (14.5,2.5) (15,2)};
    \draw [thick, black] plot [smooth, tension=2] coordinates {(9.5,-2) (10,-2.5) (10.5,-2)}; 
    \draw [thick, black] plot [smooth, tension=2] coordinates {(11,-2) (11.5,-2.5) (12,-2)};
    \draw [thick, purple] plot [smooth, tension=2] coordinates {(12.5,-2) (13,-2.5) (13.5,-2)};
    \draw [thick, black] plot [smooth, tension=2] coordinates {(14,1) (14.5,0.5) (15,1)};
    \draw [thick, black, dashed] plot [smooth, tension = 0.5] coordinates {(9.5,2) (9.5,-2)};
    \draw [thick, black, dashed] plot [smooth, tension = 0.5] coordinates {(10.5,2) (10.7,0.5) (10.9,-0.5) (11,-2)};
    \draw [thick, black, dashed] plot [smooth, tension = 0.5] coordinates {(12.5,2) (12.3,1) (10.7, -1) (10.5,-2)};
    \draw [thick, purple, dashed] plot [smooth, tension = 0.5] coordinates {(11,-1) (11.2,-1.3) (12.3,-1.7) (12.5,-2)};
    \draw [thick, purple, dashed] plot [smooth, tension = 0.5] coordinates {(12,-1) (12.2,-1.3) (13.3,-1.7) (13.5,-2)};
    \draw [thick, black, dashed] plot [smooth, tension = 0.5] coordinates {(13.5,2) (13.6,1.7) (13.9,1.3) (14,1)};
    \draw [thick, black, dashed] plot [smooth, tension = 0.5] coordinates {(14,2) (13.8,1) (12.2,-1) (12,-2)};
    \draw [thick, black, dashed] plot [smooth, tension = 0.5] coordinates {(15,2) (15,1)};
\end{tikzpicture}

    \caption{Weakly reducing $\S$ in Case 4}
    \label{fig:4BrL4}
\end{figure}
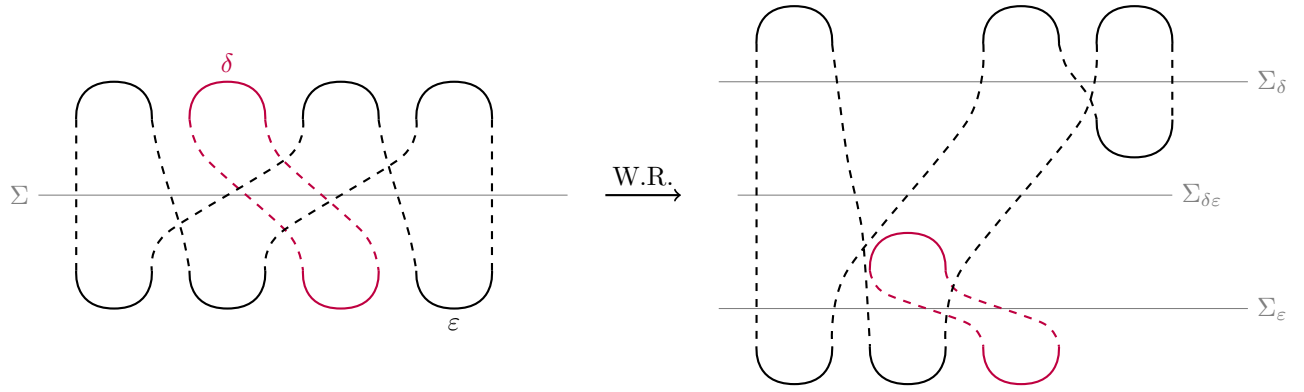

\

As in the previous cases, one of the two thick surfaces must be weakly-reducible, since $U$'s complement does not contain any essential punctured surfaces. This situation lacks the symmetry of the others, though, so we must consider separately weak reductions of $\S_\d$ and of $\S_\e$.

\underline{Case 4a:} Suppose $\S_\d$ is weakly-reducible. Let $\d'^*$ be a shadow in $\S_\d$ for the tangle opposite $\e$ such that $\d'^* \cap \e^* = \0$. We weakly reduce again and obtain a surface $\S_{\d\d'}$ splitting a punctured unknot on one side of $\S_{\d\d'\e}$. We refer to Case 1, and $\e$ is primitive with respect to $T_1$, as exhibited by a shadow diagram including $\d^*$ and $\e^*$.

\begin{figure}[H]
    \centering

\begin{tikzpicture}
    \draw [black!50] (0,1.5) -- (7,1.5) node at (0,1.5) [anchor = east] {$\S_\d$};
    \draw [black!50] (0.25,0) -- (6,0);
    \draw [black!50] (0,-1.5) -- (7,-1.5);
    \draw [thick, black] plot [smooth, tension=2] coordinates {(0.5,2) (1,2.5) (1.5,2)} node at (1,2.5) [anchor = south] {$\d'$};
    \draw [thick, purple] plot [smooth, tension=2] coordinates {(2,-1) (2.5,-0.5) (3,-1)};
    \draw [thick, black] plot [smooth, tension=2] coordinates {(3.5,2) (4,2.5) (4.5,2)};
    \draw [thick, black] plot [smooth, tension=2] coordinates {(5,2) (5.5,2.5) (6,2)};
    \draw [thick, black] plot [smooth, tension=2] coordinates {(0.5,-2) (1,-2.5) (1.5,-2)}; 
    \draw [thick, black] plot [smooth, tension=2] coordinates {(2,-2) (2.5,-2.5) (3,-2)};
    \draw [thick, purple] plot [smooth, tension=2] coordinates {(3.5,-2) (4,-2.5) (4.5,-2)};
    \draw [thick, black] plot [smooth, tension=2] coordinates {(5,1) (5.5,0.5) (6,1)}  node at (5.5,0.5) [anchor = north] {$\e$};
    \draw [thick, black, dashed] plot [smooth, tension = 0.5] coordinates {(0.5,2) (0.5,-2)};
    \draw [thick, black, dashed] plot [smooth, tension = 0.5] coordinates {(1.5,2) (1.7,0.5) (1.9,-0.5) (2,-2)};
    \draw [thick, black, dashed] plot [smooth, tension = 0.5] coordinates {(3.5,2) (3.3,1) (1.7, -1) (1.5,-2)};
    \draw [thick, purple, dashed] plot [smooth, tension = 0.5] coordinates {(2,-1) (2.2,-1.3) (3.3,-1.7) (3.5,-2)};
    \draw [thick, purple, dashed] plot [smooth, tension = 0.5] coordinates {(3,-1) (3.2,-1.3) (4.3,-1.7) (4.5,-2)};
    \draw [thick, black, dashed] plot [smooth, tension = 0.5] coordinates {(4.5,2) (4.6,1.7) (4.9,1.3) (5,1)};
    \draw [thick, black, dashed] plot [smooth, tension = 0.5] coordinates {(5,2) (4.8,1) (3.2,-1) (3,-2)};
    \draw [thick, black, dashed] plot [smooth, tension = 0.5] coordinates {(6,2) (6,1)};
    \draw [->] [thick, black] (7.5,0) -- (8.5,0) node at (8,0) [anchor = south] {W.R.};
    \draw [black!50] (9,3) -- (16,3)  node at (16,3) [anchor = west] {$\S_{\d\d'}$};
    \draw [black!50] (11,1.5) -- (14,1.5) node at (14,1.5) [anchor = west] {$\S_{\d\d'\e}$};
    \draw [black!50] (9,0) -- (14,0);
    \draw [black!50] (9,-1.5) -- (16,-1.5);
    \draw [thick, black] plot [smooth, tension=2] coordinates {(9.5,0.5) (10,1) (10.5,0.5)};
    \draw [thick, purple] plot [smooth, tension=2] coordinates {(11,-1) (11.5,-0.5) (12,-1)};
    \draw [thick, black] plot [smooth, tension=2] coordinates {(12.5,3.5) (13,4) (13.5,3.5)};
    \draw [thick, black] plot [smooth, tension=2] coordinates {(14,3.5) (14.5,4) (15,3.5)};
    \draw [thick, black] plot [smooth, tension=2] coordinates {(9.5,-2) (10,-2.5) (10.5,-2)}; 
    \draw [thick, black] plot [smooth, tension=2] coordinates {(11,-2) (11.5,-2.5) (12,-2)};
    \draw [thick, purple] plot [smooth, tension=2] coordinates {(12.5,-2) (13,-2.5) (13.5,-2)};
    \draw [thick, black] plot [smooth, tension=2] coordinates {(14,2.5) (14.5,2) (15,2.5)} node at (14.75,2.1) [anchor = north west] {$\e$};
    \draw [thick, black, dashed] plot [smooth, tension = 0.5] coordinates {(9.5,0.5) (9.5,-2)};
    \draw [thick, black, dashed] plot [smooth, tension = 0.5] coordinates {(10.5,0.5) (10.7,0.0) (10.9,-0.5) (11,-2)};
    \draw [thick, black, dashed] plot [smooth, tension = 0.5] coordinates {(12.5,3.5) (12.3,1) (10.7, -1) (10.5,-2)};
    \draw [thick, purple, dashed] plot [smooth, tension = 0.5] coordinates {(11,-1) (11.2,-1.3) (12.3,-1.7) (12.5,-2)};
    \draw [thick, purple, dashed] plot [smooth, tension = 0.5] coordinates {(12,-1) (12.2,-1.3) (13.3,-1.7) (13.5,-2)};
    \draw [thick, black, dashed] plot [smooth, tension = 0.5] coordinates {(13.5,3.5) (13.6,3.2) (13.9,2.8) (14,2.5)};
    \draw [thick, black, dashed] plot [smooth, tension = 0.5] coordinates {(14,3.5) (13.8,2.5) (12.2,-1) (12,-2)};
    \draw [thick, black, dashed] plot [smooth, tension = 0.5] coordinates {(15,3.5) (15,2.5)};
\end{tikzpicture}

    \caption{Weakly reducing $\S_\d$ in Case 4a}
    \label{fig:4BrL4a}
\end{figure}

\underline{Case 4b:} Suppose $\S_\e$ is weakly-reducible. So $\d^*$, which, by the same argument as in Case 1, is the shadow of the unique bridge disk for $\d$ intersecting $\S_\e$, is disjoint from some shadow $\e'^* \subset \S_\e$ for the tangle opposite $\d$. We weakly reduce, and see that we now have a new situation in $\S_{\e\e'}$.

\

\begin{figure}[H]
    \centering

\begin{tikzpicture}
    \draw [black!50] (0,1.5) -- (7,1.5);
    \draw [black!50] (1,0) -- (6,0);
    \draw [black!50] (0,-1.5) -- (7,-1.5) node at (0,-1.5) [anchor = east] {$\S_\e$};
    \draw [thick, black] plot [smooth, tension=2] coordinates {(0.5,2) (1,2.5) (1.5,2)};
    \draw [thick, purple] plot [smooth, tension=2] coordinates {(2,-1) (2.5,-0.5) (3,-1)} node at (2.5,-0.5) [anchor = south west] {$\d$};
    \draw [thick, black] plot [smooth, tension=2] coordinates {(3.5,2) (4,2.5) (4.5,2)};
    \draw [thick, black] plot [smooth, tension=2] coordinates {(5,2) (5.5,2.5) (6,2)};
    \draw [thick, black] plot [smooth, tension=2] coordinates {(0.5,-2) (1,-2.5) (1.5,-2)} node at (1,-2.5) [anchor = north] {$\e'$}; 
    \draw [thick, black] plot [smooth, tension=2] coordinates {(2,-2) (2.5,-2.5) (3,-2)};
    \draw [thick, purple] plot [smooth, tension=2] coordinates {(3.5,-2) (4,-2.5) (4.5,-2)};
    \draw [thick, black] plot [smooth, tension=2] coordinates {(5,1) (5.5,0.5) (6,1)};
    \draw [thick, black, dashed] plot [smooth, tension = 0.5] coordinates {(0.5,2) (0.5,-2)};
    \draw [thick, black, dashed] plot [smooth, tension = 0.5] coordinates {(1.5,2) (1.7,0.5) (1.9,-0.5) (2,-2)};
    \draw [thick, black, dashed] plot [smooth, tension = 0.5] coordinates {(3.5,2) (3.3,1) (1.7, -1) (1.5,-2)};
    \draw [thick, purple, dashed] plot [smooth, tension = 0.5] coordinates {(2,-1) (2.2,-1.3) (3.3,-1.7) (3.5,-2)};
    \draw [thick, purple, dashed] plot [smooth, tension = 0.5] coordinates {(3,-1) (3.2,-1.3) (4.3,-1.7) (4.5,-2)};
    \draw [thick, black, dashed] plot [smooth, tension = 0.5] coordinates {(4.5,2) (4.6,1.7) (4.9,1.3) (5,1)};
    \draw [thick, black, dashed] plot [smooth, tension = 0.5] coordinates {(5,2) (4.8,1) (3.2,-1) (3,-2)};
    \draw [thick, black, dashed] plot [smooth, tension = 0.5] coordinates {(6,2) (6,1)};
    \draw [->] [thick, black] (7.5,0) -- (8.5,0) node at (8,0) [anchor = south] {W.R.};
    \draw [black!50] (9,1.5) -- (16,1.5);
    \draw [black!50] (9,0) -- (14,0);
    \draw [black!50] (10,-1.5) -- (13.5,-1.5) node at (13.5,-1.5) [anchor = west] {$\S_{\d\e\e'}$};
    \draw [black!50] (9,-3) -- (14,-3) node at (14,-3) [anchor = west] {$\S_{\e\e'}$};
    \draw [thick, black] plot [smooth, tension=2] coordinates {(9.5,2) (10,2.5) (10.5,2)};
    \draw [thick, purple] plot [smooth, tension=2] coordinates {(11,-2.5) (11.5,-2) (12,-2.5)} node at (11.5,-2) [anchor = south] {$\d$};
    \draw [thick, black] plot [smooth, tension=2] coordinates {(12.5,2) (13,2.5) (13.5,2)};
    \draw [thick, black] plot [smooth, tension=2] coordinates {(14,2) (14.5,2.5) (15,2)};
    \draw [thick, black] plot [smooth, tension=2] coordinates {(9.5,-0.5) (10,-1) (10.5,-0.5)}; 
    \draw [thick, black] plot [smooth, tension=2] coordinates {(11,-3.5) (11.5,-4) (12,-3.5)};
    \draw [thick, purple] plot [smooth, tension=2] coordinates {(12.5,-3.5) (13,-4) (13.5,-3.5)};
    \draw [thick, black] plot [smooth, tension=2] coordinates {(14,1) (14.5,0.5) (15,1)};
    \draw [thick, black, dashed] plot [smooth, tension = 0.5] coordinates {(9.5,2) (9.5,-0.5)};
    \draw [thick, black, dashed] plot [smooth, tension = 0.5] coordinates {(10.5,2) (10.6,0.5) (10.7 ,-0.5) (11,-3.5)};
    \draw [thick, black, dashed] plot [smooth, tension = 0.5] coordinates {(12.5,2) (12.3,1) (10.7, 0) (10.5,-0.5)};
    \draw [thick, purple, dashed] plot [smooth, tension = 0.5] coordinates {(11,-2.5) (11.2,-2.8) (12.3,-3.2) (12.5,-3.5)};
    \draw [thick, purple, dashed] plot [smooth, tension = 0.5] coordinates {(12,-2.5) (12.2,-2.8) (13.3,-3.2) (13.5,-3.5)};
    \draw [thick, black, dashed] plot [smooth, tension = 0.5] coordinates {(13.5,2) (13.6,1.7) (13.9,1.3) (14,1)};
    \draw [thick, black, dashed] plot [smooth, tension = 0.5] coordinates {(14,2) (13.8,1) (12.5,-1) (12,-3.5)};
    \draw [thick, black, dashed] plot [smooth, tension = 0.5] coordinates {(15,2) (15,1)};
\end{tikzpicture}

 \caption{Weakly reducing $\S_{\e}$ in Case 4b}
    \label{fig:4BrL4b}
\end{figure}

The shadows in $\S_{\e\e'}$ form the union of a 1-bridge unknot with a 1-bridge trivial arc. Also, $\S_{\e\e'}$ is a 4-punctured sphere. We normalize such that $\d^*$ is a $\frac10$ arc in $\S_{\e\e'}$ and notice that there is a unique arc $\e''^*$ in $\S_{\e\e'}$ such that $\d^* \cup \e''^*$ is a 1-bridge splitting of an unknot. This means that $\e''^*$ is a $\frac10$ arc in $\S_{\e\e'}$, and so $\d^*$ is parallel to some shadow $\e''^*$ from the opposite tangle. As in Case 1, any isotopy of $\e''^*$ in $\S_{\e\e'}$ can be realized by isotopy and possible shadow slides over $\e'^*$ in $S_\e$. Likewise, we can identify a shadow for $T_2$ isotopic to $e''^*$ in $\S$ after possible slides over $\e^*$ and $\e'^*$. Thus, $\d^*$ is isotopic to a shadow of some bridge disk for $T_2$ in a shadow diagram also including $\e^*$.

\end{proof}

Equipped with Lemma \ref{lemma:Classical}, we can conclude that weakly-reducible 4-bridge trisections are nice; namely, that we can realize forgetting $\d$ on the level of shadows and the result is, indeed, a bridge trisection.

\begin{lemma}
\label{lemma:4BrNice}
   Let $\mathcal{T} = (\S;T_\a,T_\b,T_\g)$ be a weakly-reducible 4-bridge trisection for a connected surface $K \subset S^4$ with weak reducing shadows $\d^*$ and $\e^*$, where $\e^*$ is a shadow for $\e_\a \in T_\a$ and $\e_\b \in T_\b$, and $\d^*$ is a shadow for $\d \in T_\g$. Then either $\mathcal{T}$ is perturbed, or else $\d^*$ is primitive with respect to each of $T_\a$ and $T_\b$ and the result $\mathcal{T}_\d$ of forgetting $\d$ is a bridge trisection.
\end{lemma}

\begin{proof}

Note first that if any of the pairs $T_\a \cup \overline{T_\b}$, $T_\b \cup \overline{T_\g}$, or $T_\g \cup \overline{T_\g}$ has more than two components, then by Theorem 5.2 of \cite{jmmzsolids}, $K$ is unknotted and $\mathcal{T}$ is perturbed. We will suppose that $c_i$, the number of components in each pairwise union of tangles, is at most 2 for all $i$.

Notice further that since $\d^* \cap \e^* = \0$, and since $\e^*$ is a shadow for each of $T_\a$ and $T_\b$, each of the pairs $T_\b \cup \overline{T_\g}$ and $T_\g \cup \overline{T_\a}$ is a weakly-reducible bridge splitting with reducing shadows $\d^*$ and $\e^*$. Further, each of these pairs is a 4-bridge splitting of an unknot or 2-component unlink, and so in each pair, either $\d^*$ is primitive with respect to the opposite tangle, $\e^*$ is primitive with respect to $T_\g$, or else one of $\e^*$ or $\d^*$ is a shadow for some bridge disk from the opposite tangle.

Suppose first that $\e^*$ co-bounds a bridge disk in $T_\g$. Then $\e^*$ is part of a 1-bridge component of the splitting. But this would mean that $K$ has two components, as $\e^*$ is now parallel into each of the three tangles. Since $K$ is connected, this cannot occur.

Suppose instead that $\e^*$ is primitive with respect to $T_\g$. Then an arc in each of $T_\a$ and $T_\b$ is primitive with respect to $T_\g$, and so $\mathcal{T}$ is perturbed by an application of Lemma 6.2 of \cite{mz}.

Suppose now that $\d^*$ co-bounds a bridge disk in $T_\a$. If $\d^*$ also co-bounds a bridge disk in $T_\b$, then we once again have a disconnected $K$, so we must have $\d^*$ primitive with respect to $T_\b$. Once again, we see that the splitting is perturbed. By symmetry, if $\d^*$ co-bounds a bridge disk in $T_\b$ then $\mathcal{T}$ is perturbed.

Thus, if $\mathcal{T}$ is not perturbed, then $\d^*$ must be primitive with respect to each of $T_\a$ and $T_\b$.

To see that $\mathcal{T}_\d$ is a bridge trisection, we examine the Heegaard triple $(\widetilde{\S};H_\a,H_\b,H_\g)$ consisting of the double branched cover of each of the tangles $f_\d(T_\a)$, $f_\d(T_\b)$, and $f_\d(T_\g)$. Each of $f_\d(T_\b) \cup \overline{f_\d(T_\g)}$ and $f_\d(T_\g) \cup \overline{f_\d(T_\a)}$ is a deperturbation of a 4-bridge splitting of an unlink, since forgetting $\d$ can be realized on the level of the pairwise union as as pushing $\d$ across $\S$. Since $\d$ is primitive with respect to each of the opposite tangles, this deperturbs the splitting in both cases.

Further, since we suppose that each $c_i$ is at most 2, we have that each of $H_\b \cup H_\g$ and $H_\g \cup H_\a$ is either $S^3$ (in the case that the relevant unlink is in fact an unknot) or $S^1 \x S^2$ (when the unlink has two components). The third pairing, $f_\d(T_\a) \cup \overline{f_\d(T_\b)}$, is the split union of a 1-bridge component formed by the arcs for which $\e^*$ is a shadow in each of $f_\d(T_\a)$ and $f_\d(T_\b)$ with a two-bridge link $L$. Thus, its double branched $H_\a \cup H_\b$ cover is $(S^1 \x S^2) \# L(p,q)$ for some lens space $L(p,q)$.

In \cite{a&r}, each of these situations is considered in the study of Heegaard triples. The case where one of $f_\d(T_\b) \cup \overline{f_\d(T_\g)}$ and $f_\d(T_\g) \cup \overline{f_\d(T_\a)}$ is an unknot and the other a two-component unlink is impossible, per Lemma 4.2 of \cite{a&r}. In the case where each of $f_\d(T_\b) \cup \overline{f_\d(T_\g)}$ and $f_\d(T_\g) \cup \overline{f_\d(T_\a)}$ has two components, Proposition 4.3 of \cite{a&r} tells us that there is a curve $c$ in $\widetilde{\S}$ which bounds a disk in each of $H_\a$, $H_\b$, and $H_\g$. Since $\widetilde{\S}$ has genus two, $c$ must be equivariant with respect to the hyperelliptic involution. Thus, there is some shadow $\zeta^* \cup \S$ covered by $c$ which is a shadow for some arc in each of $f_\d(T_\a)$, $f_\d(T_\b)$, and $f_\d(T_\g)$. Each of $f_\d(T_\b) \cup \overline{f_\d(T_\g)}$ and $f_\d(T_\g) \cup \overline{f_\d(T_\a)}$ consists of a two-component unlink, with a single 1-bridge component.

Notice that if $\e^*$ is a shadow for some $\e_\g \in f_\d(T_\g)$, then $\e_\g \in T_\g$, as forgetting $\d$ occurs away from other arcs of $T_\g$. In this case, $K$ would be disconnected, so $\zeta^*$ must be away from $\e$, i.e., a shadow for an arc of the 2-bridge link $L$ in $f_\d(T_\a) \cup \overline{f_\d(T_\b)}$. Since $\e_\a \cup \overline{\e_\b}$ is split from $L$, shadows for $L$ will appear as $\frac{p}{q}$-slope curves in the 4-punctured sphere $(\S_\d \bs \nu(\e^*)) \cup_{\del\nu(\e^*)} D^2$. Each of the arcs $\zeta_\a \in f_\d(T_\a)$ and $\zeta_\b \in f_\d(T_\b)$ for which $\zeta^*$ is a shadow forms a 1-bridge unknot when paired with $\zeta_\g \in f_\d(T_\g)$, since each of $f_\d(T_\b) \cup \overline{f_\d(T_\g)}$ and $f_\d(T_\g) \cup \overline{f_\d(T_\a)}$ is an unlink. We now see that (after possibly normalizing $\S_\d$), $\zeta^*$ has slope $\frac10$ or $\frac01$, and thus $f_\d(T_\a) \cup \overline{f_\d(T_\b)}$ is a 3-component unlink.

In the final case, where each of $f_\d(T_\b) \cup \overline{f_\d(T_\g)}$ and $f_\d(T_\g) \cup \overline{f_\d(T_\a)}$ is an unknot, we appeal to Theorem 4.8 of \cite{a&r}. Note that our conventions differ so that their $\a$ handlebody corresponds to our $\g$ handlebody; trisections and bridge-trisections are invariant under cyclic permutation of the pieces. We now have curves $c_1,c_2 \in \widetilde{\S}$ in the double branched cover of $\mathcal{T}_\d$ such that $c_1$ bounds a disk in each of $H_\a$ and $H_\b$, $c_2$ bounds a disk in $H_\g$, and $|c_1 \cap c_2| = 1$.

Once again, each of these curves is invariants under the hyperelliptic involution, since $\widetilde{\S}$ has genus two. Thus, there are arcs $\zeta_1^*$ and $\zeta_2^*$ in $\S_\d$ that lift to $c_1$ and $c_2$ respectively. Since $c_1$ is a compressing curve in each of $H_\a$ and $H_\b$, $\zeta_1^*$ is the shadow of a bridge disk for each of $f_\d(T_\a)$ and $f_\d(T_\b)$, say for arcs $\zeta_\a$ and $\zeta_\b$. Likewise, $\zeta_2^*$ is a shadow for some arc $\zeta_\g \in f_\d(T_\g)$. Since $c_1$ and $c_2$ intersect at a single point, we may isotope them until this point is part of the branching set for the cover, meaning that $\zeta_1^*$ and $\zeta_2^*$ meet at an endpoint. Thus, $\zeta_1^*$ is primitive with respect to $f_\d(T_\g)$.

Suppose that $\zeta_\a = \e_\a$, and thus that $\zeta_\b = \e_\b$ (since $\e_\a$ and $\e_\b$ are isotopic). Then $e^*$ is primitive with respect to $f_\d(T_\g)$, but since $\d^* \cap \e^* = \0$, $\e^*$ is also primitive with respect to $T_\g$ in $\mathcal{T}$, and thus the splitting is perturbed.

Suppose instead that each of $\zeta_\a$ and $\zeta_\b$ is disjoint from $\e_\a$ and $\e$, respectively. Then there is a second 1-bridge unlinked component of $f_\d(T_\a) \cup \overline{f_\d(T_\b)}$, and so $f_\d(T_\a) \cup \overline{f_\d(T_\b)}$ is a 3-component unlink.

In each possible situation, we lift $\mathcal{T}_\d$ to the double branched cover and find compressing curves in the splitting surface $\widetilde{\S}$, of genus two. Since every curve the genus two surface is equivariant under the hyperelliptic involution, each such curve is the lift of some shadow for $\mathcal{T}_\d$, and in every case we find that $\mathcal{T}_\d$ is a bridge trisection.

\end{proof}

\section{Theorems}

When weak reduction of a $n$-bridge-trisected surface $K$ in $S^4$ results in another bridge trisection (that is, all the tangles are trivial and their pairwise unions are unlinks), the forgetting realizes one of a small family of elementary operations on surface embeddings. The first operation, deperturbation, is exactly as described in \cite{mz}; this is an isotopy of the knotted surface. The second is removal of an unknotted $\R P^2$-summand. If we can perform this operation, we must have $K = K' \# P^+$ or $K = K' \# P^-$ for a surface $K'$ with bridge number $n' < n$ and an unknotted $\R P^2$ $P^+$ or $P^-$ (with normal Euler number $2$ or $-2$, respectively). The last operation, detubing, is of particular interest. Tubing a surface is surgery along an arc intersecting the surface only at its endpoints; detubing is realized as the dual surgery along a meridian of the tube. Notably, ribbon surfaces may be described as tubings of unlinked spheres, and one straightforward way to show a surface is ribbon is to find a series of appropriate detubings that yields an unlink of spheres.

To prove our first main result, we will need to appeal to a Theorem of Scharlemann concerning genus two 1-handlebodies in 3-manifolds.

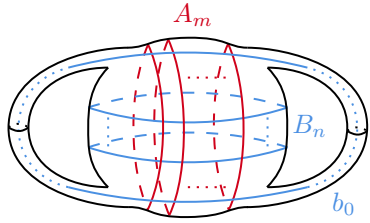
\begin{figure}[h]
    \centering

\tikzset{every picture/.style={line width=0.75pt}} 

\begin{tikzpicture}[x=0.75pt,y=0.75pt,yscale=-1,xscale=1]

\draw [color={rgb, 255:red, 208; green, 2; blue, 27 }  ,draw opacity=1 ] [dash pattern={on 4.5pt off 4.5pt}]  (90.14,88.33) .. controls (80.14,109.79) and (80.29,150.07) .. (90,171.29) ;
\draw [color={rgb, 255:red, 208; green, 2; blue, 27 }  ,draw opacity=1 ] [dash pattern={on 4.5pt off 4.5pt}]  (100.27,85.8) .. controls (90.27,107.25) and (91.05,152.97) .. (100.77,174.19) ;
\draw [color={rgb, 255:red, 208; green, 2; blue, 27 }  ,draw opacity=1 ] [dash pattern={on 4.5pt off 4.5pt}]  (130.27,87.94) .. controls (120.27,109.39) and (120.93,150.22) .. (130.64,171.44) ;
\draw  [dash pattern={on 0.84pt off 2.51pt}]  (20,130) .. controls (23.19,127.56) and (26.31,127.69) .. (30,130) ;
\draw    (190,131.78) .. controls (194,134.23) and (196.67,134.23) .. (200,131.78) ;
\draw  [dash pattern={on 0.84pt off 2.51pt}]  (190,131.78) .. controls (193.19,129.34) and (196.31,129.47) .. (200,131.78) ;
\draw [color={rgb, 255:red, 74; green, 144; blue, 226 }  ,draw opacity=1 ]   (60,120) .. controls (90,130) and (130,130) .. (160,120) ;
\draw [color={rgb, 255:red, 74; green, 144; blue, 226 }  ,draw opacity=1 ]   (60,140) .. controls (90,150) and (130,150) .. (160,140) ;
\draw [color={rgb, 255:red, 74; green, 144; blue, 226 }  ,draw opacity=1 ] [dash pattern={on 0.84pt off 2.51pt}]  (70,125) -- (70,140) ;
\draw [color={rgb, 255:red, 74; green, 144; blue, 226 }  ,draw opacity=1 ] [dash pattern={on 0.84pt off 2.51pt}]  (150,125) -- (150,140) ;
\draw [color={rgb, 255:red, 74; green, 144; blue, 226 }  ,draw opacity=1 ] [dash pattern={on 4.5pt off 4.5pt}]  (60,120) .. controls (90,110.11) and (130.31,110.08) .. (160,120) ;
\draw [color={rgb, 255:red, 74; green, 144; blue, 226 }  ,draw opacity=1 ] [dash pattern={on 4.5pt off 4.5pt}]  (60,140) .. controls (90,130.11) and (130.31,130.08) .. (160,140) ;
\draw [color={rgb, 255:red, 74; green, 144; blue, 226 }  ,draw opacity=1 ] [dash pattern={on 0.84pt off 2.51pt}]  (170,160) .. controls (190.15,154.92) and (195.08,139.85) .. (195,130) .. controls (194.92,120.15) and (190,105.08) .. (170,100) ;
\draw [color={rgb, 255:red, 74; green, 144; blue, 226 }  ,draw opacity=1 ] [dash pattern={on 0.84pt off 2.51pt}]  (50,160) .. controls (30.31,155.23) and (24.92,140) .. (25,130) .. controls (25.08,120) and (30.15,105.38) .. (50,100) ;
\draw    (20,130) .. controls (24,132.44) and (26.67,132.44) .. (30,130) ;
\draw [color={rgb, 255:red, 208; green, 2; blue, 27 }  ,draw opacity=1 ]   (90,88.33) .. controls (100.13,110.13) and (100.13,150) .. (89.86,171.29) ;
\draw [color={rgb, 255:red, 208; green, 2; blue, 27 }  ,draw opacity=1 ]   (130.13,87.94) .. controls (140.25,109.73) and (140.77,150.15) .. (130.5,171.44) ;
\draw    (70,160) .. controls (50.22,160) and (30,150) .. (30,130) .. controls (30,110) and (50,100) .. (70,100) ;
\draw    (80,170) .. controls (60.22,170) and (19.78,160) .. (20,130) .. controls (20.22,100) and (60,90) .. (80,90) ;
\draw    (70,100) .. controls (60.73,109.38) and (60,120.08) .. (60,130) .. controls (60,139.92) and (59.6,149.48) .. (70,160) ;
\draw    (150,100) .. controls (169.78,100.1) and (190.1,110) .. (190,130) .. controls (189.9,150) and (170,160.1) .. (150,160) ;
\draw    (140,90) .. controls (159.78,90.1) and (200.37,100) .. (200,130) .. controls (199.63,160) and (159.61,170.1) .. (139.61,170) ;
\draw    (149.71,160) .. controls (159.56,150.13) and (159.8,140.05) .. (159.86,130) .. controls (159.91,119.95) and (159.82,110.02) .. (150,100) ;
\draw [color={rgb, 255:red, 208; green, 2; blue, 27 }  ,draw opacity=1 ]   (100.13,85.8) .. controls (110.25,107.59) and (110.89,152.9) .. (100.63,174.19) ;
\draw    (79.81,90) .. controls (91.31,89.75) and (90,85.25) .. (110,85) .. controls (130,84.75) and (129.81,89.92) .. (140,90) ;
\draw    (80,170) .. controls (91.5,169.75) and (90,175.25) .. (110,175) .. controls (130,174.75) and (130,169.92) .. (140.19,170) ;
\draw [color={rgb, 255:red, 74; green, 144; blue, 226 }  ,draw opacity=1 ]   (50,160) .. controls (90.73,170.18) and (130,169.45) .. (170,160) ;
\draw [color={rgb, 255:red, 74; green, 144; blue, 226 }  ,draw opacity=1 ]   (50,100) .. controls (90,89.82) and (130.18,90) .. (170,100) ;
\draw [color={rgb, 255:red, 208; green, 2; blue, 27 }  ,draw opacity=1 ] [dash pattern={on 0.84pt off 2.51pt}]  (110,105) -- (130,105) ;
\draw [color={rgb, 255:red, 208; green, 2; blue, 27 }  ,draw opacity=1 ] [dash pattern={on 0.84pt off 2.51pt}]  (110,160) -- (130,160) ;

\draw (161,122.4) node [anchor=north west][inner sep=0.75pt]    {$\textcolor[rgb]{0.29,0.56,0.89}{B_{n}}$};
\draw (181,161.4) node [anchor=north west][inner sep=0.75pt]    {$\textcolor[rgb]{0.29,0.56,0.89}{b_{0}}$};
\draw (101,66.4) node [anchor=north west][inner sep=0.75pt]    {$\textcolor[rgb]{0.82,0.01,0.11}{A_{m}}$};

\end{tikzpicture}
    \caption{Two systems of curves in the boundary of a genus two 1-handlebody $N$ as in Theorem \ref{thm:schar}}
    \label{fig:schar}
\end{figure}

\begin{thm}[\cite{schar} Theorem 1.1]
\label{thm:schar}
    Let $A_m, B_n$ (with $m$ odd) be systems of curves on the boundary of a genus two 1-handlebody $N$ in an oriented 3-manifold $M$ such that $A_m$ consists of $m$ parallel curves, each bounding a disk in $N$ and separating $\del N$ into a pair of tori and such that $B_n$ consists of $b_0$, a curve running once over each 1-handle of $N$, and $n$ parallel curves disjoint from $b_0$ and bounding disks in the 0-handle of $N$. See Figure \ref{fig:schar}. Suppose that each of $A_m$ and $B_n$ bounds a planar surface in the closure of $M \bs N$. Then $a_1$ and $b_0$ (possibly after twisting around the co-cores of the 1-handles of $N$) bound a pair of disks intersecting in a single arc.
\end{thm}

Using Corollary \ref{cor:prim} and the Theorem \ref{thm:schar}, we are able to verify that when resolution of a bridge trisection under certain conditions yields another bridge trisection, we satisfy conditions laid out \cite{C&R} for one of the elementary operations described above.

\begin{thm}
\label{thm:ElemOp}

    Let $\mathcal{T} = (\S; T_\a, T_\b, T_\g)$ be bridge trisection of a surface knot $K$ in $S^4$ such that $T_\a \cup \overline{T_\b}$ has a 1-bridge component $U_\e$ disjoint from some arc $\d \in T_\g$. Suppose that $\mathcal{T}_\d$, the result of forgetting $\d$, is a bridge trisection. Then $\mathcal{T}$ is weakly reducible with weak-reducing arcs $\d$, $\e_\a$, and  $\e_\b$, where $\e_\a \cup \overline{\e_\b} = U_\e$. Further, $\mathcal{T}_\d$ is related to $\mathcal{T}$ by a detubing, a deperturbation, or the removal of an $\R P^2$ summand.

\end{thm}

\begin{proof}

To see that $\mathcal{T}$ is weakly-reducible, and that $\d$ is primitive with respect to the other two tangles, we refer to Corollary \ref{cor:prim}.

In $\mathcal{T}_\d$, the result of forgetting $\d$, we will record an arc $\d'$ such that resolving $f_\d(T_\a) \cup \overline{f_\d(T_\b)}$ along $\d'$ recovers $T_\a \cup \overline{T_\b}$. We can verify that such an arc exists by noting that there is a thin band $\mathfrak{b} \cong \d^* \x I$ in $\S$, and since $\d^*$ is primitive with respect to both $T_\a$ and $T_\b$ by Corollary \ref{cor:prim}, $\mathfrak{b}$ intersects $T_\a$ and $T_\b$ only at the endpoints of $\d^*$.

We take $\d'$ to be an arc in $\S_\d$ such that resolving $f_\d(T_\a) \cup \overline{f_\d(T_\b)}$ along $\d'$ recovers $T_\a \cup \overline{T_\b}$. Notice that $\d'$ has one endpoint in the interior of an arc $a^*$ of $f_\d(T_\a)^*$ and the other in the interior of an arc $b^*$ of $f_\d(T_\b)^*$. It is elsewhere disjoint from $f_\d(T_\a)^* \cup f_\d(T_\b)^*$.

We perturb $\d'$, maintaining the framing induced by $\S_\d$, by following isotopies for each of $a^*$ and $b^*$ that trace out their associated bridge disks. We thus obtain am arc $\rho$ transverse to $\S_\d$, intersecting it only in a single point $*$. We will view $\rho$ as $\rho_\a \cup \rho_\b$, where $\rho_\a$ is an arc in $f_\d(H_\a)$ with one endpoint at $*$ and intersecting $f_\d(T_\a)$ only at its other endpoint, and similarly for $\rho_\b \subset f_\d(H_\b)$ (recalling that $H_\a$ is the 3-ball containing $T_\a$, and likewise for $H_\b$ and $T_\b$). We then choose an arc $\rho_\g \subset f_\d(H_\g)$ transverse to $\S_\d$ and contained in a small neighborhood (in $f_\d(H_\g)$) of $*$.

We now wish to verify that the set $\{\rho_\a, \rho_\b, \rho_\g\}$ satisfies a set of conditions laid out in Definition 4.1 of \cite{C&R} for an \emph{admissible} collection of arcs. These conditions are as follows:
\begin{enumerate}
    \item Each $\rho_i$ is an increasing arc in $f_\d(H_i)$ with one endpoint at $*$ and interior disjoint from $f_\d(T_i)$.
    
    \item Each of $\rho_\a$ and $\rho_\b$ has one endpoint in $f_\d(T_\a)$ or $f_\d(T_\b)$, respectively, and $\rho_\g \cap f_\d(T_\g) = \0$.
    
    \item If both endpoints of $\rho = \rho_\a \cup \rho_\b$ lie in a single component $U$ of $f_\d(T_\a) \cup \overline{f_\d(T_\b)}$, then there is a sphere $S$ embedded in $f_\d(H_\a) \cup \overline{f_\d(H_\b)}$ satisfying
        \begin{itemize}
            \item $\rho \subset S$,
            \item $S \cap f_\d(T_\a) \cup \overline{f_\d(T_\b)}$ is an arc of $f_\d(T_\a) \cup \overline{f_\d(T_\b)}$ sharing endpoints with $\rho$, and 
            \item $S \cap \S_\d$ is a loop containing $*$ and some points of $K \cap \S_\d$.
        \end{itemize}

    \item If the endpoints of $\rho = \rho_\a \cup \rho_\b$ lie in distinct components $U_1$ and $U_2$ of $f_\d(T_\a^*) \cup \overline{f_\d(T_\b^*)}$, then there is a pair of disjoint embedded spheres $S_1$ and $S_2$ in $f_\d(H_\a) \cup \overline{f_\d(H_\b)}$ satisfying
        \begin{itemize}
            \item $S_i \cap f_\d(T_\a) \cup \overline{f_\d(T_\b)} = U_i$ for $i = 1,2$,
            \item $S_1 \cup S_2 \cap \rho = \del \rho$, and
            \item For $i = 1,2$, $S_i \cap \S_\d$ is a loop containing some points of $K \cap \S_\d$ but not $*$.
        \end{itemize}
\end{enumerate}

The first two conditions are already met by our construction, so we turn our attention to conditions (3) and (4), depending on how whether the endpoints of $\rho$ lie in the same or different components of $f_\d(T_\a) \cup \overline{f_\d(T_\b)}$.

In the case that both the endpoints of $\rho$ lie in the same component $U$ of $f_\d(T_\a) \cup \overline{f_\d(T_\b)}$, we take a tubular neighborhood $N$ in $S^3$ (given by the pairing $f_\d(H_\a) \cup \overline{f_\d(H_\b)}$) of the graph $U \cup \rho$, noting that $N$ is a 3-dimensional 1-handlebody of genus two. Indeed, since $\del \rho$ separates $U$ into two arcs, $\zeta_1$ and $\zeta_2$, we may view $N = \nu(\rho) \cup \nu(\zeta_1) \cup \nu(\zeta_2)$ as a 0-handle and two one-handles $Z_1$ and $Z_2$, with each $Z_i$ isotopic to a neighborhood of $\zeta_i$. We specify two curves on $\del N$: $r$ is a curve in $\del \nu(\rho)$ away from the attaching regions of $Z_1$ and $Z_2$ and separating them from one another. We will let $s$ be a curve isotopic to $L$ and 0-framed in $\del N$ (so $s$ runs along each of the 1-handles once). Note that since each of $r$ and $s$ is an unknot, each bounds a disk in $\overline{S^3 \bs N}$. Further, we see that $\rho$ is isotopic to an arc $\rho'$ of $r$ with endpoints $r \cap s$, as $\nu(\rho)$ is a 3-ball. Likewise, we may view $s$ as $\eta_1 \cup \eta_2$, where each $\eta_i$ is an arc with endpoints in $r$ isotopic to $\zeta_i$.

We quickly see that our system of curves satisfies the needed conditions, and we obtain disks $E$ with $r = \del E$ and $F$ with $s = \del F$ in $\overline{S^3 \bs N}$ such that $E \cap F$ is a single arc $\sigma$. Indeed, since $\del \sigma = r \cap s$, we have a pair of disks $D_1$ and $D_2$ such that for $i \in \{1,2\}$, $D_i$ is bounded by $\sigma \cup \eta_i$. There is also a disk bounded by $\sigma \cup \rho'$, which traces out an isotopy, so we also have disks $D_1'$ and $D_2'$ bounded by $\rho' \cup \eta_1$ and $\rho' \cup \eta_2$, respectively. Taking $D_1' \cup_{\rho'} D_2'$ and doubling, we find a sphere containing $s \cup \rho'$. So $s \cup \rho'$ is a trivial spatial graph in the sense of \cite{Ozawa}, meaning it is a graph embedded in $S^3$ contained in some 2-sphere. This graph is isotopic to $U \cup \rho$, which is therefore also a trivial spatial graph, and is in bridge position with respect to $\S_\d$ (by construction). By Theorem 1.1 of \cite{Ozawa}, there is a sphere $S$ containing $U \cup \rho$ intersecting $\S_\d$ in a single loop (which will of course contain all points of $(U \cup \rho) \cap \S_\d$). Thus, condition (3) of Definition 4.1 in \cite{C&R} is met.

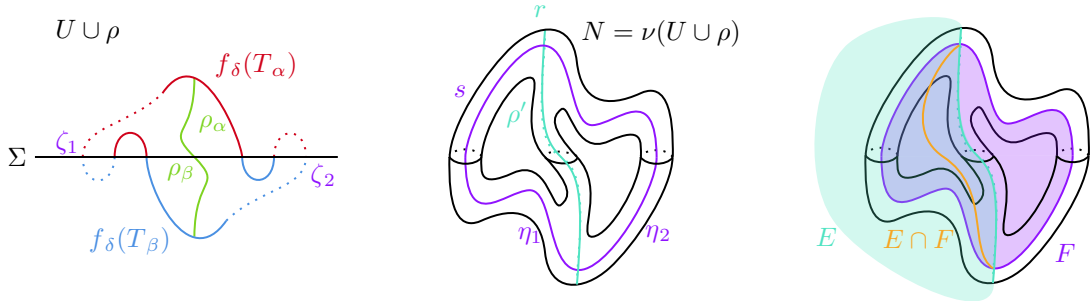
\begin{figure}[h]
    \centering

\tikzset{every picture/.style={line width=0.75pt}} 

\begin{tikzpicture}[x=0.6pt,y=0.6pt,yscale=-1,xscale=1]

\draw [color={rgb, 255:red, 80; green, 227; blue, 194 }  ,draw opacity=1 ] [dash pattern={on 0.84pt off 2.51pt}]  (360,180) .. controls (359.8,159.6) and (369.4,110.8) .. (350,100) ;
\draw [color={rgb, 255:red, 80; green, 227; blue, 194 }  ,draw opacity=1 ] [dash pattern={on 0.84pt off 2.51pt}]  (340,20) .. controls (340.6,44) and (330.2,89.6) .. (350,100) ;
\draw [color={rgb, 255:red, 126; green, 211; blue, 33 }  ,draw opacity=1 ]   (120,50) .. controls (120.38,89.02) and (100.79,80.47) .. (120,100) ;
\draw [color={rgb, 255:red, 126; green, 211; blue, 33 }  ,draw opacity=1 ]   (120,100) .. controls (138.37,119.1) and (120.38,110.98) .. (120,150) ;
\draw [color={rgb, 255:red, 208; green, 2; blue, 27 }  ,draw opacity=1 ] [dash pattern={on 0.84pt off 2.51pt}]  (170,100) .. controls (170.2,80) and (190.2,80.4) .. (190,100) ;
\draw [color={rgb, 255:red, 208; green, 2; blue, 27 }  ,draw opacity=1 ]   (100,60) .. controls (124.6,24) and (149.8,91.2) .. (150,100) ;
\draw [color={rgb, 255:red, 74; green, 144; blue, 226 }  ,draw opacity=1 ]   (140,140) .. controls (116.2,174) and (90.2,119.6) .. (90,100) ;
\draw [color={rgb, 255:red, 208; green, 2; blue, 27 }  ,draw opacity=1 ] [dash pattern={on 0.84pt off 2.51pt}]  (50,100) .. controls (50,80.24) and (79.94,79.96) .. (100,60) ;
\draw [color={rgb, 255:red, 74; green, 144; blue, 226 }  ,draw opacity=1 ] [dash pattern={on 0.84pt off 2.51pt}]  (140,140) .. controls (160,120) and (190,120) .. (190,100) ;
\draw [color={rgb, 255:red, 74; green, 144; blue, 226 }  ,draw opacity=1 ] [dash pattern={on 0.84pt off 2.51pt}]  (70,100) .. controls (69.4,119.6) and (49.8,119.6) .. (50,99.76) ;
\draw    (280,100) .. controls (280,110) and (278.25,149.75) .. (300,140) .. controls (321.75,130.25) and (337.22,134.33) .. (340,160) .. controls (342.78,185.67) and (363.89,184.33) .. (380,170) .. controls (396.11,155.67) and (419.75,120.25) .. (420,100) ;
\draw [color={rgb, 255:red, 144; green, 19; blue, 254 }  ,draw opacity=1 ]   (290,100.6) .. controls (290.25,149.1) and (314.75,121.6) .. (330,120.85) .. controls (345.25,120.1) and (341.25,168.6) .. (360,170.85) .. controls (378.75,173.1) and (409.86,120.99) .. (410,100.14) ;
\draw [color={rgb, 255:red, 208; green, 2; blue, 27 }  ,draw opacity=1 ]   (70,100) .. controls (69.8,79.6) and (89.8,80) .. (90,100) ;
\draw [color={rgb, 255:red, 74; green, 144; blue, 226 }  ,draw opacity=1 ]   (150,100) .. controls (150.2,118.8) and (169.8,119.6) .. (170,100) ;
\draw    (20,100) -- (210,100) ;
\draw    (300,100) .. controls (300,110) and (298.25,119.75) .. (320,110) .. controls (341.75,100.25) and (344.25,134.75) .. (350,130) .. controls (355.75,125.25) and (357.25,118.75) .. (340,100) ;
\draw    (360,100) .. controls (378.75,119.75) and (359.25,151.25) .. (370,150) .. controls (380.75,148.75) and (399.25,117.75) .. (400,100) ;
\draw    (400.15,99.99) .. controls (400.11,89.99) and (401.83,80.23) .. (380.11,90.06) .. controls (358.4,99.9) and (355.77,65.41) .. (350.04,70.18) .. controls (344.31,74.95) and (342.83,81.46) .. (360.15,100.14) ;
\draw    (340.15,100.22) .. controls (321.33,80.54) and (340.71,48.96) .. (329.96,50.26) .. controls (319.22,51.55) and (300.84,82.62) .. (300.15,100.37) ;
\draw    (280,100) .. controls (283.8,105.6) and (296.2,105.6) .. (300.15,100.37) ;
\draw    (340,100.37) .. controls (346.6,105.6) and (358.6,104.8) .. (360.15,100.74) ;
\draw    (399.85,100) .. controls (403.65,105.6) and (416.05,105.6) .. (420,100.37) ;
\draw  [dash pattern={on 0.84pt off 2.51pt}]  (420,100) .. controls (416.22,94.38) and (403.82,94.33) .. (399.85,99.54) ;
\draw  [dash pattern={on 0.84pt off 2.51pt}]  (360,99.37) .. controls (354.2,96.4) and (343.8,95.2) .. (339.85,98.91) ;
\draw  [dash pattern={on 0.84pt off 2.51pt}]  (300.15,99.48) .. controls (296.38,93.86) and (283.98,93.81) .. (280,99.02) ;
\draw [color={rgb, 255:red, 80; green, 227; blue, 194 }  ,draw opacity=1 ]   (340,20) .. controls (335.36,79.82) and (339.83,89.67) .. (350,100) .. controls (360.17,110.33) and (365,119.82) .. (360,180) ;
\draw    (420.15,99.91) .. controls (420.11,89.91) and (421.71,50.15) .. (400,59.99) .. controls (378.29,69.82) and (362.8,65.8) .. (359.92,40.14) .. controls (357.05,14.48) and (335.94,15.9) .. (319.88,30.29) .. controls (303.83,44.69) and (280.33,80.2) .. (280.15,100.45) ;
\draw [color={rgb, 255:red, 144; green, 19; blue, 254 }  ,draw opacity=1 ]   (410,100.14) .. controls (409.56,51.63) and (385.17,79.22) .. (369.92,80.03) .. controls (354.67,80.84) and (358.49,32.32) .. (339.73,30.15) .. controls (320.97,27.97) and (290.06,79.74) .. (290,100.6) ;
\draw [color={rgb, 255:red, 80; green, 227; blue, 194 }  ,draw opacity=1 ] [dash pattern={on 0.84pt off 2.51pt}]  (621,179.01) .. controls (620.8,158.61) and (630.4,109.81) .. (611,99.01) ;
\draw [color={rgb, 255:red, 80; green, 227; blue, 194 }  ,draw opacity=1 ] [dash pattern={on 0.84pt off 2.51pt}]  (601,19.01) .. controls (601.6,43.01) and (591.2,88.61) .. (611,99.01) ;
\draw    (541,99.01) .. controls (541,109.01) and (539.25,148.76) .. (561,139.01) .. controls (582.75,129.26) and (598.22,133.34) .. (601,159.01) .. controls (603.78,184.67) and (624.89,183.34) .. (641,169.01) .. controls (657.11,154.67) and (680.75,119.26) .. (681,99.01) ;
\draw [color={rgb, 255:red, 144; green, 19; blue, 254 }  ,draw opacity=1 ][fill={rgb, 255:red, 144; green, 19; blue, 254 }  ,fill opacity=0.25 ]   (551,99.6) .. controls (551.25,148.11) and (575.75,120.61) .. (591,119.86) .. controls (606.25,119.11) and (602.25,167.61) .. (621,169.86) .. controls (639.75,172.11) and (670.86,120) .. (671,99.14) ;
\draw    (561,99.01) .. controls (561,109.01) and (559.25,118.76) .. (581,109.01) .. controls (602.75,99.26) and (605.25,133.76) .. (611,129.01) .. controls (616.75,124.26) and (618.25,117.76) .. (601,99.01) ;
\draw    (621,99.01) .. controls (639.75,118.76) and (620.25,150.26) .. (631,149.01) .. controls (641.75,147.76) and (660.25,116.76) .. (661,99.01) ;
\draw    (661.15,98.99) .. controls (661.11,88.99) and (662.83,79.24) .. (641.11,89.07) .. controls (619.4,98.9) and (616.77,64.41) .. (611.04,69.18) .. controls (605.31,73.96) and (603.83,80.46) .. (621.15,99.15) ;
\draw    (601.15,99.22) .. controls (582.33,79.54) and (601.71,47.97) .. (590.96,49.26) .. controls (580.22,50.55) and (561.84,81.62) .. (561.15,99.38) ;
\draw    (541,99.01) .. controls (544.8,104.61) and (557.2,104.61) .. (561.15,99.38) ;
\draw    (601,99.38) .. controls (607.6,104.61) and (619.6,103.81) .. (621.15,99.75) ;
\draw    (660.85,99.01) .. controls (664.65,104.61) and (677.05,104.61) .. (681,99.38) ;
\draw  [dash pattern={on 0.84pt off 2.51pt}]  (681,99.01) .. controls (677.22,93.39) and (664.82,93.33) .. (660.85,98.55) ;
\draw  [dash pattern={on 0.84pt off 2.51pt}]  (621,98.37) .. controls (615.2,95.41) and (604.8,94.21) .. (600.85,97.92) ;
\draw  [dash pattern={on 0.84pt off 2.51pt}]  (561.15,98.48) .. controls (557.38,92.87) and (544.98,92.81) .. (541,98.02) ;
\draw    (681.15,98.91) .. controls (681.11,88.91) and (682.71,49.16) .. (661,58.99) .. controls (639.29,68.83) and (623.8,64.8) .. (620.92,39.15) .. controls (618.05,13.49) and (596.94,14.9) .. (580.88,29.3) .. controls (564.83,43.69) and (541.33,79.2) .. (541.15,99.45) ;
\draw [color={rgb, 255:red, 144; green, 19; blue, 254 }  ,draw opacity=1 ][fill={rgb, 255:red, 144; green, 19; blue, 254 }  ,fill opacity=0.25 ]   (671,99.14) .. controls (670.56,50.63) and (646.17,78.23) .. (630.92,79.04) .. controls (615.67,79.84) and (619.49,31.33) .. (600.73,29.15) .. controls (581.97,26.97) and (551.06,78.74) .. (551,99.6) ;
\draw  [color={rgb, 255:red, 80; green, 227; blue, 194 }  ,draw opacity=0.25 ][fill={rgb, 255:red, 80; green, 227; blue, 194 }  ,fill opacity=0.25 ] (531,29.01) .. controls (552.75,11.38) and (600.5,8.38) .. (601,19.01) .. controls (601.5,29.63) and (591,78.88) .. (611,99.01) .. controls (631,119.13) and (620.75,148.13) .. (621,179.01) .. controls (621.25,209.88) and (541,169.01) .. (521,139.01) .. controls (501,109.01) and (509.25,46.63) .. (531,29.01) -- cycle ;
\draw [color={rgb, 255:red, 80; green, 227; blue, 194 }  ,draw opacity=1 ]   (601,19.01) .. controls (596.36,78.82) and (600.83,88.67) .. (611,99.01) .. controls (621.17,109.34) and (626,118.82) .. (621,179.01) ;
\draw [color={rgb, 255:red, 245; green, 166; blue, 35 }  ,draw opacity=1 ]   (600,30) .. controls (590,34.75) and (560,80.13) .. (590,100) .. controls (620,119.88) and (604.82,158.55) .. (621,169.86) ;

\draw (131,32.4) node [anchor=north west][inner sep=0.75pt]    {$\textcolor[rgb]{0.82,0.01,0.11}{f}\textcolor[rgb]{0.82,0.01,0.11}{_{\delta }}\textcolor[rgb]{0.82,0.01,0.11}{(}\textcolor[rgb]{0.82,0.01,0.11}{T}\textcolor[rgb]{0.82,0.01,0.11}{_{\alpha }}\textcolor[rgb]{0.82,0.01,0.11}{)}$};
\draw (52,142.4) node [anchor=north west][inner sep=0.75pt]    {$\textcolor[rgb]{0.29,0.56,0.89}{f}\textcolor[rgb]{0.29,0.56,0.89}{_{\delta }}\textcolor[rgb]{0.29,0.56,0.89}{(}\textcolor[rgb]{0.29,0.56,0.89}{T}\textcolor[rgb]{0.29,0.56,0.89}{_{\beta }}\textcolor[rgb]{0.29,0.56,0.89}{)}$};
\draw (1,92.4) node [anchor=north west][inner sep=0.75pt]    {$\Sigma $};
\draw (121,72.4) node [anchor=north west][inner sep=0.75pt]    {$\textcolor[rgb]{0.49,0.83,0.13}{\rho }\textcolor[rgb]{0.49,0.83,0.13}{_{\alpha }}$};
\draw (101,102.4) node [anchor=north west][inner sep=0.75pt]    {$\textcolor[rgb]{0.49,0.83,0.13}{\rho }\textcolor[rgb]{0.49,0.83,0.13}{_{\beta }}$};
\draw (31,81.4) node [anchor=north west][inner sep=0.75pt]    {$\textcolor[rgb]{0.56,0.07,1}{\zeta _{1}}$};
\draw (191,102.4) node [anchor=north west][inner sep=0.75pt]    {$\textcolor[rgb]{0.56,0.07,1}{\zeta }\textcolor[rgb]{0.56,0.07,1}{_{2}}$};
\draw (281,53.4) node [anchor=north west][inner sep=0.75pt]    {$\textcolor[rgb]{0.56,0.07,1}{s}$};
\draw (331,3.4) node [anchor=north west][inner sep=0.75pt]    {$\textcolor[rgb]{0.31,0.89,0.76}{r}$};
\draw (314,63.4) node [anchor=north west][inner sep=0.75pt]    {$\textcolor[rgb]{0.31,0.89,0.76}{\rho '}$};
\draw (321,142.4) node [anchor=north west][inner sep=0.75pt]    {$\textcolor[rgb]{0.56,0.07,1}{\eta _{1}}$};
\draw (401,141.4) node [anchor=north west][inner sep=0.75pt]    {$\textcolor[rgb]{0.56,0.07,1}{\eta }\textcolor[rgb]{0.56,0.07,1}{_{2}}$};
\draw (31,12.4) node [anchor=north west][inner sep=0.75pt]    {$U\cup \rho \ $};
\draw (361,12.4) node [anchor=north west][inner sep=0.75pt]    {$N=\nu ( U\cup \rho ) \ $};
\draw (508,142.4) node [anchor=north west][inner sep=0.75pt]    {$\textcolor[rgb]{0.31,0.89,0.76}{E}$};
\draw (658,152.4) node [anchor=north west][inner sep=0.75pt]    {$\textcolor[rgb]{0.56,0.07,1}{F}$};
\draw (550,144) node [anchor=north west][inner sep=0.75pt]    {$\textcolor[rgb]{0.96,0.65,0.14}{E\cap F}$};

\end{tikzpicture}

    \caption{Disks bounded by curves in $N$}
    \label{fig:ScharDisks}
\end{figure}

Consider now the case where $\rho$ connects distinct components $U_1$ and $U_2$ of $f_\d(T_\a) \cup \overline{f_\d(T_\b)}$, corresponding to condition (4) of Definition 4.1 in \cite{C&R}. If $U_1 \cup U_2$ is in 2-bridge position with respect to $\S_\d$ (so each of $U_1$ and $U_2$ is in 1-bridge position), then the union $D_1$ of bridge disks for each of the $\a$ and $\b$ arcs of $U_1$ is a disk spanning $U_1$ that can be doubled to form a sphere containing $U_1$, intersecting $\S_\d$ in a single arc, and disjoint from $\rho$. We produce a similar sphere for $U_2$, and see that in the minimal-bridge case, condition (4) is satisfied.

Suppose now that we have some perturbation of $U_1 \cup \rho \cup U_2$ (which will necessarily be a perturbation of either $U_1$ or $U_2$, as we suppose that $\rho$ always intersects $\S_\d$ in a single point). Suppose that prior to the perturbation, $U_1$ and $U_2$ each lies in a sphere, $S_1$ and $S_2$, respectively, such that each $S_i$ intersects $\S$ in a single loop and is disjoint from the interior of $\rho$, and such that $S_1 \cap S_2 = \0$. We will take a set of shadows for $f_\d(T_\a)$ arising from $S_1$ and $S_2$ in the following ways: take a disk $B_a$ in the relevant $S_i$ cobounded by an arc $a$ of $f_\d(T_\a)$ and an arc $a^*$ in $\S_\d$ to be a bridge disk for $a$ if $B_a$ is disjoint from $f_\d(T_\a)$ in its interior. Otherwise, push the interior of $B_a$ slightly off $S_i$ to obtain a bridge disk $B_a'$ for $a$. Do this for each arc of $f_\d(T_\a)$, pushing off in the same direction each time and ensuring that all bridge disks chosen are pairwise disjoint. Then perform an analogous process for each arc of $f_\d(T_\b)$ (pushing off in the opposite direction to that used for $f_\d(T_\a)$). We now see that the union of the shadows obtained for each of $U_1$ and $U_2$ forms a loop isotopic to the loops $\S_\d \cap S_1$ and $\S_\d \cap S_2$, respectively, in $\S_\d$, and that these loops are disjoint from one another and from shadows for $\rho$.

We will suppose that the perturbation happens from the $\a$ side, meaning that it is specified by a descending arc $\xi$ with one endpoint in the interior of an arc $a$ of $f_\d(T_\a)$ (say in $U_1$), the other in $\S_\d$, and elsewhere disjoint from both tangles. An analogous argument will apply to perturbations from the $\b$ side. Note that we can project $\xi$ onto $\S_\d$ to obtain a shadow $\xi^*$ with one endpoint in the interior of a shadow $a^*$ for $a$ and the other endpoint away from all shadows. Since this endpoint is away from $\a$ shadows, we may slide each $\a$ shadow intersecting $\xi^*$ in its interior along $\xi^*$ and over the endpoint, removing all such intersections. This corresponds to isotopy in $f_\d(T_\a)$ of the bridge disks, and we can view it as $\xi^*$ emerging from its endpoint in $a^*$ and pushing other $\a$ shadows before it. We then perform the perturbation (on the shadow level) by summing $a^*$ with a surface-framed band with core $\xi^*$, breaking the new arc near the endpoint of $\xi^*$ not in $a^*$, and connecting the new endpoints via a new $\b$ shadow $b^*$. Finally, for any intersection of $\xi^*$ with a $\b$ shadow, we slide that shadow over $b^*$, removing the intersection. An example of this process is outlined on the level of shadows in Figure \ref{fig:ShadowPert}. This corresponds to sliding bridge disks for $f_\d(T_\b)$ over that for the new arc $b$.

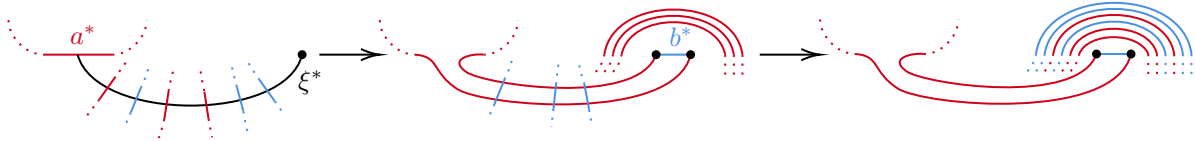
\begin{figure}[h]
    \centering

\tikzset{every picture/.style={line width=0.75pt}} 

\begin{tikzpicture}[x=0.65pt,y=0.65pt,yscale=-1,xscale=1]

\draw [color={rgb, 255:red, 208; green, 2; blue, 27 }  ,draw opacity=1 ]   (30,60) -- (70,60) ;
\draw    (50,60) .. controls (60.5,99.17) and (170.5,99.5) .. (180,60) ;
\draw  [fill={rgb, 255:red, 0; green, 0; blue, 0 }  ,fill opacity=1 ] (178,60) .. controls (178,58.9) and (178.9,58) .. (180,58) .. controls (181.1,58) and (182,58.9) .. (182,60) .. controls (182,61.1) and (181.1,62) .. (180,62) .. controls (178.9,62) and (178,61.1) .. (178,60) -- cycle ;
\draw [color={rgb, 255:red, 208; green, 2; blue, 27 }  ,draw opacity=1 ]   (62.6,85.9) -- (72,72.7) ;
\draw [color={rgb, 255:red, 208; green, 2; blue, 27 }  ,draw opacity=1 ]   (125.8,97.1) -- (123.4,81.5) ;
\draw [color={rgb, 255:red, 74; green, 144; blue, 226 }  ,draw opacity=1 ]   (80.8,93.5) -- (86.8,78.7) ;
\draw [color={rgb, 255:red, 208; green, 2; blue, 27 }  ,draw opacity=1 ]   (100.2,97.3) -- (102.8,80.5) ;
\draw [color={rgb, 255:red, 74; green, 144; blue, 226 }  ,draw opacity=1 ]   (164.2,87) -- (155,74.2) ;
\draw [color={rgb, 255:red, 74; green, 144; blue, 226 }  ,draw opacity=1 ]   (146,92.8) -- (142.2,78.6) ;
\draw [color={rgb, 255:red, 208; green, 2; blue, 27 }  ,draw opacity=1 ] [dash pattern={on 0.84pt off 2.51pt}]  (70,60) .. controls (80.2,60.1) and (90,49.9) .. (90,40) ;
\draw [color={rgb, 255:red, 208; green, 2; blue, 27 }  ,draw opacity=1 ] [dash pattern={on 0.84pt off 2.51pt}]  (10,40) .. controls (10,49.9) and (20.2,59.9) .. (30,60) ;
\draw [color={rgb, 255:red, 208; green, 2; blue, 27 }  ,draw opacity=1 ] [dash pattern={on 0.84pt off 2.51pt}]  (62.6,85.9) .. controls (60,89.5) and (58.2,91.5) .. (54.8,96.1) ;
\draw [color={rgb, 255:red, 208; green, 2; blue, 27 }  ,draw opacity=1 ] [dash pattern={on 0.84pt off 2.51pt}]  (77.8,64.7) .. controls (75.2,68.3) and (75.4,68.1) .. (72,72.7) ;
\draw [color={rgb, 255:red, 208; green, 2; blue, 27 }  ,draw opacity=1 ] [dash pattern={on 0.84pt off 2.51pt}]  (100.2,97.3) .. controls (99.2,101.7) and (98.6,103.9) .. (97,108.3) ;
\draw [color={rgb, 255:red, 208; green, 2; blue, 27 }  ,draw opacity=1 ] [dash pattern={on 0.84pt off 2.51pt}]  (105.2,69.9) .. controls (104.2,73.9) and (103.6,75.7) .. (102.8,80.5) ;
\draw [color={rgb, 255:red, 208; green, 2; blue, 27 }  ,draw opacity=1 ] [dash pattern={on 0.84pt off 2.51pt}]  (125.8,97.1) .. controls (126.4,100.9) and (126.8,102.5) .. (127.8,107.9) ;
\draw [color={rgb, 255:red, 208; green, 2; blue, 27 }  ,draw opacity=1 ] [dash pattern={on 0.84pt off 2.51pt}]  (121.2,69.7) .. controls (121.8,74.1) and (122.6,76.9) .. (123.4,81.5) ;
\draw [color={rgb, 255:red, 74; green, 144; blue, 226 }  ,draw opacity=1 ] [dash pattern={on 0.84pt off 2.51pt}]  (146,92.8) .. controls (146.53,96.85) and (146.98,99.55) .. (148.38,104.75) ;
\draw [color={rgb, 255:red, 74; green, 144; blue, 226 }  ,draw opacity=1 ] [dash pattern={on 0.84pt off 2.51pt}]  (139.4,67.4) .. controls (140.2,71) and (140.8,73.4) .. (142.2,78.6) ;
\draw [color={rgb, 255:red, 74; green, 144; blue, 226 }  ,draw opacity=1 ] [dash pattern={on 0.84pt off 2.51pt}]  (164.2,87) .. controls (166.6,90.2) and (167.6,91.8) .. (170,95.4) ;
\draw [color={rgb, 255:red, 74; green, 144; blue, 226 }  ,draw opacity=1 ] [dash pattern={on 0.84pt off 2.51pt}]  (149,65) .. controls (151.2,68.4) and (152.4,70.8) .. (155,74.2) ;
\draw [color={rgb, 255:red, 74; green, 144; blue, 226 }  ,draw opacity=1 ] [dash pattern={on 0.84pt off 2.51pt}]  (80.8,93.5) .. controls (79.2,97.9) and (78.4,99.3) .. (76.4,102.7) ;
\draw [color={rgb, 255:red, 74; green, 144; blue, 226 }  ,draw opacity=1 ] [dash pattern={on 0.84pt off 2.51pt}]  (90.8,68.5) .. controls (89.2,72.9) and (88.8,75.3) .. (86.8,78.7) ;
\draw [color={rgb, 255:red, 208; green, 2; blue, 27 }  ,draw opacity=1 ] [dash pattern={on 0.84pt off 2.51pt}]  (285,60) .. controls (295.2,60.1) and (305,49.9) .. (305,40) ;
\draw [color={rgb, 255:red, 208; green, 2; blue, 27 }  ,draw opacity=1 ] [dash pattern={on 0.84pt off 2.51pt}]  (225,40) .. controls (225,49.9) and (235.2,59.9) .. (245,60) ;
\draw [color={rgb, 255:red, 208; green, 2; blue, 27 }  ,draw opacity=1 ]   (245,60) .. controls (260,60.8) and (254.67,70.17) .. (270,80) .. controls (285.33,89.83) and (395.33,99.83) .. (405,60) ;
\draw [color={rgb, 255:red, 208; green, 2; blue, 27 }  ,draw opacity=1 ]   (285,60) .. controls (278,59.17) and (264.33,62.5) .. (275,70) .. controls (285.67,77.5) and (374.33,90.5) .. (385,60) ;
\draw [color={rgb, 255:red, 74; green, 144; blue, 226 }  ,draw opacity=1 ]   (385,60) -- (405,60) ;
\draw  [fill={rgb, 255:red, 0; green, 0; blue, 0 }  ,fill opacity=1 ] (383,60) .. controls (383,58.9) and (383.9,58) .. (385,58) .. controls (386.1,58) and (387,58.9) .. (387,60) .. controls (387,61.1) and (386.1,62) .. (385,62) .. controls (383.9,62) and (383,61.1) .. (383,60) -- cycle ;
\draw  [fill={rgb, 255:red, 0; green, 0; blue, 0 }  ,fill opacity=1 ] (403,60) .. controls (403,58.9) and (403.9,58) .. (405,58) .. controls (406.1,58) and (407,58.9) .. (407,60) .. controls (407,61.1) and (406.1,62) .. (405,62) .. controls (403.9,62) and (403,61.1) .. (403,60) -- cycle ;
\draw [color={rgb, 255:red, 208; green, 2; blue, 27 }  ,draw opacity=1 ]   (365,60) .. controls (370.15,34.92) and (419.85,34.92) .. (425,60) ;
\draw [color={rgb, 255:red, 208; green, 2; blue, 27 }  ,draw opacity=1 ]   (360,60) .. controls (365.08,30) and (424.92,30) .. (430,60) ;
\draw [color={rgb, 255:red, 208; green, 2; blue, 27 }  ,draw opacity=1 ]   (355,60) .. controls (360,25.23) and (430,25.08) .. (435,60) ;
\draw [color={rgb, 255:red, 208; green, 2; blue, 27 }  ,draw opacity=1 ] [dash pattern={on 0.84pt off 2.51pt}]  (425,60) -- (425,75) ;
\draw [color={rgb, 255:red, 208; green, 2; blue, 27 }  ,draw opacity=1 ] [dash pattern={on 0.84pt off 2.51pt}]  (430,60) -- (430,75) ;
\draw [color={rgb, 255:red, 208; green, 2; blue, 27 }  ,draw opacity=1 ] [dash pattern={on 0.84pt off 2.51pt}]  (435,60) -- (435,75) ;
\draw [color={rgb, 255:red, 208; green, 2; blue, 27 }  ,draw opacity=1 ] [dash pattern={on 0.84pt off 2.51pt}]  (360,70) .. controls (363.42,66.17) and (364.08,64) .. (365,60) ;
\draw [color={rgb, 255:red, 208; green, 2; blue, 27 }  ,draw opacity=1 ] [dash pattern={on 0.84pt off 2.51pt}]  (350,70) .. controls (353.42,66.17) and (354.08,64) .. (355,60) ;
\draw [color={rgb, 255:red, 208; green, 2; blue, 27 }  ,draw opacity=1 ] [dash pattern={on 0.84pt off 2.51pt}]  (355,70) .. controls (358.42,66.17) and (359.08,64) .. (360,60) ;
\draw    (190,60) -- (223,60) ;
\draw [shift={(225,60)}, rotate = 180] [color={rgb, 255:red, 0; green, 0; blue, 0 }  ][line width=0.75]    (10.93,-3.29) .. controls (6.95,-1.4) and (3.31,-0.3) .. (0,0) .. controls (3.31,0.3) and (6.95,1.4) .. (10.93,3.29)   ;
\draw [color={rgb, 255:red, 208; green, 2; blue, 27 }  ,draw opacity=1 ] [dash pattern={on 0.84pt off 2.51pt}]  (540,59.62) .. controls (550.2,59.72) and (560,49.52) .. (560,39.62) ;
\draw [color={rgb, 255:red, 208; green, 2; blue, 27 }  ,draw opacity=1 ] [dash pattern={on 0.84pt off 2.51pt}]  (480,39.63) .. controls (480,49.53) and (490.2,59.53) .. (500,59.63) ;
\draw [color={rgb, 255:red, 208; green, 2; blue, 27 }  ,draw opacity=1 ]   (500,59.63) .. controls (515,60.43) and (509.67,69.79) .. (525,79.63) .. controls (540.33,89.46) and (650.33,99.46) .. (660,59.63) ;
\draw [color={rgb, 255:red, 208; green, 2; blue, 27 }  ,draw opacity=1 ]   (540,59.63) .. controls (533,58.79) and (519.33,62.13) .. (530,69.63) .. controls (540.67,77.13) and (629.33,90.13) .. (640,59.63) ;
\draw [color={rgb, 255:red, 74; green, 144; blue, 226 }  ,draw opacity=1 ]   (640,59.63) -- (660,59.63) ;
\draw  [fill={rgb, 255:red, 0; green, 0; blue, 0 }  ,fill opacity=1 ] (638,59.63) .. controls (638,58.52) and (638.9,57.63) .. (640,57.63) .. controls (641.1,57.63) and (642,58.52) .. (642,59.63) .. controls (642,60.73) and (641.1,61.63) .. (640,61.63) .. controls (638.9,61.63) and (638,60.73) .. (638,59.63) -- cycle ;
\draw  [fill={rgb, 255:red, 0; green, 0; blue, 0 }  ,fill opacity=1 ] (658,59.63) .. controls (658,58.52) and (658.9,57.63) .. (660,57.63) .. controls (661.1,57.63) and (662,58.52) .. (662,59.63) .. controls (662,60.73) and (661.1,61.63) .. (660,61.63) .. controls (658.9,61.63) and (658,60.73) .. (658,59.63) -- cycle ;
\draw [color={rgb, 255:red, 208; green, 2; blue, 27 }  ,draw opacity=1 ]   (630,59.63) .. controls (635,46.79) and (664.33,46.46) .. (670,59.63) ;
\draw [color={rgb, 255:red, 208; green, 2; blue, 27 }  ,draw opacity=1 ]   (625,59.63) .. controls (629.67,40.13) and (670.33,40.13) .. (675,59.63) ;
\draw [color={rgb, 255:red, 74; green, 144; blue, 226 }  ,draw opacity=1 ]   (620,59.63) .. controls (625,34.46) and (675,34.79) .. (680,59.63) ;
\draw [color={rgb, 255:red, 208; green, 2; blue, 27 }  ,draw opacity=1 ]   (615,59.63) .. controls (620,29.46) and (680,29.46) .. (685,59.63) ;
\draw [color={rgb, 255:red, 74; green, 144; blue, 226 }  ,draw opacity=1 ]   (610,59.63) .. controls (615.33,24.79) and (685.33,24.79) .. (690,59.63) ;
\draw [color={rgb, 255:red, 74; green, 144; blue, 226 }  ,draw opacity=1 ]   (605,59.63) .. controls (610,20.46) and (690,19.79) .. (695,59.63) ;
\draw [color={rgb, 255:red, 208; green, 2; blue, 27 }  ,draw opacity=1 ] [dash pattern={on 0.84pt off 2.51pt}]  (670,59.63) -- (670,74.63) ;
\draw [color={rgb, 255:red, 208; green, 2; blue, 27 }  ,draw opacity=1 ] [dash pattern={on 0.84pt off 2.51pt}]  (675,59.63) -- (675,74.63) ;
\draw [color={rgb, 255:red, 208; green, 2; blue, 27 }  ,draw opacity=1 ] [dash pattern={on 0.84pt off 2.51pt}]  (685,59.63) -- (685,74.63) ;
\draw [color={rgb, 255:red, 74; green, 144; blue, 226 }  ,draw opacity=1 ] [dash pattern={on 0.84pt off 2.51pt}]  (680,59.63) -- (680,74.63) ;
\draw [color={rgb, 255:red, 74; green, 144; blue, 226 }  ,draw opacity=1 ] [dash pattern={on 0.84pt off 2.51pt}]  (690,59.63) -- (690,74.63) ;
\draw [color={rgb, 255:red, 74; green, 144; blue, 226 }  ,draw opacity=1 ] [dash pattern={on 0.84pt off 2.51pt}]  (695,59.63) -- (695,74.63) ;
\draw [color={rgb, 255:red, 208; green, 2; blue, 27 }  ,draw opacity=1 ] [dash pattern={on 0.84pt off 2.51pt}]  (625,69.63) .. controls (628.42,65.79) and (629.08,63.63) .. (630,59.63) ;
\draw [color={rgb, 255:red, 208; green, 2; blue, 27 }  ,draw opacity=1 ] [dash pattern={on 0.84pt off 2.51pt}]  (620,69.63) .. controls (623.42,65.79) and (624.08,63.63) .. (625,59.63) ;
\draw [color={rgb, 255:red, 208; green, 2; blue, 27 }  ,draw opacity=1 ] [dash pattern={on 0.84pt off 2.51pt}]  (610,69.63) .. controls (613.42,65.79) and (614.08,63.63) .. (615,59.63) ;
\draw [color={rgb, 255:red, 74; green, 144; blue, 226 }  ,draw opacity=1 ] [dash pattern={on 0.84pt off 2.51pt}]  (615,69.63) .. controls (618.42,65.79) and (619.08,63.63) .. (620,59.63) ;
\draw [color={rgb, 255:red, 74; green, 144; blue, 226 }  ,draw opacity=1 ] [dash pattern={on 0.84pt off 2.51pt}]  (605,69.63) .. controls (608.42,65.79) and (609.08,63.63) .. (610,59.63) ;
\draw [color={rgb, 255:red, 74; green, 144; blue, 226 }  ,draw opacity=1 ] [dash pattern={on 0.84pt off 2.51pt}]  (600,69.63) .. controls (603.42,65.79) and (604.08,63.63) .. (605,59.63) ;
\draw [color={rgb, 255:red, 74; green, 144; blue, 226 }  ,draw opacity=1 ]   (291.4,88.6) -- (297.4,73.8) ;
\draw [color={rgb, 255:red, 74; green, 144; blue, 226 }  ,draw opacity=1 ]   (346.2,91.85) -- (343.43,76.97) ;
\draw [color={rgb, 255:red, 74; green, 144; blue, 226 }  ,draw opacity=1 ]   (325.55,91.36) -- (327.08,77.07) ;
\draw [color={rgb, 255:red, 74; green, 144; blue, 226 }  ,draw opacity=1 ] [dash pattern={on 0.84pt off 2.51pt}]  (325.55,91.36) .. controls (325,95.36) and (324.82,97.18) .. (324.45,101.91) ;
\draw [color={rgb, 255:red, 74; green, 144; blue, 226 }  ,draw opacity=1 ] [dash pattern={on 0.84pt off 2.51pt}]  (328.4,65.75) .. controls (327.8,70.25) and (327.5,72.35) .. (327.08,77.07) ;
\draw [color={rgb, 255:red, 74; green, 144; blue, 226 }  ,draw opacity=1 ] [dash pattern={on 0.84pt off 2.51pt}]  (346.2,91.85) .. controls (346.6,95.65) and (347.1,97.16) .. (348,100.86) ;
\draw [color={rgb, 255:red, 74; green, 144; blue, 226 }  ,draw opacity=1 ] [dash pattern={on 0.84pt off 2.51pt}]  (340.48,66.19) .. controls (341.33,70.38) and (342.29,73.24) .. (343.43,76.97) ;
\draw [color={rgb, 255:red, 74; green, 144; blue, 226 }  ,draw opacity=1 ] [dash pattern={on 0.84pt off 2.51pt}]  (291.4,88.6) .. controls (289.8,93) and (289,94.4) .. (287,97.8) ;
\draw [color={rgb, 255:red, 74; green, 144; blue, 226 }  ,draw opacity=1 ] [dash pattern={on 0.84pt off 2.51pt}]  (301.4,63.6) .. controls (299.8,68) and (299.4,70.4) .. (297.4,73.8) ;
\draw    (445,60) -- (478,60) ;
\draw [shift={(480,60)}, rotate = 180] [color={rgb, 255:red, 0; green, 0; blue, 0 }  ][line width=0.75]    (10.93,-3.29) .. controls (6.95,-1.4) and (3.31,-0.3) .. (0,0) .. controls (3.31,0.3) and (6.95,1.4) .. (10.93,3.29)   ;

\draw (44,41.4) node [anchor=north west][inner sep=0.75pt]    {$\textcolor[rgb]{0.82,0.01,0.11}{a^{*}}$};
\draw (176,67.4) node [anchor=north west][inner sep=0.75pt]    {$\xi ^{*}$};
\draw (391,42.4) node [anchor=north west][inner sep=0.75pt]    {$\textcolor[rgb]{0.29,0.56,0.89}{b^{*}}$};

\end{tikzpicture}
    \caption{Shadows for a perturbation}
    \label{fig:ShadowPert}
\end{figure}

In the end, we obtain a shadow diagram for the perturbed bridge splitting which forms two disjoint (crossingless) loops, each disjoint from shadows for $\rho$ (which go through the same process as arcs of $U_1$ and $U_2$ during the perturbation). We take a disk $D_1$ bounded by the shadows $U_1^*$ for $U_1$ and push its interior off in the $\a$-direction. Near the relevant bridge disk $B_a$ for each arc $a$ of $f_\d(T_\a)_\xi \cap U_1$, we isotope the interior of the disk to contain $a$. By pushing the interior of another copy of $D_1$ in the $\b$-direction and following the isotopy traced out by the bridge disks for $f_\d(T_\b)_\xi \cap U_1$, we obtain a sphere $S_1'$ containing $U_1$, intersecting $\S_\d$ in a single loop, and disjoint from $\rho$. An analogous process gets us $S_2'$, a sphere containing $U_2$, intersecting $\S_\d$ in a single loop, and disjoint from both $S_1'$ and $\rho$.

Thus, we see that the minimal splitting of $U_1 \cup \rho \cup U_2$ meets condition (4), and that perturbation preserves the condition. What remains is to prove that every splitting of $U_1 \cup \rho \cup U_2$ is a perturbation of the minimal one.

Recall from \cite{Ozawa} that for any trivial (i.e. unknotted, in the sense of a graph) graph in $S^3$ in bridge position with respect to a splitting surface $\S$, there is a sphere $S$ containing the graph and intersecting $\S$ in a single loop. In the case of $U_1 \cup \rho \cup U_2$ in non-minimal position, we see that every arc in $f_\d(T_\a)$ co-bounds a bigon in $S$ with some arc of $\S_\d$. We take an innermost such bigon bounded by $a \cup a^*$ for an arc $a$ of $f_\d(T_\a)$ and a shadow $a^*$ for the same arc, and we isotope $a$ along the bigon, across $\S_\d$, and into the $\b$ side of the splitting.

We may not see immediately that the resulting splitting consists of two trivial tangles. However, if we carry out the process outlined above of finding bridge disks for both $f_\d(T_\a)$ and $f_\d(T_\b)$ traced out by $S$, we see that, since each bridge disk for $f_\d(T_\b)$ is parallel to a disk in $S$, we have no interior intersections between $a^*$ and $\b$ shadows. So $a$ is primitive with respect to $f_\d(T_\b)$, and so this isotopy is in fact a deperturbation. The deperturbed splitting consists of two trivial tangles, and we could recover $f_\d(T_\a) \cup \rho \cup f_\d(T_\b)$ by a perturbation. Thus, $U_1 \cup \rho \cup U_2$ is either in 2-bridge position with respect to $\S_\d$ or else perturbed, and in either case, the splitting meets condition (4) of Definition 4.1 from \cite{C&R}.

Since the collection $\{\rho_\a,\rho_\b,\rho_\g\}$ with the given framings is admissible in the sense of \cite{C&R}, and since $\mathcal{T}$ is the result of resolving the three tangles of $\mathcal{T}_\d$ along $\rho_\a$, $\rho_\b$, and $\rho_\g$ yields exactly $\mathcal{T}$ by our choice of $\rho = \rho_\a \cup \rho_\b$, Theorem 4.10 of \cite{C&R} tells us that $\mathcal{T}$ is the result of perturbing, connect summing with an $\R P^2$, or tubing on $\mathcal{T}_\d$, and thus the reverse direction is a detubing, the removal of an $\R P^2$ summand, or a deperturbation.

\end{proof}

In the 4-bridge case, all weakly-reducible trisections meet the conditions laid out above, and so weakly reducing such a splitting must realize an elementary operation. For non-trivial knotted spheres specifically, we can say further that the forgetting map realizes a detubing to a a trivial embedding of two spheres.

\begin{prop}
\label{prop:4Br}

    Let $\mathcal{T}$ be a weakly-reducible 4-bridge trisection of a 2-knot $K$ in $S^4$ with reducing arcs $\d \in T_\g^*$ and $\e$ parallel in each of $T_\a^*$ and $T_\b^*$. Then either the splitting is perturbed or the result of forgetting $\d$ detubes $K$ to a bridge trisection $\mathcal{T}'$ of an unlink of two $S^2$'s.
    
\end{prop}

\begin{proof}

From Lemma \ref{lemma:4BrNice} and Theorem \ref{thm:ElemOp}, we can conclude that forgetting $\d$ from $\mathcal{T}$ induces a deperturbation or detubing of $K$. We will suppose that $\mathcal{T}$ is not perturbed, and so forgetting $\d$ must be a detubing. This increases Euler characteristic by two and increases the number of components by at most one, so the result of the detubing must be a link $L$ of two spheres. Since $\chi(L) = -b + c_1 + c_2 + c_3$ (where $b$ is the bridge number and each $c_i$ the patch number of a pair of tangles) and $b = 3$, we must have at least one $c_i = 3$. By Theorem 5.1 of, \cite{jmmzsolids}, we see this must be an unlink.

\end{proof}

A sphere created by tubing an unlink of two spheres once is exactly a ribbon sphere of fusion number one, so we can conclude our main result after confirming that 4-bridge weakly-reducible trisections of non-spherical surfaces are trivial and showing that every such ribbon sphere arises in the desired form.

\begin{thm}
\label{thm:4Br}
A non-trivially knotted surface $K \subset S^4$ admits a 4-bridge weakly-reducible trisection if and only if it is a ribbon sphere with fusion number one.
\end{thm}

\begin{proof}

By Lemma \ref{lemma:4BrNice} and Theorem \ref{thm:ElemOp}, we see that forgetting $\d$ gives rise to a 3-bridge trisection for a surface link $L$, where $K$ is the result of perturbing, adding an unknotted $\R P^2$, or adding a tube to $L$. If $L$ is connected, then it is unknotted, (see Theorem 1.8 of \cite{mz}). Further, the addition of a single tube or $\R P^2$ summand to an unknotted surface in $S^4$ is trivial, as a result of Lemma 3 of \cite{B&S}, and so $K$ is unknotted whenever $L$ is connected.

In the case where $L$ is disconnected, we first notice that forgetting $\d$ must have detubed $K$, as tubing is the only of the three elementary operations that may change the number of components of a surface. Since we have detubed once, $L$ has exactly two components. Further, since $L$ is in 3-bridge position and has each patch number at least two, $L$ is an unlink of either two $S^2$'s or one $S^2$ and one $\R P^2$ (of either sign) (by Theorem 5.2 of \cite{jmmzsolids}). If $L$ contains an $\R P^2$, we examine the pair $f_\d(T_\a) \cup \overline{f_\d(T_\b)}$, which has exactly two components.

In \cite{MTZ}, Meier, Thompson, and Zupan give a procedure for realizing tubing on the level of shadows for a bridge trisection. If an arc $\eta$ has endpoints in the interiors of, say, an $\a$ shadow $a$ and a $\b$ shadow $b$ in different components of a shadow diagram, and $\eta$ is elsewhere disjoint from all $\a$ and $\b$ shadows, then we can tube along $\eta$, as in Figure \ref{fig:tube}, by replacing each of $a$ and $b$ with two shadows of the same color  with endpoints near $\eta$ connected by an arc dual to $\eta$. This will correspond to a tubing of the surface per Lemma 3.3 of \cite{MTZ}.

\begin{figure}[h]
    \centering

\tikzset{every picture/.style={line width=0.75pt}} 

\begin{tikzpicture}[x=0.75pt,y=0.75pt,yscale=-1,xscale=1]

\draw [color={rgb, 255:red, 126; green, 211; blue, 33 }  ,draw opacity=1 ]   (390,100) -- (390,140) ;
\draw [color={rgb, 255:red, 208; green, 2; blue, 27 }  ,draw opacity=1 ]   (330,80) .. controls (340.14,99) and (360.43,100.14) .. (390,100) ;
\draw [color={rgb, 255:red, 74; green, 144; blue, 226 }  ,draw opacity=1 ]   (450,80) .. controls (440.14,99.86) and (419.86,100.43) .. (390,100) ;
\draw [color={rgb, 255:red, 74; green, 144; blue, 226 }  ,draw opacity=1 ]   (450,160) .. controls (439.86,139.86) and (419.86,140.43) .. (390,140) ;
\draw [color={rgb, 255:red, 208; green, 2; blue, 27 }  ,draw opacity=1 ]   (330,160) .. controls (340.14,139.86) and (360.43,140.14) .. (390,140) ;
\draw [color={rgb, 255:red, 74; green, 144; blue, 226 }  ,draw opacity=1 ]   (90,80) -- (70,60) ;
\draw [color={rgb, 255:red, 74; green, 144; blue, 226 }  ,draw opacity=1 ]   (90,160) -- (70,180) ;
\draw [color={rgb, 255:red, 74; green, 144; blue, 226 }  ,draw opacity=1 ] [dash pattern={on 0.84pt off 2.51pt}]  (70,60) -- (50,40) ;
\draw [color={rgb, 255:red, 74; green, 144; blue, 226 }  ,draw opacity=1 ] [dash pattern={on 0.84pt off 2.51pt}]  (70,180) -- (50,200) ;
\draw [color={rgb, 255:red, 208; green, 2; blue, 27 }  ,draw opacity=1 ]   (90,80) .. controls (100.8,84.8) and (110.18,101.78) .. (110.32,119.46) .. controls (110.45,137.14) and (100,154.8) .. (90,160) ;
\draw [color={rgb, 255:red, 208; green, 2; blue, 27 }  ,draw opacity=1 ]   (210,80) -- (230,60) ;
\draw [color={rgb, 255:red, 208; green, 2; blue, 27 }  ,draw opacity=1 ]   (210,160) -- (230,180) ;
\draw [color={rgb, 255:red, 208; green, 2; blue, 27 }  ,draw opacity=1 ] [dash pattern={on 0.84pt off 2.51pt}]  (230,60) -- (250,40) ;
\draw [color={rgb, 255:red, 208; green, 2; blue, 27 }  ,draw opacity=1 ] [dash pattern={on 0.84pt off 2.51pt}]  (230,180) -- (250,200) ;
\draw [color={rgb, 255:red, 74; green, 144; blue, 226 }  ,draw opacity=1 ]   (210,80) .. controls (200.25,85.12) and (190.16,102.64) .. (190,120) .. controls (189.84,137.36) and (199.75,155.13) .. (210,160) ;
\draw  [fill={rgb, 255:red, 0; green, 0; blue, 0 }  ,fill opacity=1 ] (87,80) .. controls (87,78.34) and (88.34,77) .. (90,77) .. controls (91.66,77) and (93,78.34) .. (93,80) .. controls (93,81.66) and (91.66,83) .. (90,83) .. controls (88.34,83) and (87,81.66) .. (87,80) -- cycle ;
\draw  [fill={rgb, 255:red, 0; green, 0; blue, 0 }  ,fill opacity=1 ] (87,160) .. controls (87,158.34) and (88.34,157) .. (90,157) .. controls (91.66,157) and (93,158.34) .. (93,160) .. controls (93,161.66) and (91.66,163) .. (90,163) .. controls (88.34,163) and (87,161.66) .. (87,160) -- cycle ;
\draw  [fill={rgb, 255:red, 0; green, 0; blue, 0 }  ,fill opacity=1 ] (207,80) .. controls (207,78.34) and (208.34,77) .. (210,77) .. controls (211.66,77) and (213,78.34) .. (213,80) .. controls (213,81.66) and (211.66,83) .. (210,83) .. controls (208.34,83) and (207,81.66) .. (207,80) -- cycle ;
\draw  [fill={rgb, 255:red, 0; green, 0; blue, 0 }  ,fill opacity=1 ] (207,160) .. controls (207,158.34) and (208.34,157) .. (210,157) .. controls (211.66,157) and (213,158.34) .. (213,160) .. controls (213,161.66) and (211.66,163) .. (210,163) .. controls (208.34,163) and (207,161.66) .. (207,160) -- cycle ;
\draw [color={rgb, 255:red, 189; green, 16; blue, 224 }  ,draw opacity=1 ]   (110,120) -- (190,120) ;
\draw [color={rgb, 255:red, 74; green, 144; blue, 226 }  ,draw opacity=1 ]   (330,80) -- (310,60) ;
\draw [color={rgb, 255:red, 74; green, 144; blue, 226 }  ,draw opacity=1 ]   (330,160) -- (310,180) ;
\draw [color={rgb, 255:red, 74; green, 144; blue, 226 }  ,draw opacity=1 ] [dash pattern={on 0.84pt off 2.51pt}]  (310,60) -- (290,40) ;
\draw [color={rgb, 255:red, 74; green, 144; blue, 226 }  ,draw opacity=1 ] [dash pattern={on 0.84pt off 2.51pt}]  (310,180) -- (290,200) ;
\draw [color={rgb, 255:red, 208; green, 2; blue, 27 }  ,draw opacity=1 ]   (450,80) -- (470,60) ;
\draw [color={rgb, 255:red, 208; green, 2; blue, 27 }  ,draw opacity=1 ]   (450,160) -- (470,180) ;
\draw [color={rgb, 255:red, 208; green, 2; blue, 27 }  ,draw opacity=1 ] [dash pattern={on 0.84pt off 2.51pt}]  (470,60) -- (490,40) ;
\draw [color={rgb, 255:red, 208; green, 2; blue, 27 }  ,draw opacity=1 ] [dash pattern={on 0.84pt off 2.51pt}]  (470,180) -- (490,200) ;
\draw  [fill={rgb, 255:red, 0; green, 0; blue, 0 }  ,fill opacity=1 ] (327,80) .. controls (327,78.34) and (328.34,77) .. (330,77) .. controls (331.66,77) and (333,78.34) .. (333,80) .. controls (333,81.66) and (331.66,83) .. (330,83) .. controls (328.34,83) and (327,81.66) .. (327,80) -- cycle ;
\draw  [fill={rgb, 255:red, 0; green, 0; blue, 0 }  ,fill opacity=1 ] (327,160) .. controls (327,158.34) and (328.34,157) .. (330,157) .. controls (331.66,157) and (333,158.34) .. (333,160) .. controls (333,161.66) and (331.66,163) .. (330,163) .. controls (328.34,163) and (327,161.66) .. (327,160) -- cycle ;
\draw  [fill={rgb, 255:red, 0; green, 0; blue, 0 }  ,fill opacity=1 ] (447,80) .. controls (447,78.34) and (448.34,77) .. (450,77) .. controls (451.66,77) and (453,78.34) .. (453,80) .. controls (453,81.66) and (451.66,83) .. (450,83) .. controls (448.34,83) and (447,81.66) .. (447,80) -- cycle ;
\draw  [fill={rgb, 255:red, 0; green, 0; blue, 0 }  ,fill opacity=1 ] (447,160) .. controls (447,158.34) and (448.34,157) .. (450,157) .. controls (451.66,157) and (453,158.34) .. (453,160) .. controls (453,161.66) and (451.66,163) .. (450,163) .. controls (448.34,163) and (447,161.66) .. (447,160) -- cycle ;
\draw  [fill={rgb, 255:red, 0; green, 0; blue, 0 }  ,fill opacity=1 ] (387,100) .. controls (387,98.34) and (388.34,97) .. (390,97) .. controls (391.66,97) and (393,98.34) .. (393,100) .. controls (393,101.66) and (391.66,103) .. (390,103) .. controls (388.34,103) and (387,101.66) .. (387,100) -- cycle ;
\draw  [fill={rgb, 255:red, 0; green, 0; blue, 0 }  ,fill opacity=1 ] (387,140) .. controls (387,138.34) and (388.34,137) .. (390,137) .. controls (391.66,137) and (393,138.34) .. (393,140) .. controls (393,141.66) and (391.66,143) .. (390,143) .. controls (388.34,143) and (387,141.66) .. (387,140) -- cycle ;
\draw    (230,120) -- (308,120) ;
\draw [shift={(310,120)}, rotate = 180] [color={rgb, 255:red, 0; green, 0; blue, 0 }  ][line width=0.75]    (10.93,-3.29) .. controls (6.95,-1.4) and (3.31,-0.3) .. (0,0) .. controls (3.31,0.3) and (6.95,1.4) .. (10.93,3.29)   ;

\draw (144.33,123.4) node [anchor=north west][inner sep=0.75pt]  [color={rgb, 255:red, 189; green, 16; blue, 224 }  ,opacity=1 ]  {$\eta $};
\draw (91,112.4) node [anchor=north west][inner sep=0.75pt]    {$\textcolor[rgb]{0.82,0.01,0.11}{a}$};
\draw (201,112.4) node [anchor=north west][inner sep=0.75pt]    {$\textcolor[rgb]{0.29,0.56,0.89}{b}$};

\end{tikzpicture}
    \caption{A tubing on the level of shadows}
    \label{fig:tube}
\end{figure}
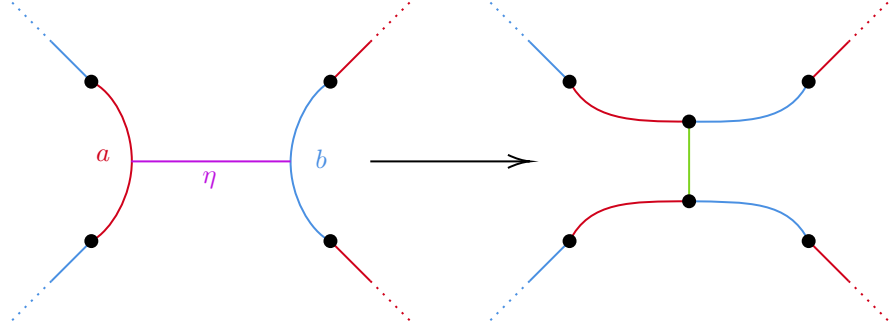

From Lemma 2.4 of \cite{goeritz}, we know that the genus two splitting surface of $S^2 \x S^1$ has a unique non-separating reducing curve $c$. Since it is a curve on a genus two surface, $c$ is equivariant under the hyper-elliptic involution, and so when we view $S^2 \x S^1$ as the double cover of $S^3 = f_\d(H_\a) \cup \overline{f_\d(H_\b)}$ branched over $f_\d(T_\a) \cup \overline{f_\d(T_\b)}$, $c$ covers a unique curve in $\S_\d$ corresponding to the 2-bridge component of $f_\d(T_\a) \cup \overline{f_\d(T_\b)}$. Since we know that $\e_\a$ and $\e_\b$ share a shadow $\e^*$, a shadow diagram for $f_\d(T_\a) \cup \overline{f_\d(T_\b)}$ is isotopic to that in Figure \ref{fig:RP2andS2}.

Notice that in this case, there is a unique curve separating the two components of the shadow diagram. As such, any arc connecting the interior of an $\a$ shadow from one component to the interior of a $\b$ shadow from the other may be isotoped to intersect this splitting curve once. Thus, tubing along any such arc is equivalent to a trivial tubing; i.e. a connect sum.
 
If forgetting $\d$ gives rise to a link $L$ of two spheres, then since $K$ is obtained by tubing these two spheres together, $K$ must be a ribbon sphere.

\begin{figure}[h]
    \centering

\tikzset{every picture/.style={line width=0.75pt}} 

\begin{tikzpicture}[x=0.75pt,y=0.75pt,yscale=-1,xscale=1]

\draw [color={rgb, 255:red, 208; green, 2; blue, 27 }  ,draw opacity=1 ]   (269.88,80) .. controls (243.18,98.6) and (243.58,142.6) .. (269.88,160) ;
\draw [color={rgb, 255:red, 74; green, 144; blue, 226 }  ,draw opacity=1 ]   (269.88,80) .. controls (296.78,98) and (296.18,142.2) .. (269.88,160) ;
\draw [color={rgb, 255:red, 74; green, 144; blue, 226 }  ,draw opacity=1 ]   (170,80) -- (170,160) ;
\draw [color={rgb, 255:red, 74; green, 144; blue, 226 }  ,draw opacity=1 ]   (90,80) -- (90,160) ;
\draw [color={rgb, 255:red, 208; green, 2; blue, 27 }  ,draw opacity=1 ]   (90,160) -- (170,160) ;
\draw [color={rgb, 255:red, 208; green, 2; blue, 27 }  ,draw opacity=1 ]   (90,80) -- (170,80) ;
\draw  [fill={rgb, 255:red, 0; green, 0; blue, 0 }  ,fill opacity=1 ] (88,80) .. controls (88,78.9) and (88.9,78) .. (90,78) .. controls (91.1,78) and (92,78.9) .. (92,80) .. controls (92,81.1) and (91.1,82) .. (90,82) .. controls (88.9,82) and (88,81.1) .. (88,80) -- cycle ;
\draw  [fill={rgb, 255:red, 0; green, 0; blue, 0 }  ,fill opacity=1 ] (168,80) .. controls (168,78.9) and (168.9,78) .. (170,78) .. controls (171.1,78) and (172,78.9) .. (172,80) .. controls (172,81.1) and (171.1,82) .. (170,82) .. controls (168.9,82) and (168,81.1) .. (168,80) -- cycle ;
\draw  [fill={rgb, 255:red, 0; green, 0; blue, 0 }  ,fill opacity=1 ] (88,160) .. controls (88,158.9) and (88.9,158) .. (90,158) .. controls (91.1,158) and (92,158.9) .. (92,160) .. controls (92,161.1) and (91.1,162) .. (90,162) .. controls (88.9,162) and (88,161.1) .. (88,160) -- cycle ;
\draw  [fill={rgb, 255:red, 0; green, 0; blue, 0 }  ,fill opacity=1 ] (168,160) .. controls (168,158.9) and (168.9,158) .. (170,158) .. controls (171.1,158) and (172,158.9) .. (172,160) .. controls (172,161.1) and (171.1,162) .. (170,162) .. controls (168.9,162) and (168,161.1) .. (168,160) -- cycle ;
\draw  [fill={rgb, 255:red, 0; green, 0; blue, 0 }  ,fill opacity=1 ] (267.88,80) .. controls (267.88,78.9) and (268.77,78) .. (269.88,78) .. controls (270.98,78) and (271.88,78.9) .. (271.88,80) .. controls (271.88,81.1) and (270.98,82) .. (269.88,82) .. controls (268.77,82) and (267.88,81.1) .. (267.88,80) -- cycle ;
\draw  [fill={rgb, 255:red, 0; green, 0; blue, 0 }  ,fill opacity=1 ] (267.88,160) .. controls (267.88,158.9) and (268.77,158) .. (269.88,158) .. controls (270.98,158) and (271.88,158.9) .. (271.88,160) .. controls (271.88,161.1) and (270.98,162) .. (269.88,162) .. controls (268.77,162) and (267.88,161.1) .. (267.88,160) -- cycle ;
\draw [color={rgb, 255:red, 126; green, 211; blue, 33 }  ,draw opacity=1 ]   (170,120) .. controls (190.17,119.67) and (200.2,50.2) .. (270,50) .. controls (339.8,49.8) and (350.17,80.33) .. (350,120) .. controls (349.83,159.67) and (340.2,189.8) .. (270,190) .. controls (199.8,190.2) and (199.4,119.8) .. (250,120) ;
\draw   (230,120) .. controls (230,86.86) and (247.91,60) .. (270,60) .. controls (292.09,60) and (310,86.86) .. (310,120) .. controls (310,153.14) and (292.09,180) .. (270,180) .. controls (247.91,180) and (230,153.14) .. (230,120) -- cycle ;

\end{tikzpicture}
    \caption{$\a$ and $\b$ shadows for an unlinked $\R P^2$ and $S^2$ with separating curve and a trivial tubing arc}
    \label{fig:RP2andS2}
\end{figure}
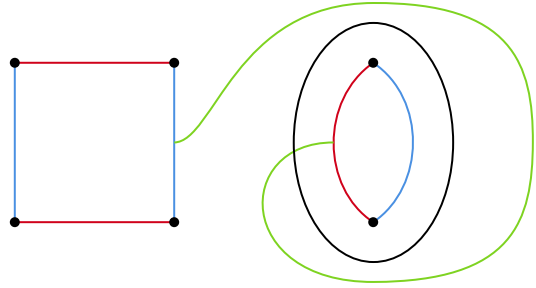

Consider now any ribbon sphere $R$ of fusion number one. Such a sphere is obtained by adding a single tube, or boundary of a 3-dimensional 1-handle, to an unlink of two $S^2$'s. The isotopy class of this tube is determined by its winding around each of the spheres, as described in \cite{Marumoto}, and encoded by a list $p_1,q_1,p_2,q_2,...,p_m,q_m$ of integers. In fact, we can easily manipulate the diagram in Figure 4 of \cite{Marumoto} to obtain a bridge trisection of any such surface.

\begin{figure}[h]
    \centering

\tikzset{every picture/.style={line width=0.75pt}} 

\begin{tikzpicture}[x=0.75pt,y=0.75pt,yscale=-1,xscale=1]

\draw    (95,50) .. controls (94.67,56.33) and (160.34,64.67) .. (169.98,65) .. controls (179.62,65.32) and (182.29,67.89) .. (183.85,68.91) ;
\draw    (95,70) .. controls (94.8,63) and (160,74.8) .. (170,75) .. controls (180,75.2) and (183.9,83.6) .. (172.5,85) ;
\draw    (95,245) .. controls (95.1,274.99) and (72.42,275) .. (60,275) .. controls (47.58,275) and (25.1,274.99) .. (25,245) .. controls (24.9,215.01) and (24.79,79.99) .. (25,50) .. controls (25.21,20.01) and (47.37,19.83) .. (60,20) .. controls (72.63,20.17) and (95.1,20.26) .. (95,50) ;
\draw    (165.18,49.98) .. controls (165.09,19.98) and (187.58,19.99) .. (200,20) .. controls (212.42,20.01) and (234.92,20.01) .. (235,50) .. controls (235.08,79.99) and (235.23,215.01) .. (235,245) .. controls (234.77,274.99) and (212.63,275.18) .. (200,275) .. controls (187.37,274.82) and (164.87,274.74) .. (165,245) ;
\draw    (130,60) .. controls (132,62.09) and (132,68.11) .. (130,70) ;
\draw    [dash pattern={on 1pt off 2pt}] (130,60) .. controls (127.81,63.38) and (128.13,67.56) .. (130,70) ;
\draw    (183.85,68.91) .. controls (188.05,71.69) and (201.2,88.35) .. (172.5,94.99) ;
\draw    (165.18,49.98) -- (165.26,60.47) ;
\draw    (165,80) -- (165,100) ;
\draw   (175,89.99) .. controls (175,87.24) and (173.88,85) .. (172.5,85) .. controls (171.12,85) and (170,87.24) .. (170,89.99) .. controls (170,92.75) and (171.12,94.99) .. (172.5,94.99) .. controls (173.88,94.99) and (175,92.75) .. (175,89.99) -- cycle ;
\draw   (160.3,95) .. controls (160.3,92.24) and (159.11,90) .. (157.65,90) .. controls (156.19,90) and (155,92.24) .. (155,95) .. controls (155,97.76) and (156.19,100) .. (157.65,100) .. controls (159.11,100) and (160.3,97.76) .. (160.3,95) -- cycle ;
\draw    (157.65,100) .. controls (153.79,102) and (161.07,103.57) .. (165,105) .. controls (168.93,106.43) and (203.7,123.21) .. (172.5,135) ;
\draw    (172.5,125.01) .. controls (181.7,122.61) and (169.07,116.86) .. (165,115) .. controls (160.93,113.14) and (125.7,100) .. (157.65,90) ;
\draw   (175,130.01) .. controls (175,127.25) and (173.88,125.01) .. (172.5,125.01) .. controls (171.12,125.01) and (170,127.25) .. (170,130.01) .. controls (170,132.76) and (171.12,135) .. (172.5,135) .. controls (173.88,135) and (175,132.76) .. (175,130.01) -- cycle ;
\draw    (165,120) -- (165,145) ;
\draw   (160,135) .. controls (160,132.24) and (158.81,130) .. (157.35,130) .. controls (155.89,130) and (154.7,132.24) .. (154.7,135) .. controls (154.7,137.76) and (155.89,140) .. (157.35,140) .. controls (158.81,140) and (160,137.76) .. (160,135) -- cycle ;
\draw    (87.65,150) .. controls (80.36,145.71) and (88.36,145.29) .. (95,145) .. controls (101.64,144.71) and (140.81,142.61) .. (157.35,140) ;
\draw    (87.65,160) .. controls (74.01,159.57) and (65.79,136.71) .. (90,135) .. controls (114.21,133.29) and (145.79,132.23) .. (157.35,130) ;
\draw    (95,70) -- (95,130) ;
\draw   (90.3,155) .. controls (90.3,152.24) and (89.11,150) .. (87.65,150) .. controls (86.19,150) and (85,152.24) .. (85,155) .. controls (85,157.76) and (86.19,160) .. (87.65,160) .. controls (89.11,160) and (90.3,157.76) .. (90.3,155) -- cycle ;
\draw    (95,150) -- (95,245) ;
\draw    (165,145) -- (165,210) ;
\draw   (105.3,160) .. controls (105.3,157.24) and (104.11,155) .. (102.65,155) .. controls (101.19,155) and (100,157.24) .. (100,160) .. controls (100,162.76) and (101.19,165) .. (102.65,165) .. controls (104.11,165) and (105.3,162.76) .. (105.3,160) -- cycle ;
\draw    (102.65,155) .. controls (108.21,157.14) and (118.79,159.29) .. (125,160) ;
\draw    (102.65,165) .. controls (108.21,167.14) and (118.79,169.29) .. (125,170) ;
\draw  [dash pattern={on 0.84pt off 2.51pt}]  (125,165) .. controls (139.29,167.29) and (175.21,173.71) .. (175,180) .. controls (174.79,186.29) and (85.21,195.71) .. (85,205) .. controls (84.79,214.29) and (122.43,215.29) .. (135,220) ;
\draw    (135,215) .. controls (144.3,218.5) and (165.1,225.9) .. (165,210) ;
\draw    (135,225) .. controls (144.3,228.5) and (165.08,231.38) .. (165,245) ;
\draw    (151.53,220) .. controls (153.53,222.09) and (153.53,228.11) .. (151.53,230) ;
\draw   [dash pattern={on 1pt off 2pt}] (151.53,220) .. controls (149.34,223.37) and (149.65,227.56) .. (151.53,230) ;
\draw    (25,115) .. controls (40.29,130.14) and (80,130.14) .. (95,115) ;
\draw  [dash pattern={on 0.84pt off 2.51pt}]  (25,115) .. controls (40,99.86) and (80.29,100.14) .. (95,115) ;
\draw    (165,153.64) .. controls (180.29,168.79) and (220,168.79) .. (235,153.64) ;
\draw  [dash pattern={on 0.84pt off 2.51pt}]  (165,153.64) .. controls (180,138.5) and (220.29,138.79) .. (235,153.64) ;
\draw    (364.88,50.02) .. controls (364.54,56.35) and (430.21,64.69) .. (439.85,65.02) .. controls (449.49,65.34) and (452.17,67.91) .. (453.73,68.93) ;
\draw    (364.88,70.02) .. controls (364.68,63.02) and (429.88,74.82) .. (439.88,75.02) .. controls (449.88,75.22) and (451.4,83.6) .. (440,85) ;
\draw    (364.88,245.02) .. controls (364.98,275.01) and (342.29,275.02) .. (329.88,275.02) .. controls (317.46,275.02) and (294.98,275.01) .. (294.88,245.02) .. controls (294.77,215.03) and (294.67,80.01) .. (294.88,50.02) .. controls (295.09,20.03) and (317.24,19.85) .. (329.88,20.02) .. controls (342.51,20.19) and (364.98,20.28) .. (364.88,50.02) ;
\draw    (435.05,50) .. controls (434.97,20) and (457.46,20.01) .. (469.88,20.02) .. controls (482.29,20.03) and (504.79,20.03) .. (504.88,50.02) .. controls (504.96,80.01) and (505.11,215.03) .. (504.88,245.02) .. controls (504.65,275.01) and (482.51,275.2) .. (469.88,275.02) .. controls (457.24,274.84) and (434.75,274.76) .. (434.88,245.02) ;
\draw    (453.73,68.93) .. controls (457.93,71.71) and (468.7,88.36) .. (440,95) ;
\draw    (435.05,50) -- (435.14,60.49) ;
\draw    (434.88,80.02) -- (434.88,100.02) ;
\draw    (430,100) .. controls (426.14,102) and (430.95,103.59) .. (434.88,105.02) .. controls (438.81,106.45) and (471.2,123.21) .. (440,135) ;
\draw    (440,125) .. controls (449.2,122.6) and (438.95,116.88) .. (434.88,115.02) .. controls (430.81,113.16) and (398.05,100) .. (430,90) ;
\draw    (434.88,120.02) -- (434.88,145.02) ;
\draw    (360,150) .. controls (352.71,145.71) and (358.36,145.29) .. (365,145) .. controls (371.64,144.71) and (413.46,142.61) .. (430,140) ;
\draw    (360,160) .. controls (340.33,154.83) and (335.79,136.71) .. (360,135) .. controls (384.21,133.29) and (418.44,132.23) .. (430,130) ;
\draw    (364.88,70.02) -- (364.88,130.02) ;
\draw    (364.88,150.02) -- (364.88,245.02) ;
\draw    (434.88,145.02) -- (434.88,210.02) ;
\draw    (370,155) .. controls (375.56,157.14) and (388.66,159.31) .. (394.88,160.02) ;
\draw    (370,165) .. controls (375.56,167.14) and (388.66,169.31) .. (394.88,170.02) ;
\draw  [dash pattern={on 0.84pt off 2.51pt}]  (394.88,165.02) .. controls (409.16,167.31) and (445.09,173.73) .. (444.88,180.02) .. controls (444.66,186.31) and (355.09,195.73) .. (354.88,205.02) .. controls (354.66,214.31) and (392.31,215.31) .. (404.88,220.02) ;
\draw    (404.88,215.02) .. controls (414.18,218.52) and (434.98,225.92) .. (434.88,210.02) ;
\draw    (404.88,225.02) .. controls (414.18,228.52) and (434.95,231.4) .. (434.88,245.02) ;

\draw (190,95) node [anchor=north west][inner sep=0.75pt]    {$p_{1}$};
\draw (55,140) node [anchor=north west][inner sep=0.75pt]    {$q_{1}$};

\end{tikzpicture}
    \caption{A ribbon 2-knot and an equatorial ribbon 1-knot}
    \label{fig:1fusion}
\end{figure}
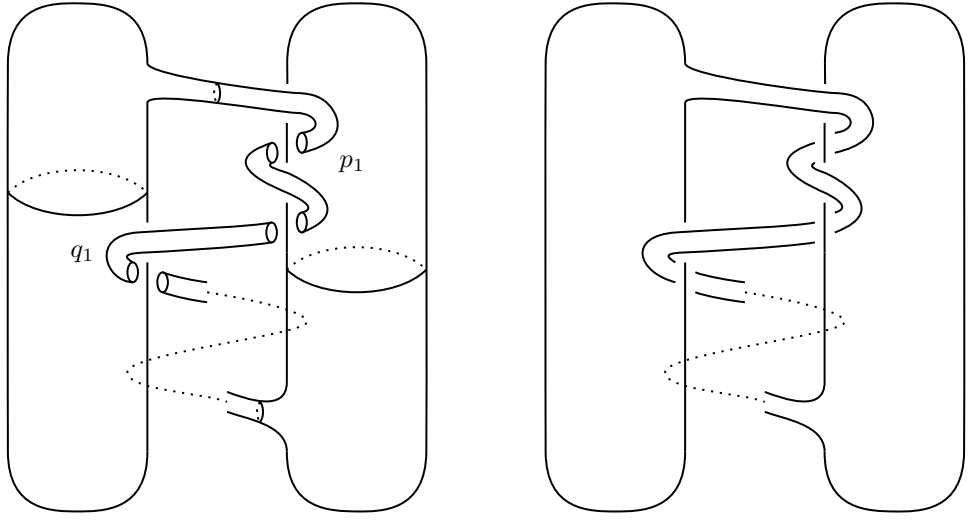

As seen in Figure \ref{fig:1fusion}, we take the given ball-tube presentation of $R$ and find an equatorial 1-knot $R'$ appearing as a disk-band presentation in 3-bridge position (with respect to some splitting surface), consisting of an unlink of two $S^1$'s and a band connecting them. We isotope the splitting surface such that one of the tangles in the splitting contains no crossings. We now resolve along a band $\mathfrak{b}$ 0-framed by the splitting surface and with endpoints on opposite sides of the ribbon band near one of the circles. We record a parallel copy of $\mathfrak{b}$ to one side of the resolution and a band $\mathfrak{b}'$ dual to $\mathfrak{b}$, whose resolution recovers $R'$. We then perturb between these two recorded bands. The resulting link (ignoring for a moment $\mathfrak{b}$ and $\mathfrak{b}'$) is a 2-component unlink in 4-bridge position, viewed as the union of two trivial tangles $T_1$ and $T_2$, where $T_2$ contains no crossings.

After pushing both bands off the splitting surface into $T_2$, we have three trivial tangles, $T_1$, $T_2$, and $T_3 := T_{2\mathfrak{b}\mathfrak{b}'}$, the result of resolving $T_2$ along both bands. Since $R'_\mathfrak{b}$ is an unlink, each of $T_1 \cup T_2$ and $T_1 \cup T_3$ is an unlink (also with two components). Since $T_2 \cup T_3$ contains no crossings, it is likewise unlinked, and thus the three tangles form a bridge trisection. Refer to Section 3 of \cite{mz} for a detailed discussion of this procedure for generating bridge trisections.

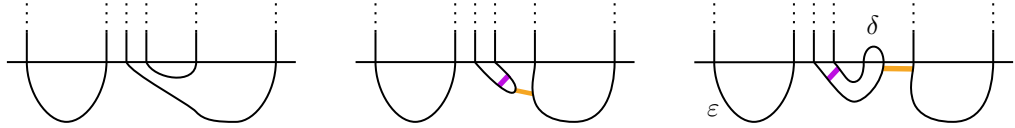
\begin{figure}[h]
    \centering

\tikzset{every picture/.style={line width=0.75pt}} 

\begin{tikzpicture}[x=0.75pt,y=0.75pt,yscale=-1,xscale=1]

\draw  [color={rgb, 255:red, 245; green, 166; blue, 35 }  ,draw opacity=1 ][fill={rgb, 255:red, 245; green, 166; blue, 35 }  ,fill opacity=1 ] (285.41,68.53) -- (293.82,70.26) -- (293.56,71.54) -- (285.15,69.81) -- cycle ;
\draw  [color={rgb, 255:red, 189; green, 16; blue, 224 }  ,draw opacity=1 ][fill={rgb, 255:red, 189; green, 16; blue, 224 }  ,fill opacity=1 ] (276.47,66.18) -- (281.11,61.21) -- (282.34,62.36) -- (277.71,67.34) -- cycle ;
\draw    (40,55) .. controls (40,69.92) and (50,85.09) .. (60,85) .. controls (70,84.91) and (79.83,70.25) .. (80,55) ;
\draw    (30,55) -- (175,55) ;
\draw    (90,55) .. controls (92.5,60.93) and (121.22,75.78) .. (125,80) .. controls (128.78,84.22) and (135,85.09) .. (145,85) .. controls (155,84.91) and (164.83,70.25) .. (165,55) ;
\draw    (100,55) .. controls (103.67,63.56) and (125,67.11) .. (125,55) ;
\draw    (215,55) .. controls (215,69.92) and (225,85.09) .. (235,85) .. controls (245,84.91) and (254.83,70.25) .. (255,55) ;
\draw    (205,55) -- (345,55) ;
\draw    (295,55) .. controls (295.22,60.22) and (292.11,67) .. (295,75) .. controls (295.89,77.78) and (300.11,85.11) .. (315,85) .. controls (329.89,84.89) and (334.83,70.25) .. (335,55) ;
\draw    (274.79,54.91) .. controls (281.4,60.83) and (287,66.6) .. (284.9,69.75) .. controls (282.79,72.91) and (272.8,64.42) .. (264.79,54.91) ;
\draw  [color={rgb, 255:red, 245; green, 166; blue, 35 }  ,draw opacity=1 ][fill={rgb, 255:red, 245; green, 166; blue, 35 }  ,fill opacity=1 ] (470,57.35) -- (484.75,57.56) -- (484.72,59.29) -- (469.98,59.08) -- cycle ;
\draw  [color={rgb, 255:red, 189; green, 16; blue, 224 }  ,draw opacity=1 ][fill={rgb, 255:red, 189; green, 16; blue, 224 }  ,fill opacity=1 ] (442,62.61) -- (446.64,57.64) -- (447.88,58.79) -- (443.24,63.76) -- cycle ;
\draw    (385,55.09) .. controls (385,70) and (395,85.09) .. (405,85) .. controls (415,84.91) and (424.83,70.34) .. (425,55.09) ;
\draw    (375,55.09) -- (535,55) ;
\draw    (485,55) .. controls (485.22,60.22) and (482.11,67) .. (485,75) .. controls (485.89,77.78) and (490.11,85.11) .. (505,85) .. controls (519.89,84.89) and (524.83,70.25) .. (525,55) ;
\draw    (470,55) .. controls (470,65) and (460.36,75.36) .. (454.9,74.84) .. controls (449.43,74.32) and (443,64.51) .. (435,55) ;
\draw    (445,55) .. controls (450,60.45) and (451.45,65.18) .. (455,65) .. controls (458.55,64.82) and (460.18,60.27) .. (460,55) ;
\draw    (460,55) .. controls (460.18,44.64) and (469.82,44.64) .. (470,55) ;
\draw    (40,55) -- (40,40) ;
\draw    (80,55) -- (80,40) ;
\draw    (90,55) -- (90,40) ;
\draw    (100,55) -- (100,40) ;
\draw    (125,55) -- (125,40) ;
\draw    (165,55) -- (165,40) ;
\draw  [dash pattern={on 0.84pt off 2.51pt}]  (40,40) -- (40,25) ;
\draw  [dash pattern={on 0.84pt off 2.51pt}]  (80,40) -- (80,25) ;
\draw  [dash pattern={on 0.84pt off 2.51pt}]  (90,40) -- (90,25) ;
\draw  [dash pattern={on 0.84pt off 2.51pt}]  (100,40) -- (100,25) ;
\draw  [dash pattern={on 0.84pt off 2.51pt}]  (125,40) -- (125,25) ;
\draw  [dash pattern={on 0.84pt off 2.51pt}]  (165,40) -- (165,25) ;
\draw    (215,55) -- (215,40) ;
\draw    (255,55) -- (255,40) ;
\draw    (265,55) -- (265,40) ;
\draw    (275,55) -- (275,40) ;
\draw    (295,55) -- (295,40) ;
\draw    (335,55) -- (335,40) ;
\draw  [dash pattern={on 0.84pt off 2.51pt}]  (215,40) -- (215,25) ;
\draw  [dash pattern={on 0.84pt off 2.51pt}]  (255,40) -- (255,25) ;
\draw  [dash pattern={on 0.84pt off 2.51pt}]  (265,40) -- (265,25) ;
\draw  [dash pattern={on 0.84pt off 2.51pt}]  (275,40) -- (275,25) ;
\draw  [dash pattern={on 0.84pt off 2.51pt}]  (295,40) -- (295,25) ;
\draw  [dash pattern={on 0.84pt off 2.51pt}]  (335,40) -- (335,25) ;
\draw    (385,55) -- (385,40) ;
\draw    (425,55) -- (425,40) ;
\draw    (435,55) -- (435,40) ;
\draw    (445,55) -- (445,40) ;
\draw    (485,55) -- (485,40) ;
\draw    (525,55) -- (525,40) ;
\draw  [dash pattern={on 0.84pt off 2.51pt}]  (385,40) -- (385,25) ;
\draw  [dash pattern={on 0.84pt off 2.51pt}]  (425,40) -- (425,25) ;
\draw  [dash pattern={on 0.84pt off 2.51pt}]  (435,40) -- (435,25) ;
\draw  [dash pattern={on 0.84pt off 2.51pt}]  (445,40) -- (445,25) ;
\draw  [dash pattern={on 0.84pt off 2.51pt}]  (485,40) -- (485,25) ;
\draw  [dash pattern={on 0.84pt off 2.51pt}]  (525,40) -- (525,25) ;

\draw (380,75) node [anchor=north west][inner sep=0.75pt]    {$\varepsilon $};
\draw (460,30) node [anchor=north west][inner sep=0.75pt]    {$\delta $};

\end{tikzpicture}
    \caption{Producing a bridge trisection from a 1-fusion ribbon 1-knot}
    \label{fig:BandPert}
\end{figure}

As seen in Figure \ref{fig:BandPert}, the perturbation near the bands, which is entirely local and does not depend on the specific 2-knot being trisected, gives rise to a pair of arcs $\d \in T_1$ and $\e$ parallel in each of $T_2$ and $T_3$, so the bridge trisection we obtain is weakly-reducible.

\end{proof}

Finally, we answer a natural next question: ``Are there any strongly-irreducible 4-bridge trisections?" in the following Corollary.

\begin{cor}
    \label{cor:2twist}
    The 2-twist spin of the trefoil admits a strongly-irreducible 4-bridge trisection.
\end{cor}

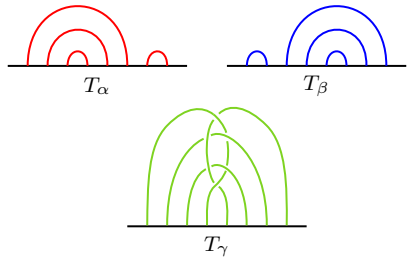
\begin{figure}[h!]
    \centering

\tikzset{every picture/.style={line width=0.75pt}} 

\begin{tikzpicture}[x=0.75pt,y=0.75pt,yscale=-1,xscale=1]

\draw    (20,80) -- (110,80) ;
\draw [color={rgb, 255:red, 255; green, 0; blue, 0 }  ,draw opacity=1 ]   (90,80) .. controls (90.09,73.22) and (94.95,71.15) .. (97.86,73.8) .. controls (99.14,74.96) and (100.04,77.03) .. (100,80) ;
\draw [color={rgb, 255:red, 255; green, 0; blue, 0 }  ,draw opacity=1 ]   (50,80) .. controls (50.13,70.25) and (60.13,70.25) .. (60,80) ;
\draw [color={rgb, 255:red, 255; green, 0; blue, 0 }  ,draw opacity=1 ]   (30,80) .. controls (30.13,40) and (80.13,40) .. (80,80) ;
\draw [color={rgb, 255:red, 255; green, 0; blue, 0 }  ,draw opacity=1 ]   (40,80) .. controls (40,56) and (70.5,55.5) .. (70,80) ;
\draw    (130,80) -- (220,80) ;
\draw [color={rgb, 255:red, 0; green, 0; blue, 255 }  ,draw opacity=1 ]   (140,80) .. controls (140.13,70.25) and (150.13,70.25) .. (150,80) ;
\draw [color={rgb, 255:red, 0; green, 0; blue, 255 }  ,draw opacity=1 ]   (179.99,80) .. controls (180.12,70.25) and (190.12,70.25) .. (189.99,80) ;
\draw [color={rgb, 255:red, 0; green, 0; blue, 255 }  ,draw opacity=1 ]   (170,80) .. controls (170,56) and (200.5,55.5) .. (200,80) ;
\draw    (80,160.5) -- (170,160.5) ;
\draw [color={rgb, 255:red, 0; green, 0; blue, 255 }  ,draw opacity=1 ]   (160,80) .. controls (160.13,40) and (210.13,40) .. (210,80) ;
\draw [color={rgb, 255:red, 126; green, 211; blue, 33 }  ,draw opacity=1 ]   (90,160) .. controls (90.22,140.33) and (90.11,120.05) .. (100,109.95) .. controls (109.89,99.84) and (118,100.37) .. (123.16,104.3) ;
\draw [color={rgb, 255:red, 126; green, 211; blue, 33 }  ,draw opacity=1 ]   (160,160) .. controls (160.22,140.11) and (159.89,119.95) .. (150,110) .. controls (140.11,100.05) and (131.22,98.89) .. (126,104) ;
\draw [color={rgb, 255:red, 126; green, 211; blue, 33 }  ,draw opacity=1 ]   (123.54,138.63) .. controls (120.79,135.13) and (119.88,125.21) .. (119.88,122.88) .. controls (119.88,120.54) and (119.79,110.96) .. (123.38,107.13) ;
\draw [color={rgb, 255:red, 126; green, 211; blue, 33 }  ,draw opacity=1 ]   (130,130) .. controls (130.64,126.57) and (130.91,121.51) .. (130.13,117.79) ;
\draw [color={rgb, 255:red, 126; green, 211; blue, 33 }  ,draw opacity=1 ]   (120,160) .. controls (120.15,155) and (119.04,144.79) .. (123.71,140.54) .. controls (128.38,136.29) and (129.13,133.46) .. (129.46,132.63) ;
\draw [color={rgb, 255:red, 126; green, 211; blue, 33 }  ,draw opacity=1 ]   (130,159.96) .. controls (130.15,154.96) and (131.29,145.07) .. (125.93,140.57) ;
\draw [color={rgb, 255:red, 126; green, 211; blue, 33 }  ,draw opacity=1 ]   (122,114.89) .. controls (129.68,110.47) and (150.21,125.31) .. (150,160.11) ;
\draw [color={rgb, 255:red, 126; green, 211; blue, 33 }  ,draw opacity=1 ]   (122.27,130.46) .. controls (127.16,127.11) and (140.1,134.88) .. (139.89,160) ;
\draw [color={rgb, 255:red, 126; green, 211; blue, 33 }  ,draw opacity=1 ]   (119.47,131.84) .. controls (119.16,132.05) and (110.21,135) .. (110,160) ;
\draw [color={rgb, 255:red, 126; green, 211; blue, 33 }  ,draw opacity=1 ]   (118.32,116.58) .. controls (116.84,117.11) and (100,125) .. (100,160) ;
\draw [color={rgb, 255:red, 126; green, 211; blue, 33 }  ,draw opacity=1 ]   (129.68,113.53) .. controls (128.74,110.05) and (125.37,105.84) .. (123.16,104.3) ;

\draw (65,90) node  [font=\footnotesize]  {$T_{\alpha }$};
\draw (175,90) node  [font=\footnotesize]  {$T_{\beta }$};
\draw (125,170.5) node  [font=\footnotesize]  {$T_{\gamma }$};

\end{tikzpicture}
    \caption{A diagram for a 4-bridge trisection of the 2-twist spun trefoil}
    \label{fig:2twist}
\end{figure}

\begin{proof}
    A 4-bridge trisection for the 2-twist spun trefoil is known (see Figure \ref{fig:2twist}). Since we know from \cite{Cochran} that the 2-twist spin of the trefoil is non-ribbon, Theorem \ref{thm:4Br} tells us that this bridge trisection is strongly irreducible.
\end{proof}

\printbibliography

\end{document}